\documentclass[10pt]{amsart}

\usepackage{enumerate}
\usepackage{amssymb, amsmath}
\usepackage{color}

\newtheorem{theorem}{Theorem}
\newtheorem{proposition}[theorem]{Proposition}
\newtheorem{lemma}[theorem]{Lemma}
\newtheorem{corollary}[theorem]{Corollary}
\newtheorem{question}{Question}
\newtheorem{claim}[theorem]{Claim}

\newtheorem{definition}{Definition}

\newtheorem{remark}{Remark}

\newcommand \bN{\Bbb{N}}
\newcommand \bZ{\Bbb{Z}}

\newcommand \eps{\varepsilon}

\newcommand \cM{\mathcal{M}}

\newcommand \restrict{\upharpoonright}

\newcommand{\JX}{EC}
\newcommand{\JXn}{ECn}

\newcommand{\Ymn}{\mathcal{Y}^-_n}

\newcommand{\Ymz}{\mathcal{Y}^-_0}

\newcommand{\Ymnm}{\mathcal{Y}^-_{n-1}}
\newcommand{\bbk}{\mathbf{k}}

\usepackage{amsmath}

\begin{document}

\title[Limsup is needed in two definitions of topological entropy I]{On the role of limsup in the definition of topological entropy via spanning or separation numbers.
Part I: Basic examples}
\author{Winfried Just and Ying Xin, Ohio University}

\date{\today}

\begin{abstract}
The notion of topological entropy can be conceptualized in terms of the number of forward trajectories that are
distinguishable at resolution~$\eps$ within~$T$ time units.
It can then be formally defined as a limit of a limit superior that involves either covering numbers, or separation numbers, or
spanning numbers. If covering numbers are used, the limit superior reduces to a limit.  While it has been
generally believed that the latter may not necessarily be the case when the definition is based on separation or spanning numbers,
no actual counterexamples appear to have been previously known.  Here we fill this gap in
the literature by constructing such counterexamples.
\end{abstract}

\maketitle

\section{Introduction}\label{sec:intro}

\subsection{Terminology, motivation, and main result}\label{subsec:NotationEtc}
A discrete dynamical system is a pair $(X, F)$, where $X$ is a metric space with distance function~$D$ and $F: X \rightarrow X$ is continuous.

For a given dynamical system $(X, F)$ with distance function $D$, we consider the Bowen-Dinaburg definition of topological entropy. This was first published by E. Dinaburg
in~\cite{Dinaburg}, where the author attributes the idea to unpublished work of Kolmogorov. Slightly later, but independently, the same definition was introduced and studied by R. Bowen
in~\cite{Bowen, BowenA}.\footnote{The equivalence between~\eqref{eqn:h-def} and the definition of Adler, Konheim, and McAndrew~\cite{AKM} is attributed in~\cite{Bowen} to~\cite{BowenA}, where it is
Remark (4.6). However, \cite{Bowen} has somewhat more material on this notion than~\cite{BowenA} and actually appeared slightly earlier.}
It conceptualizes topological entropy $h$ as
\begin{equation}\label{eqn:h-def}
h = h(X, F) =  \lim_{\eps \rightarrow 0^+}  \limsup_{T \rightarrow \infty} \frac{\ln N_T(\eps,D)}{T},
\end{equation}
where $N_T(\eps,D)$ measures the complexity of the system at resolution $\eps$
within $T$ steps. One can interpret~$N_T(\eps,D)$ as a covering
number~$cov(X, \eps, D_T)$, a separation number~$sep(X, \eps, D_T)$, or a spanning number~$span(X, \eps, D_T)$, as in the definitions below.

Fix any $\eps > 0$, metric space~$(X,D)$ and continuous $F: X \rightarrow X$.

\noindent
\begin{itemize}
\item  A cover~$\mathcal{U}$ of~$X$
 will be called an \emph{$\eps$-cover of~$X$}, if  for each~$U \in \mathcal{U}$ and all $x, y \in U$ the
 inequality $D(x,y) < \eps$ holds.
Then the \emph{covering number}~$cov(X, \eps, D)$ is the minimum size of an $\eps$-cover of~$X$.

\item A set of points $x_1, x_2, \dots , x_n \in X$ is said to be \emph{$\eps$-separated}, if $D(x_i, x_j) \geq \eps$ for all $1 \leq i \neq j \leq n$.
The \emph{separation number} $sep(X, \eps, D)$ is the  maximum size of an $\eps$-separated subset of~$X$.

\item A set of points $x_1, x_2, \dots , x_n \in X$ is called an \emph{$\eps$-spanning set of $X$}, if for all $x \in X$ there exists $1 \leq i \leq n$ such that $D(x_i, x) < \eps$.
The \emph{spanning number} $span(X, \eps, D)$ is the minimum size of an $\eps$-spanning subset of~$X$.

\item When~$(X, D)$ is compact (or just totally bounded), then $cov(X,\eps, D)$, \\
$sep(X,\eps, D)$, and  $span(X, \eps, D)$ always exist and are finite.

\item When $X$ is implied by the context, we will use the simplified notations~$cov(\eps, D)$,
$sep(\eps, D)$, and  $span(\eps, D)$.

\item For $x, y \in X$ and an integer $T \geq 1$ one can define
\begin{equation*}
D_T(x, y) = \max\{D(F^t(x), F^t(y)):  \ t \in \{0, 1, \dots , T-1\}\}.
\end{equation*}
These functions are metrics on~$X$; for compact~$X$ they are equivalent to~$D$.

\item A subset $A \subset X$ is \emph{$(T, \eps)$-separated} if it is $\eps$-separated with respect to~$D_T$, and is
\emph{$(T, \eps)$-spanning} if it is $\eps$-spanning with respect to~$D_T$.
\end{itemize}

We will also use the following conventions in our notation:
\begin{itemize}
\item The size of a finite set~$A$ will be denoted by~$|A|$.
\item A positive integer~$n$ will be identified with the set $\{0, \ldots , n-1\}$.  In particular, ${}^T\{0,1\}$ is the set of all functions
from $T = \{0, \dots, T-1\}$ into~$\{0,1\}$.
\item In contrast, $[n] = \{1, \dots , n\}$.
\item $diam(X, D)$ will denote the diameter of~$X$ with respect to~$D$.
\item The symbol~$\sigma$ will always denote the shift operator.
\end{itemize}

The following lemma collects some well-known relevant results.

\begin{lemma}\label{lem:BasicFacts}
Let $(X,D)$ be a compact metric space, $F: X \rightarrow X$ continuous. Then for any~$\eps > 0$:
\begin{equation}\label{ineq:covsepspan}
cov(X,\eps, D) \geq sep(X,\eps, D) \geq span(X,\eps, D) \geq cov(X,2\eps,D).
\end{equation}

\begin{equation}\label{eqn:subadd}
\forall T_1, T_2 > 0 \ \ln cov(\eps, D_{T_1+T_2}) \leq \ln cov(\eps, D_{T_1}) + \ln cov(\eps, D_{T_2}).
\end{equation}

\begin{equation}\label{eqn:coveq}
\liminf_{T \rightarrow \infty} \ \frac{\ln cov(\eps, D_T)}{T} = \limsup_{T \rightarrow \infty} \ \frac{\ln cov(\eps, D_T)}{T}.
\end{equation}
\end{lemma}

The inequalities~\eqref{ineq:covsepspan} imply that it doesn't matter which version of $N_T(\eps, D)$ we use in the
definition~\eqref{eqn:h-def} of topological entropy.
Equation~\eqref{eqn:coveq} follows from the \emph{subadditivity property}~\eqref{eqn:subadd} (see, for
example, Lemma 3.1.5 of~\cite{KH} or Section~2.1 of~\cite{DownarowiczBook} for a detailed discussion of properties related to subadditivity). It implies
that if~$h$ is defined in terms of covering numbers, then the $\limsup_{T \rightarrow \infty}$
in~\eqref{eqn:h-def} can be replaced by~$\lim_{T \rightarrow \infty}$.
Now the question naturally arises:

\begin{question}\label{q:motivating}
Is it true that for every system on any compact metric space~$(X, D)$ and every given $\eps > 0$ the following equalities hold?
\begin{equation}\label{eqn:sepeq?}
\liminf_{T \rightarrow \infty} \ \frac{\ln sep(\eps, D_T)}{T} = \limsup_{T \rightarrow \infty} \ \frac{\ln sep(\eps, D_T)}{T},
\end{equation}
\begin{equation}\label{eqn:spaneq?}
\liminf_{T \rightarrow \infty} \ \frac{\ln span(\eps, D_T)}{T} = \limsup_{T \rightarrow \infty} \ \frac{\ln span(\eps, D_T)}{T}.
\end{equation}
\end{question}

By Proposition~\ref{prop:subshift=} below, these equalities will hold when $(X,F)$ is a subshift system.  However, in general
the inequalities in~\eqref{ineq:covsepspan} may be strict,
and the analogue of the subadditivity property~\ref{eqn:subadd} may fail when~$cov$ is replaced by~$sep$
or by~$span$.  Thus it has long been widely believed that the answer to both parts of
Question~\ref{q:motivating} is negative.\footnote{For example, the last paragraph on page~164 of~\cite{DownarowiczBook} and the remark that follows
Lemma 3.1.5 of~\cite{KH} (bottom of page~109) suggest as much.}  However, as far as we could determine,
no actual counterexamples were previously known.

For the reasons outlined above, the answer to Question~1 appears not to be of much practical relevance
for calculating~$h(X,F)$ of any particular system~$(X,F)$. However, given the fundamental importance of the
concept of topological entropy, it is certainly unsatisfactory from a
theoretical point of view that the necessity of using $\limsup_{T \rightarrow \infty}$ in certain versions
of its definition has not so far been substantiated by actual counterexamples.
The main goal of this paper is to fill this  gap in the literature by proving the following result:

\begin{theorem}\label{thm:main}
There exists a system $(X^-, F)$ with a metric~$D$ on~$X^-$ such that:

\smallskip

\noindent
(i) $X^-$ is compact wrt~$D$ and $F: X^- \rightarrow X^-$ is a homeomorphism.

\smallskip

\noindent
(ii) For some $\eps > 0$ we have
\begin{equation}\label{eqn:not-eq}
\liminf_{T \rightarrow \infty} \frac{\ln span (X^-, \eps, D_T)}{T} < \limsup_{T \rightarrow \infty} \frac{\ln span (X^-, \eps, D_T)}{T}.
\end{equation} and
\begin{equation}\label{eqn:not-eqsep}
\liminf_{T \rightarrow \infty} \frac{\ln sep (X^-, \eps, D_T)}{T} < \limsup_{T \rightarrow \infty} \frac{\ln sep (X^-, \eps, D_T)}{T}.
\end{equation}

\smallskip

\noindent
(iii) $h(X^-, F) < \infty$.
\end{theorem}

\subsection{Some related results and open problems}\label{subsec:problems}

While Theorem~\ref{thm:main} gives a complete negative answer to~Question~\ref{q:motivating}, it also raises several
related problems.  The most natural perhaps is whether Equation~\eqref{eqn:sepeq?} implies
Equation~\eqref{eqn:spaneq?}, or \emph{vice versa.}

As long as we focus on a single~$\eps >0$, then our example for Theorem~\ref{thm:main}
can be easily modified to an example where~\eqref{eqn:sepeq?} fails but~\eqref{eqn:spaneq?}
holds: As the~$\eps$ of Theorem~\ref{thm:main} is also the diameter of~$X^-$  in our
construction, one can simply add a fixed point~$x^*$ of~$F$ to the space~$X^-$ and make its distance from
all other points~$\frac{\eps}{2}$.  This operation will not alter $sep(X^-, \eps, D_T)$, but it will
make $span(X^-, \eps, D_T) = 1$ for
all~$T >0$.  The question becomes more interesting if we interpret it with the existential quantifier for
failure of~\eqref{eqn:sepeq?} and the universal quantifier for~\eqref{eqn:spaneq?}.  The following result shows that even under this interpretation
Equation~\eqref{eqn:spaneq?} does not imply  Equations~\eqref{eqn:sepeq?}.

\begin{theorem}\label{thm:main-seponly}
There exist systems $(X, F)$ and $(W, F\restrict W)$ with a metric~$D$ \newline
on~$X \supset W$ such that:

\smallskip

\noindent
(i) $X, W$ are compact wrt~$D$ and  $F: X \rightarrow X$ as well as $F\restrict W: W \rightarrow W$ are homeomorphisms.

\smallskip

\noindent
(ii) For some $\eps > 0$ we have
\begin{equation}\label{eqn:not-eq-both}
\begin{split}
\liminf_{T \rightarrow \infty} \frac{\ln sep (X, \eps, D_T)}{T} &< \limsup_{T \rightarrow \infty} \frac{\ln sep (X, \eps, D_T)}{T},\\
\liminf_{T \rightarrow \infty} \frac{\ln sep (W, \eps, D_T)}{T} &< \limsup_{T \rightarrow \infty} \frac{\ln sep (W, \eps, D_T)}{T}.
\end{split}
\end{equation}

\smallskip

\noindent
(iii) For all $\delta > 0$,
\begin{equation}\label{eqn:eqspan=delta}
    \liminf_{T\rightarrow\infty} \frac{\ln span(W,\delta, D_T)}{T} = \limsup_{T\rightarrow\infty} \frac{\ln span(W,\delta, D_T)}{T}.
\end{equation}
Specifically, for some~$\delta^*$ with $0 < \delta^* < \eps$
\begin{itemize}
\item[(iiia)] If $\delta > \delta^*$,
then $\lim_{T\rightarrow\infty} \frac{\ln span(X,\delta, D_T)}{T} = \lim_{T\rightarrow\infty} \frac{\ln span(W,\delta, D_T)}{T} = 0.$
\item[(iiib)] If $\delta \leq \delta^*$,
then $\lim_{T\rightarrow\infty} \frac{\ln span(W,\delta, D_T)}{T} = \ln{2}.$
\end{itemize}

\smallskip

\noindent
(iv) $h(X, F) = \infty$ while $h(W, F\restrict W) = \ln 2$.

\smallskip

\noindent
(v) The system $(X, F)$ is not topologically transitive, while the system $(W, F \restrict W)$ is topologically transitive.
\end{theorem}

\begin{remark}\label{rem:On-Thm-seponly}
It is quite possible that the analogue of~\eqref{eqn:eqspan=delta}  holds for all~$\delta$ also in the system~$(X, F)$.  By~(iiia) this is true for $\delta > \delta^*$. However, for $\delta \leq \delta^*$ the calculations of~$span(X, \delta, D_T)$ become very tedious. Side-stepping them  by
considering the restriction of the system to a certain forward-invariant subset~$W$ of~$X$ provided the added bonus of a topologically transitive example with finite entropy.

Let us also mention that the state space~$X^-$ constructed in the proof Theorem~\ref{thm:main} is a subspace of the space~$X$ of Theorem~\ref{thm:main-seponly}.  The function~$F$ of the former theorem is the restriction of the function~$F$ of the latter to~$X^-$.  However, the metrics~$D$ are subtly different, although constructed according to the same general definition of what we call \JX-metrics. We use the same letter for them to streamline arguments that rely exclusively on their shared properties.
\end{remark}

Similarly, \eqref{eqn:spaneq?} may fail for some resolution~$\eps$, while~\eqref{eqn:sepeq?} holds for all resolutions:

\begin{theorem}\label{thm:main-spanonly}
There exists a system $(Z, H)$ with a metric~$\rho$ on~$Z$ such that:

\smallskip

\noindent
(i) $Z$ is compact wrt~$\rho$ and $H: Z \rightarrow Z$ is a homeomorphism.

\smallskip

\noindent
(ii) For some $\eps > 0$ we have
\begin{equation}\label{eqn:not-eq-spanonly}
\liminf_{T \rightarrow \infty} \frac{\ln span (Z, \eps, \rho_T)}{T} < \limsup_{T \rightarrow \infty} \frac{\ln span (Z, \eps, \rho_T)}{T}.
\end{equation}

\smallskip

\noindent
(iii) For all $\delta > 0$ we have
\begin{equation}\label{eqn:eqsep}
\liminf_{T \rightarrow \infty} \frac{\ln sep (Z, \delta, \rho_T)}{T} = \limsup_{T \rightarrow \infty} \frac{\ln sep (Z, \delta, \rho_T)}{T}.
\end{equation}

\smallskip

\noindent
(iv) $h(Z,H) < \infty$.
\end{theorem}

A related and very natural question is whether one could produce a system where the
equalities~\eqref{eqn:sepeq?} and~\eqref{eqn:spaneq?} fail for arbitrarily small~$\eps$. This question
was brought to our attention by B. Hasselblatt~\cite{Hasselblatt}.  The existence of such systems follows from
Theorems~\ref{thm:main}--\ref{thm:main-spanonly}:

\begin{corollary}\label{corol:main}
There exist systems $(X^-, F)$, $(W, F)$, $(Z,H)$ with metrics~$D$ and~$\rho$ as in Theorems~\ref{thm:main}--\ref{thm:main-spanonly} such that

\begin{itemize}
\item[(i)] Parts~(i) of
Theorems~\ref{thm:main}--\ref{thm:main-spanonly} hold.
\item[(ii)]  The inequalities in parts~(ii) of these theorems hold whenever~$\eps$ is of
the form~$\eps = 3^{-n}$ for some~$n \in \bN$.
\item[(iii3)] Equality~\eqref{eqn:eqspan=delta} of Theorem~\ref{thm:main-seponly}(iii) holds for all $\delta > 0$ in $(W,F)$.
\item[(iii4)] Equality~\eqref{eqn:eqsep} of Theorem~\ref{thm:main-spanonly}(iii) holds for all $\delta > 0$ in $(Z,H)$.
\end{itemize}
\end{corollary}

Our proof of Theorem~\ref{thm:main} required a very specialized construction, and one might ask whether
Equation~\eqref{eqn:sepeq?} and~\eqref{eqn:spaneq?} would necessarily hold in ``natural'' dynamical systems,
that is, under additional assumptions about the system. Let us discuss here just three such natural assumptions.

Let~$X$ be
a set of one- or two-sided sequences $x = (x_n)$ of symbols from a finite alphabet~$A$, closed both in the topological sense and under the \emph{subshift operator~$\sigma$}, which maps a sequence~$x = (x_n)$ to $\sigma(x) = (y_n) = (x_{n+1})$. Such spaces are called \emph{subshifts,} and we will use the phrase \emph{subshift systems} for the corresponding pairs~$(X, \sigma)$.    A standard metric on a subshift~$X$
can be defined in a slightly fanciful way as
\begin{equation}\label{eqn:D-subshift}
D(x,y) = d\left(x_{\Delta(x,y)}, y_{\Delta(x,y)}\right)k^{-\Delta(x,y)},
\end{equation}
where $k > 1$ is an integer, $d$ is the discrete metric on~$A$ that takes only values~0 and~1, and
$\Delta(x,y)$ is the first~$n$ where $x_n \neq y_n$.

\begin{proposition}\label{prop:subshift=}
Let $(X, \sigma)$ be a  subshift system with a standard metric. Then:
\begin{equation}\label{eq:covsepspan}
\forall\, T > 0 \quad cov(X,\eps, D_T) =  sep(X,\eps, D_T) = span(X,\eps, D_T).
\end{equation}

In particular,
both~\eqref{eqn:sepeq?} and~\eqref{eqn:spaneq?} hold for all $\eps > 0$ in~$(X, \sigma)$.
\end{proposition}

Proposition~\ref{prop:subshift=} is a well-known elementary result.
But what if $d$ is a different metric on~$A$ and takes more than~2 values?  Then for sufficiently large~$k$
Equation~\eqref{eqn:D-subshift}  still defines an equivalent compact metric on any~$X$ as above, but the
equalities $cov(X, \eps, D_T) = sep(X, \eps, D_T) = span(X, \eps, D_T)$ will in general no longer hold.  We don't know
the answer to the following:

\begin{question}\label{q:near-subshift}
Does there exist a subshift system~$(X, \sigma)$ such that for some choice of the metric~$d$ on the underlying alphabet~$A$
and~$k$ at least one of~\eqref{eqn:sepeq?},\eqref{eqn:spaneq?} will fail for some $\eps > 0$ if~$D$
is defined as in~\eqref{eqn:D-subshift}?
\end{question}

Another nice class of dynamical systems are the \emph{minimal} ones, that is,
systems where every point has a dense
(forward) orbit.  Minimality is a stronger property than \emph{topological transitivity,} which only requires that
there exists at least one point with a dense (forward) orbit.

\begin{question}\label{q:minimal}
Can analogues of our Theorems~\ref{thm:main}--\ref{thm:main-spanonly} be obtained for minimal systems?
\end{question}

Question~\ref{q:minimal} remains open. It was suggested to us by B. Weiss~\cite{Weiss}, who had also brought to our attention the
question whether the systems in these examples could be topologically transitive.   We have now shown that this is true in the case
in our Theorem~\ref{thm:main-seponly}.
However, at the time of this writing we still don't know whether there exist topologically transitive
examples as in our Theorems~\ref{thm:main} and~\ref{thm:main-spanonly}. We plan on addressing this
question in~\cite{PartII}.

The following fascinating question was also suggested to us by B. Weiss~\cite{Weiss}:

\begin{question}\label{q:reverse}
Consider $(X,F)$ with~$(X, d)$ a compact metric space. Is there always a metric~$D$ on~$X$ that is equivalent to~$d$,
 for which~$\lim_{T \rightarrow \infty} \ \frac{\ln sep(\eps, D_T)}{T}$
and/or $\lim_{T \rightarrow \infty} \ \frac{\ln span(\eps, D_T)}{T}$
exist for all~$\eps > 0$?
\end{question}

While we have restricted our attention here to discrete-time dynamical systems, the metrics~$D_T$ and
resulting definitions of topological entropy can be adapted to the study of flows, where $T$ can take on
arbitrary positive real values.

\begin{question}\label{q:flows}
Do any of Equations~\eqref{eqn:sepeq?},\eqref{eqn:spaneq?} always hold for
(differentiable) flows on (finite-dimensional) compact manifolds?
\end{question}

In preliminary explorations of Questions~\ref{q:motivating} we had proved that for some
differentiable flows on one-dimensional compact manifolds the
analogues of~\eqref{eqn:subadd} for separation and spanning numbers  can dramatically fail, in a way similar
to Lemma~\ref{lem:no-subadd} of Section~\ref{sec:finite}; see~\cite{finite, PartII}. Such examples were constructed by starting from
finite dynamical systems as in Lemma~\ref{lem:no-subadd}
and converting them into differentiable flows on unions of circles. It remains unclear whether our examples for
Theorems~\ref{thm:main}--\ref{thm:main-spanonly} admit similar conversions to higher-dimensional manifolds.

\subsection{Organization of the remainder of this note}\label{subsec:organization}

The purpose of this preprint is to give a complete  presentation of the proofs of
Theorems~\ref{thm:main}--\ref{thm:main-spanonly} so as to create a verifiable record of all details.
The exposition will be shortened and streamlined in a journal version.

The remainder of this note is
 organized as follows:  In Section~\ref{sec:eps} we give a derivation of Corollary~\ref{corol:main}.
In Section~\ref{sec:finite} we sketch a proof that the analogue of~\eqref{eqn:subadd}  for $sep(\eps, d_T)$ and for
$span(\eps, d_T)$ may fail even for systems with a finite state space. The purpose of this section is to introduce some important ingredients of the proofs of Theorems~\ref{thm:main} and~\ref{thm:main-seponly} in a simplified context. This lemma will also be used in
\cite{PartII} for the construction described in the discussion of Question~\ref{q:flows}.  In
Section~\ref{sec:overview} we give an outline of the proofs of
Theorems~\ref{thm:main} and~\ref{thm:main-seponly}, which is followed by the detailed constructions in
Sections~\ref{sec:earlysetup} through~\ref{sec:prove-main}. Finally, in Section~\ref{sec:prove-spanonly} we
show how this construction can be modified to obtain
Theorem~\ref{thm:main-spanonly}.

For the convenience of the reader, as an appendix we include an index with pointers to
the places where we define the many conditions and other important notions that will be referenced
throughout our arguments.

\section{Proof of Corollary~\ref{corol:main}}\label{sec:eps}

Let $(Y, G)$ denote any of the systems $(X^-, F), (X, F), (W, F\restrict W), (Z, H)$ of Theorems~\ref{thm:main}--\ref{thm:main-spanonly}, and let $\mu$ denote any of the  corresponding metrics.
Moreover, let $\eps$ be as in these theorems. We will assume that
\begin{equation}\label{eqn:diam}
diam(Y) = \eps = 1.
\end{equation}

The first equation of \eqref{eqn:diam} actually holds for all the systems constructed in our proofs
of Theorems~\ref{thm:main}--\ref{thm:main-spanonly}. But even without
going into details of these constructions we can see that assuming~\eqref{eqn:diam}  does not lead to any loss of
generality. If the diameter were greater than $\eps$,
we could redefine the metric as $\mu^* = \min \{\mu, \eps\}$.  For the second equation in~\eqref{eqn:diam}, we can
just scale the diameter to 1. Then  all properties
specified in Theorems~\ref{thm:main},~\ref{thm:main-seponly}, or~\ref{thm:main-spanonly} will continue to hold
in~$(Y, G)$ with respect to the equivalent modified metric.

\medskip

Fix any sequence $(\gamma_n)_{n \in \bN}$ of positive reals such that
\begin{equation}\label{eqn:eps-rapid}
\forall n \in \bN \ \ 2\gamma_{n+1} < \gamma_n \leq 1.
\end{equation}

For example, choosing  $\gamma_n = 3^{-n}$ will work for~\eqref{eqn:eps-rapid}.

\medskip

We will show here how one can produce
$X^-, W, Z, F, H, D,\rho$ as in Corollary~\ref{corol:main}. In this argument, we will convert~$(Y, G)$ with metric~$\mu$ that is chosen as described above into  a system~$(U, F)$ with metric~$d$ on~$U$
that retains property~(i) of Theorems~\ref{thm:main}--\ref{thm:main-spanonly} and  Equalities~\eqref{eqn:eqspan=delta} of Theorem~\ref{thm:main-seponly}(iii)
or~\eqref{eqn:eqsep} of Theorem~\ref{thm:main-spanonly} if applicable. Moreover, the inequalities or inequality in part~(ii) of the relevant theorem will
hold in~$(U,F)$ with respect to the metric~$d$ for all choices of~$\eps$ as specified in part~(ii) of Corollary~\ref{corol:main}.
The corollary itself then follows by renaming $(U, F)$ and~$d$ back to the original labels for the objects that we started with in our construction.
The conversion of~$(Y, G)$ with metric~$\mu$ into~$(U,F)$ with metric~$d$ will be accomplished by applying the following
operations to $(Y, G)$ and~$\mu$:

\begin{itemize}
\item \emph{Scaling:}

For each~$n$ we produce a copy $Y_n$ of $Y$ and a metric~$\nu^n$ on~$Y_n$ such that
$\nu^n(x, x') = \gamma_n\mu(x, x')$.  The function~$G_n$ will be the same as~$G$, but on the
copy~$Y_n$ of~$Y$.  This assures that for all~$n, T \in \bN$ with $T > 0$ and~$\delta > 0$:
\begin{equation}\label{eqn:sepspan-copy}
\begin{split}
sep\left(Y_n, \delta, \nu^n_T\right) &= sep\left(Y, \frac{\delta}{\gamma_n}, \mu_T\right), \\
span\left(Y_n, \delta, \nu^n_T\right) &= span\left(Y, \frac{\delta}{\gamma_n}, \mu_T\right).
\end{split}
\end{equation}

 \item \emph{Amplifying:}

 For each $n$ we choose a finite alphabet~$A_n$ of suitable size (see~\eqref{eqn:An-size} below)
 and let $P_n$ be the product of~$Y_n \times {}^\bZ A_n$ with the corresponding full two-sided shift.
 A metric~$d^n$ on $P_n$ will be defined by:
 \begin{equation*}\label{eqn:DAn}
 d^n((x, a), (x', a')) = \max\{\nu^n(x, x'),   \gamma_n 2^{-\Delta(a,a')}\},
 \end{equation*}
 where $\Delta(a, a')$  marks the first (under a suitable enumeration of~$\bZ$) place where the sequences~$a, a'$ differ and is defined so
 that for all~$L \in \bN$:
 {\small
 \begin{equation}\label{eqn:delta-prop}
 \Delta(a,a') \leq L \ \Leftrightarrow \
 a\restrict \left( -\left\lfloor {\frac{L}{2}}\right\rfloor,  \dots, \left\lfloor {\frac{L+1}{2}}\right\rfloor\right) \neq a'\restrict \left( -\left\lfloor {\frac{L}{2}}\right\rfloor,  \dots, \left\lfloor {\frac{L+1}{2}}\right\rfloor\right).
 \end{equation}
}

We let $G^+_n$ be the product map of~$G_n$ and the shift operator~$\sigma$.

For every~$\gamma_0 \geq \delta > 0$, we define $n(\delta) = \max \{n: \ \gamma_n \geq \delta\}$.

\begin{proposition}\label{prop:spansepYmPm}
Assume~$\gamma_0 \geq \delta > 0$. For each $0 \leq m \leq n(\delta)$, there exists a constant $L_m(\delta) \in \bN$ such that for all~$T \in \bN$:
\begin{equation}\label{eqn:spansepYmPm}
\begin{split}
    span(P_m, \delta, d_T^m) &= |A_m|^{T+L_m(\delta)}span(Y_m, \delta, \nu_T^m),\\
    sep(P_m, \delta, d_T^m) &= |A_m|^{T+L_m(\delta)}sep(Y_m, \delta, \nu_T^m).
\end{split}
\end{equation}
\end{proposition}

\noindent
\textbf{Proof:}  Fix any $0 \leq m \leq n(\delta)$ and let $L_m(\delta) > 0$ be such that
\begin{equation}\label{eqn:Lm-delta}
\forall \Delta \in \bN \qquad \Delta \leq L_m(\delta) - 1\ \ \Leftrightarrow\ \ \gamma_m2^{-\Delta} \geq \delta.
\end{equation}

Let $T > 0$ be fixed throughout the remainder of this proof.

For $a \in {}^\bZ A_n$ define
\begin{equation*}
    \varphi_a = a\restrict \left( -\left\lfloor {\frac{L_m(\delta)-1}{2}}\right\rfloor,  \dots, T+ \left\lfloor {\frac{L_m(\delta)}{2}}\right\rfloor\right).
\end{equation*}

Let $R_m = \{\varphi_a: \ a \in  {}^\bZ A_m\}$.  Then $|R_m| =  |A_m|^{T+L_m(\delta)}$.

Note that by~\eqref{eqn:delta-prop} and~\eqref{eqn:Lm-delta}, for all $a, a' \in {}^\bZ A_m$ we have:
\begin{equation}\label{eqn:varphia-Delta}
\begin{split}
\gamma_m2^{-\Delta(a, a')} \geq \delta \quad \Leftrightarrow \quad \Delta(a,a') < L_m(\delta) \quad &\Rightarrow \quad \varphi_a \neq \varphi_{a'},\\
\exists 0 \leq t < T \ \  \gamma_m2^{-\Delta(\sigma^t(a), \sigma^t(a'))} \geq \delta \quad &\Leftrightarrow \quad \varphi_a \neq \varphi_{a'}.
\end{split}
\end{equation}

Fix a subset  $A \subset {}^{\bZ}A_m$ of \emph{representatives} such that
\begin{equation*}
|A| = |A_m|^{T+L_m(\delta)} \qquad  \mbox{and} \qquad  \forall \varphi \in R_m\, \exists a \in A \quad \varphi_a  = \varphi.
\end{equation*}

Let $Q$ be any $(T, \delta)$-spanning set of $Y_m$ of minimal size.   Then for all $(x, a) \in P_m$, there exists $(x', a') \in Q \times A$  such that:

\smallskip

\begin{itemize}
\item $\nu_T^m(x,x') < \delta$ by the choice of~$Q$.
\item $\varphi_a = \varphi_{a'}$,
so that $\gamma_m2^{-\Delta(\sigma^t(a), \sigma^t(a'))} < \delta$ for all $0 \leq t < T$ by~\eqref{eqn:varphia-Delta}.
\end{itemize}

\smallskip

Thus, $d_T^m((x, a), (x', a')) < \delta$, and it follows that $Q\times A$ is $(T,\delta)$-spanning in $P_m$.  We have shown that
\begin{equation}\label{ineq:span-leq-prod}
    span(P_m, \delta, d_T^m) \leq |Q \times A| = |A_m|^{T+L_m(\delta)}span(Y_m,\delta,\nu_T^m).
\end{equation}

Similarly, let $S$ be any $(T, \delta)$-separated set of $Y_m$ of maximal size.   Then for all $(x, a) \neq (x', a') \in S \times A$
 we either have $x \neq x'$ or $a \neq a'$. In both cases the inequality $d_T^m((x,a),(x',a')) \geq \delta$ holds:
\begin{itemize}
\item If $x \neq x'$, then $d_T^m((x,a),(x',a')) \geq \nu_T^m(x,x') \geq \delta$ by the choice of~$S$.
\item If $a \neq a'$, then $\varphi_a \neq \varphi_{a'}$ by the choice of~$A$, and the inequality $d_T^m((x,a),(x',a')) \geq \delta$ follows from~\eqref{eqn:varphia-Delta}.
\end{itemize}

\smallskip

It follows that $S\times A$ is $(T,\delta)$-separated in $P_m$.  We have shown that
\begin{equation}\label{ineq:sep-geq-prod}
    sep(P_m, \delta, d_T^m) \geq |S \times A| = |A_m|^{T+L_m(\delta)}sep(Y_m,\delta,\nu_T^m).
\end{equation}

\smallskip

On the other hand, assume towards a contradiction that
\begin{equation*}
    sep(P_m, \delta, d_T^m) >  |A_m|^{T+L_m(\delta)}sep(Y_m,\delta,\nu_T^m).
\end{equation*}
That is, there exists $P \subset P_m$ of size $|P| > |A_m|^{T+L_m(\delta)}sep(Y_m,\delta,\nu_T^m)$  that is $(T, \delta)$-separated.  Then by the Pigeonhole Principle there exist $\varphi \in R_m$ and $B \subset P$, that is still $(T,\delta)$-separated in $P_m$, satisfying
\begin{equation}\label{eqn:B-phi}
\begin{split}
&|B| > sep(Y_m,\delta,\nu_T^m),\\
&\forall (x,a) \in B\ \ \varphi_a = \varphi.
\end{split}
\end{equation}

By the first line of~\eqref{eqn:B-phi} there exist $(x,a) \neq (x', a') \in B$ such that $\nu_T^m(x, x') < \delta$, and  the
second line of~\eqref{eqn:B-phi} implies together with~\eqref{eqn:varphia-Delta} that
$d_T^n((x,a),(x',a')) < \delta$, which contradicts our assumption.  Together with~\eqref{ineq:sep-geq-prod}, this implies the equality
\begin{equation*}
    sep(P_m, \delta, d_T^m) = |A_m|^{T+L_m(\delta)}sep(Y_m,\delta,\nu_T^m).
\end{equation*}

\smallskip

Similarly, let $P \subset P_m$ be $(T, \delta)$-spanning in~$P_m$. For each $\varphi \in R_m$, let
\begin{equation*}
B_\varphi = \{(x, a) \in P: \ \varphi_a = \varphi\}.
\end{equation*}

Note that if $d_T^m((x,a),(x',a')) < \delta$ for some~$(x,a) \in P$ and $(x', a') \in P_m$, then it follows from the definition of~$d^m_T$
and~\eqref{eqn:varphia-Delta}
that~$\nu_T^m(x,x') < \delta$ and~$(x,a) \in B_{\varphi_{a'}}$.  In particular, for each~$\varphi \in R_m$, the set $\{x: \  \exists a \in {}^\bZ A_m \
(x,a) \in B_{\varphi}\}$ must be $(T, \delta)$-spanning in~$Y_m$.  Since the sets~$B_\varphi$ are pairwise disjoint, this implies the inequality
\begin{equation*}
    span(P_m, \delta, d_T^m) \geq |A_m|^{T+L_m(\delta)}span(Y_m,\delta,\nu_T^m),
\end{equation*}
which together with~\eqref{ineq:span-leq-prod} implies
\begin{equation*}
    span(P_m, \delta, d_T^m) = |A_m|^{T+L_m(\delta)}span(Y_m,\delta,\nu_T^m). \qquad \Box
\end{equation*}

In view of~\eqref{eqn:spansepYmPm}, we obtain the following:

\begin{corollary}\label{corol:prod-sepspan}
Assume~$\gamma_0 \geq \delta > 0$. For each $0 \leq m \leq n(\delta)$, there exists a constant $L_m(\delta) \in \bN$ such that for all~$T >0$: {\small
 \begin{equation*}\label{eqn:sepspan-A}
 \begin{split}
\frac{\ln sep\left(P_m, \delta, d^m_T\right)}{T} &=
\frac{\ln  \left(|A_m|^{T+L_m(\delta)} sep\left(Y, \frac{\delta}{\gamma_m}, \mu_T\right)\right)}{T} \\
&=
\ln |A_m| + \frac{L_m(\delta)}{T}\ln |A_m| + \frac{\ln sep\left(Y, \frac{\delta}{\gamma_m}, \mu_T\right)}{T}, \\
\frac{\ln span\left(P_m ,\delta, d^m_T\right)}{T} &=
\frac{\ln  \left(|A_m|^{T+L_m(\delta)}  span\left(Y, \frac{\delta}{\gamma_m}, \mu_T\right)\right)}{T}  \\
&=
\ln |A_m| + \frac{L_m(\delta)}{T}\ln |A_m| + \frac{\ln span\left(Y, \frac{\delta}{\gamma_m}, \mu_T\right)}{T}.
\end{split}
\end{equation*} }
 \end{corollary}

 \medskip

 \item \emph{Combining:}

 Wlog we may assume that the sets $P_n$ are pairwise disjoint.  We define:

 \smallskip

 \begin{itemize}
 \item $U = \{x^*\} \cup \bigcup_{n \in \bN} P_n$.
 \item $U_n = \bigcup_{m = 0}^n P_m$.
 \item $U^n = \{x^*\} \cup \bigcup_{m = n+1}^\infty P_m$.
 \item $F(x) = G^+_n(x)$ for $x \in P_n$, and $F(x^*) = x^*$.
 \item $d(x, x') = d^n(x, x')$ if $x, x' \in P_n$ for some~$n$.
 \item $d(x, x^*) = d(x^*, x) = 2\gamma_n$ if $x \in P_n$ for some~$n$.
 \item $d(x, x') = 2\gamma_k$, where $k = \min \{n, m\}$,  if $x \in P_n, x' \in P_m$ for some~$n, m$ with $n \neq m$.
 \item $d(x^*, x^*) = 0$.
 \end{itemize}
\end{itemize}
\medskip

\begin{proposition}\label{prop:combining}
For each $n, T \in \bN$ with $T > 0$ the following hold:

\smallskip

\begin{itemize}
\item[(a)] $d$ is a metric on~$U$ that makes~$U$ compact.

\smallskip

\item[(b)] $F$ is a homeomorphism.

\smallskip

\item[(c)] For $\delta > \gamma_0$ we have
\begin{equation}\label{eqn:s-in-12}
\begin{split}
&span\left(U, \delta, d_T\right) \in \{1,2\} \quad \mbox{and} \quad  sep\left(U, \delta, d_T\right) \in \{1,2\},\\
&\lim_{T \rightarrow \infty} \frac{\ln span (U, \delta, d_T)}{T} \ = \
\lim_{T \rightarrow \infty} \frac{\ln sep (U, \delta , d_T)}{T} \ = \ 0.
\end{split}
\end{equation}

\smallskip

\item[(d)] For every~$\gamma_0 \geq \delta > 0$, there exist $\xi \in \{1, 2\}$ and nonnegative integers
$(L_m(\delta))_{m=0}^{n(\delta)}$ such that
\begin{equation}\label{eqn:span++}
\begin{split}
span\left(U, \delta, d_T\right) &= \xi + \sum_{m =0}^{n(\delta)} span\left(P_m, \delta, d^m_T\right)\\
&= \xi +  \sum_{m =0}^{n(\delta)} |A_{m}|^{L_m(\delta)+T} span\left(Y, \frac{\delta}{\gamma_m}, \mu_T\right).
\end{split}
\end{equation}

\smallskip

\item[(e)] For every~$\gamma_0 \geq \delta > 0$, there exist $\xi \in \{1, 2\}$ and nonnegative integers
$(L_m(\delta))_{m=0}^{n(\delta)} \subset \bN$ such that
\begin{equation}\label{eqn:sep++}
\begin{split}
sep\left(U, \delta, d_T\right) &= \xi + \sum_{m =0}^{n(\delta)} sep\left(P_m, \delta, d^m_T\right)\\
&= \xi +  \sum_{m =0}^{n(\delta)} |A_{m}|^{L_m(\delta)+T} sep\left(Y, \frac{\delta}{\gamma_m}, \mu_T\right).
\end{split}
\end{equation}
\end{itemize}
\end{proposition}

\medskip

\noindent
\textbf{Proof:} (a) We verify that $d$  has the defining properties of a metric:

\smallskip

\begin{itemize}
\item Reflexivity:  For all $x \in U$, if $x = x^*$, then $d(x, x) = 0$ by definition.  \newline If $x \in P_n$ for some $n \in \bN$, then $d(x, x) = d^n(x, x) = 0$ since $d^n$ is a metric on $P_n$.
\item Positive definiteness:  For all $x \neq x' \in U$, if $x, x' \in P_n$ for some $n \in \bN$, then $d(x, x') = d^n(x, x') > 0$ as $d^n$ is a metric on $P_n$.  \newline
Otherwise, $d(x, x') = 2\gamma_k > 0$ for some $k \in \bN$.
\item Symmetry:  $d(x, x') = d(x', x)$ for all $x, x' \in U$ follows directly from the definition of $d$.
\item The Triangle Inequality:  For all pairwise distinct $x, x', x'' \in U$, we aim to show that
\begin{equation*}
    d(x, x') + d(x, x'') \geq d(x', x'').
\end{equation*}
It suffices to consider the following cases:
\begin{itemize}
\item $x, x', x'' \in P_n$ for some $n \in \bN$.
\begin{equation*}
\begin{split}
    d(x, x') + d(x,x'') &= d^n(x, x') + d^n(x, x'')\\
    &\geq d^n(x', x'')\\
    &= d(x', x'').
\end{split}
\end{equation*}
\item $x \in P_n$ and $x', x'' \in P_m$ for some $n \neq m \in \bN$.\newline
Let $k = \min\{n, m\}$.  Then
\begin{equation*}
    d(x, x') + d(x, x'') = 4\gamma_k \geq \gamma_m \geq d^m(x', x'') = d(x', x'').
\end{equation*}
\item $x' \in P_n$ and $x, x'' \in P_m$ for some $n \neq m \in \bN$.\newline
Let $k = \min\{n, m\}$.  Then
\begin{equation*}
    d(x, x') + d(x, x'') = 2\gamma_k + d^m(x, x'') \geq 2\gamma_k = d(x', x'').
\end{equation*}
\item $x \in P_i$, $x' \in P_j$, and $x'' \in P_k$ for some pairwise distinct $i, j, k \in \bN$.\newline
Let $h = \min\{i, j, k\}$.  Then
\begin{equation*}
    d(x,x') + d(x, x'') \geq 2\gamma_h  \geq d(x', x'').
\end{equation*}
\item $x = x^*$ and $x', x'' \in P_n$ for some $n \in \bN$.
\begin{equation*}
    d(x, x') + d(x, x'') = 4\gamma_n \geq \gamma_n \geq d^n(x', x'') = d(x', x'').
\end{equation*}
\item $x = x^*$ and $x' \in P_n$, $x'' \in P_m$ for some $n \neq m \in \bN$.\newline
Let $k = \min\{n, m\}$.  Then
\begin{equation*}
    d(x,x') + d(x, x'') = 2\gamma_n + 2\gamma_m \geq 2\gamma_k  = d(x', x'').
\end{equation*}
\item $x' = x^*$ and $x, x'' \in P_n$ for some $n \in \bN$.
\begin{equation*}
    d(x, x') + d(x, x'') = 2\gamma_n + d(x, x'') \geq 2\gamma_n = d(x', x'').
\end{equation*}
\item $x' = x^*$ and $x \in P_n$, $x'' \in P_m$ for some $n \neq m \in \bN$.\newline
Let $k = \min\{n, m\}$.  Then
\begin{equation*}
    d(x, x') + d(x, x'') = 2\gamma_n + 2\gamma_k \geq 2\gamma_m = d(x', x'').
\end{equation*}
\end{itemize}
\end{itemize}
To show compactness of $U$, consider any infinite sequence $(x_k)_{k\in\bN} \subset U$. \newline
If there exist $(x_{k_j})_{j\in\bN} \subset (x_k)_{k\in\bN}$ and $n \in \bN$ such that $(x_{k_j})_{j\in\bN} \subset P_n$, then by the compactness of $Y_n$ and ${}^{\bZ}A_n$, hence the compactness $P_n$, the infinite sequence $(x_{k_j})_{j\in\bN}$ has an infinite subsequence in $P_n$ that converges to some point in $P_n$.

Otherwise, we can pick an increasing sequence  $(n_j)_{j\in\bN} \subset \bN$ such that for each $j \in \bN$, we have
$(x_k)_{k\in\bN} \cap P_{n_j} \neq \emptyset$.  For each $j \in \bN$, choose an $x_{k_j} \in (x_k)_{k\in\bN} \cap P_{n_j}$.  The resulting subsequence $(x_{k_j})_{j\in\bN}$ of  $(x_k)_{k\in\bN}$ satisfies
\begin{equation*}
    \lim_{j\rightarrow\infty}d(x^*, x_{k_j}) = \lim_{j\rightarrow\infty} 2\gamma_{n_j} = 0.
\end{equation*}

Thus, $(x_{k_j})_{j\in\bN}$ converges to $x^* \in U$, and we conclude that $U$ is compact.

Note that topologically, $U$ is the one-point compactification of \ $\bigcup_{n \in \bN} P_n$.

\medskip

\noindent
(b)  Since $U$ is compact, it suffices to show that $F$ is a continuous bijection on $U$.
\begin{itemize}
\item $F$ is one-to-one:  For all $x \neq x' \in U$, consider the following three cases:
\begin{itemize}
\item $x = x^*$ and $x' \in P_n$ for some $n \in \bN$.
\begin{equation*}
\begin{split}
    &F(x) = x^* \notin P_m\ \mbox{for any}\ m \in \bN,\ \ F(x') = G_n^+(x') \in P_n.\\
    &\mbox{Thus} \ F(x) \neq F(x').
    \end{split}
\end{equation*}
\item $x = (y, a),\ x' = (y', a') \in P_n$ for some $n \in \bN$. \newline
Then, $y \neq y'$ or $a \neq a'$.  Since both $G_n$ and $\sigma$ are one-to-one,
\begin{equation*}
\begin{split}
    F(x) &= G_n^+(x) = (G_n(y), \sigma(a))\\
    &\neq (G_n(y'), \sigma(a')) = G_n^+(x') = F(x').
\end{split}
\end{equation*}
\item $x \in P_n$ and $x' \in P_m$ for some $n \neq m \in \bN$.
\begin{equation*}
\begin{split}
    &F(x) = G_n^+(x) \in P_n,\ \ F(x') = G_m^+(x') \in P_m.\\
     &\mbox{Thus,}\ \ F(x) \neq F(x').
\end{split}
\end{equation*}
\end{itemize}
\item $F$ is onto:  For any $x \in U$, if $x = x^*$, we have $F(x^*) = x^*$. \newline
If $x = (y, a) \in P_n$ for some $n \in \bN$, there exist $y' \in Y_n$ and $a' \in {}^{\bZ}A_n$ such that $G_n(y') = y$ and $\sigma(a') = a$, since $G_n$ and $\sigma$ are onto functions.  Thus, we have $x' = (y', a') \in P_n$ such that $F(x') = G_n^+(x') = (G_n(y'), \sigma(a')) = (y, a) = x$.  Therefore, $F$ maps~$U$ onto~$U$.
\item $F$ is continuous:  For any $\gamma > 0$, there exists $k \in \bN$ such that $\gamma_k < \gamma$.  \newline
Fix such a $k$.

By construction, for each $n \in N$ the function $G_n^+$ is continuous.  Thus for each $n \in \bN$ we can pick $\zeta_n > 0$ such that
\begin{equation*}
    d^n(G_n^+(x), G_n^+(x')) < \gamma \
    \mbox{whenever}\ d^n(x,x') < \zeta_n.
\end{equation*}
Let $\zeta = \min\{\zeta_n: 0 \leq n < k\}$ and let $\eta = \min\{\gamma_k, \zeta\}$.\newline
Then for all $x \neq x' \in U$ with $d(x, x') < \eta$,
\begin{itemize}
\item if $x = x^*$, $x' \in P_n$ for some $n \in \bN$, or $x \in P_n$, $x' \in P_m$ for some $n \neq m \in \bN$, by the definition of $d$,
\begin{equation*}
    d(F(x), F(x')) = d(x, x') < \eta \leq \gamma_k < \gamma.
\end{equation*}
\item if $x, x' \in P_n$ for some $n \geq k$,
\begin{equation*}
    d(F(x), F(x')) \leq \gamma_n \leq \gamma_k < \gamma.
\end{equation*}
\item if $x, x' \in P_n$ for some $0 \leq n < k$,
\begin{equation*}
    d(x, x') = d^n(x, x') < \eta \leq \zeta \leq \zeta_n.
\end{equation*}
It follows that \ $d(F(x), F(x')) = d^n(G_n^+(x), G_n^+(x')) < \gamma$.
\end{itemize}
Thus, $F$ is continuous on $U$.
\end{itemize}

We conclude that $F$ is a homeomorphism.

\medskip

\noindent
(c) As the first line of~\eqref{eqn:s-in-12} implies the second and $span(U, \delta, d_T) \leq  sep(U, \delta, d_T)$
by~\eqref{ineq:covsepspan} of Lemma~\ref{lem:BasicFacts}, it suffices to show that for any $\delta > \gamma_0$ and $T > 0$,
\begin{equation*}
    sep(U, \delta, d_T) \leq 2.
\end{equation*}

Assume towards a contradiction that there exist pairwise distinct $x, x', x'' \in U$ such that the set
$\{x, x', x''\}$ is $(T, \delta)$-separated.  It suffices to consider the following cases:
\begin{itemize}
\item $x, x' \in P_n$ for some $n \in \bN$.  Then
\begin{equation*}
    d_T(x, x') \leq \gamma_n \leq \gamma_0 < \delta.
\end{equation*}
\item $x = x^*$, $x' \in P_n$ and $x'' \in P_m$ for some $n < m \in \bN$.  Then
\begin{equation*}
    d_T(x, x'') = 2\gamma_m < \gamma_0 < \delta.
\end{equation*}
\item $x \in P_i$, $x' \in P_j$ and $x'' \in P_k$ for some $i < j < k \in \bN$.  Then
\begin{equation*}
    d_T(x', x'') = 2\gamma_j < \gamma_0 < \delta.
\end{equation*}
\end{itemize}
Thus $\{x, x', x''\}$ is not $(T, \delta)$-separated.  The inequality $sep(U, \delta, d_T) \leq 2$ follows.

\medskip

\noindent
(d) Fix any $\gamma_0 \geq \delta > 0$ and $T > 0$. Notice that the second equality in~\eqref{eqn:span++} follows from the first in view of Corollary~\ref{corol:prod-sepspan}.
For the proof of this first equality in~\eqref{eqn:span++}, let~$\xi$ denote $span(U^{n(\delta)}, \delta, d_T)$.

 Consider $(T, \delta)$-spanning sets $S_m \subset P_m$ for $0 \leq m \leq n(\delta)$ with respect to the corresponding metrics~$d^m$, and a
$(T, \delta)$-spanning set $S^{n(\delta)} \subset U^{n(\delta)}$ with respect to the metric~$d$. Then the definition of~$d$ implies that the union
$S = S^{n(\delta)} \cup \bigcup_{m = 0}^{n(\delta)} S_m$ is $(T, \delta)$-spanning in~$U$ with respect to~$d$, and by considering spanning sets of minimal size we obtain the inequality
\begin{equation*}
span\left(U, \delta, d_T\right) \leq \xi + \sum_{m =0}^{n(\delta)} span\left(P_m, \delta, d^m_T\right).
\end{equation*}

Conversely, let $S \subset U$ be $(T, \delta)$-spanning with respect to the metric~$d$.

Let~$S^{n(\delta)} = S \cap U^{n(\delta)}$, and for $0 \leq m \leq n(\delta)$, let~$S_m = S \cap U_m$.

Note that $\gamma_{n(\delta)+1} < \delta \leq \gamma_{n(\delta)}$.
Thus for all $0 \leq m \leq n(\delta)$,
\begin{equation*}
    \forall x \in P_m, x' \in U\backslash P_m,\ \ d_T(x, x') \geq 2\gamma_m \geq 2\gamma_{n(\delta)} > \delta.
\end{equation*}

It follows that each of the sets $S_m, S^{n(\delta)}$ is $(T, \delta)$-spanning with respect to the metric~$d$, in $U_m$ and $U^{n(\delta)}$, respectively.  Hence by the definition of~$d$,
each of the sets $S_m$ is also $(T, \delta)$-spanning
with respect to~$d^m$. This proves the inequality
\begin{equation*}
span\left(U, \delta, d_T\right) \geq \xi + \sum_{m =0}^{n(\delta)} span\left(P_m, \delta, d^m_T\right).
\end{equation*}

It remains to show that $\xi = span(U^{n(\delta)}, \delta, d_T) \in \{1, 2\}$. This follows from the inequalities
$1 \leq span(U^{n(\delta)}, \delta, d_T) \leq sep(U^{n(\delta)}, \delta, d_T)$ that hold in every dynamical system and the following observation:
\begin{equation}\label{eqn:xi-in-12}
sep(U^{n(\delta)}, \delta, d_T) \leq 2.
\end{equation}

To see why~\eqref{eqn:xi-in-12} holds, consider any subset $S^- \subset U^{n(\delta)}$ that is $(T, \delta)$-separated with respect to~$d$.
Then by the definition of~$d$ we must have
\begin{equation*}\label{ineq:maxdist-in-S-}
\forall x, x' \in S^- \quad   d_T(x,x') \leq 2\gamma_{n(\delta)+1}.
\end{equation*}

It follows that if $\delta > 2\gamma_{n(\delta) + 1}$, then we even get $sep(U^{n(\delta)}, \delta, d_T) = 1$.

If $\gamma_{n(\delta)+1} < \delta \leq 2\gamma_{n(\delta)+1}$, then we can use the observation that
\begin{equation*}
\begin{split}
&\forall x, x' \in U^{n(\delta)+1}\quad   d_T(x,x') \leq 2\gamma_{n(\delta)+2} < \gamma_{n(\delta)+1} < \delta,\\
&\forall x, x' \in P_{n(\delta)+1}\quad    d_T(x,x') \leq \gamma_{n(\delta)+1} < \delta.
\end{split}
\end{equation*}

Thus, $S^-$ can contain at most one element of $U^{n(\delta)+1}$ and at most one element of $P_{n(\delta)+1}$,
and~\eqref{eqn:xi-in-12} follows.

\medskip

\noindent
(e) This proof is analogous to the one for part~(d). Fix any $\gamma_0 \geq \delta > 0$ and $T > 0$, and notice that the second equality in~\eqref{eqn:span++} follows from the first in view of Corollary~\ref{corol:prod-sepspan}.
For the proof of this first equality in~\eqref{eqn:span++}, let $\xi$ denote $sep(U^{n(\delta)}, \delta, d_T)$. By~\eqref{eqn:xi-in-12}, $\xi \in \{1,2\}$.

 Consider $(T, \delta)$-separated sets $S_m \subset P_m$ for $0 \leq m \leq n(\delta)$ with respect to the corresponding metrics~$d^m$, and a
$(T, \delta)$-separated set $S^{n(\delta)} \subset U^{n(\delta)}$ with respect to the metric~$d$. Then the definition of~$d$ implies that the union
$S = S^{n(\delta)} \cup \bigcup_{m = 0}^{n(\delta)} S_m$ is $(T, \delta)$-separated in~$U$ with respect to~$d$, and by considering separated sets of maximal size we obtain the inequality
\begin{equation*}
sep\left(U, \delta, d_T\right) \geq \xi + \sum_{m =0}^{n(\delta)} sep\left(P_m, \delta, d^m_T\right).
\end{equation*}

Conversely, let $S \subset U$ be $(T, \delta)$-separated with respect to the metric~$d$.

Let~$S^{n(\delta)} = S \cap U^{n(\delta)}$, and for $0 \leq m \leq n(\delta)$, let~$S_m = S \cap U_m$.
By the definition of~$d$,  each of the sets $S_m$ is $(T, \delta)$-separated with respect to the metric~$d_m$, and $|S_m|$ cannot exceed
$sep(P_m, \delta, d^m_T)$. Similarly, $|S^{n(\delta)}| \leq sep(U^{n(\delta)}, \delta, d_T) = \xi$.  This proves the inequality
\begin{equation*}
sep\left(U, \delta, d_T\right) \leq \xi + \sum_{m =0}^{n(\delta)} sep\left(P_m, \delta, d^m_T\right). \qquad \Box
\end{equation*}

\bigskip

Recall that we assumed $0 < h(Y,G) < \infty$. Moreover, by definition
\begin{equation*}
\begin{split}
&\forall \delta < \delta^+ \, \forall T > 0 \ \ span(Y, \delta, \mu_T) \geq span(Y, \delta^+, \mu_T),\\
&\forall \delta < \delta^+ \, \forall T > 0 \ \ sep(Y, \delta, \mu_T) \geq sep(Y, \delta^+, \mu_T),\\
&\lim_{\delta \rightarrow 0^+} \limsup_{T \rightarrow \infty} \frac{\ln span(Y, \delta, \mu_T)}{T} =
\lim_{\delta \rightarrow 0^+} \limsup_{T \rightarrow \infty} \frac{\ln sep(Y, \delta, \mu_T)}{T} = h(Y, G).
\end{split}
\end{equation*}

In particular, for all~$\beta > 0$ and all sufficiently large~$T$ the following inequalities will hold:
\begin{equation}\label{eqn:sepspan-approx-h}
span(Y, \beta, \mu_T) \leq \left(2e^{h(Y, G)}\right)^T \quad \mbox{and} \quad sep(Y, \beta, \mu_T) \leq \left(2e^{h(Y, G)}\right)^T.
\end{equation}

\smallskip

Let $0 \leq \delta \leq \gamma_0$.
For all sufficiently large~$T$ we get from~\eqref{eqn:span++},~\eqref{eqn:sep++},  and~\eqref{eqn:sepspan-approx-h}:
\begin{equation*}\label{eqn:An-estimates}
\begin{split}
|A_{n(\delta)}|^{L_{n(\delta)}(\delta)+T} span\left(Y, \frac{\delta}{\gamma_{n(\delta)}}, \mu_T\right) &\leq span\left(U, \delta, d_T\right)\\
&\leq 2+ |A_{n(\delta)}|^{L_{n(\delta)}(\delta)+T} span\left(Y, \frac{\delta}{\gamma_{n(\delta)}}, \mu_T\right)\\
&\quad + \sum_{m =0}^{n(\delta)-1} |A_{m}|^{L_m(\delta)+T} \left(2 e^{h(Y, G)}\right)^T,\\
|A_{n(\delta)}|^{L_{n(\delta)}(\delta)+T} sep\left(Y, \frac{\delta}{\gamma_{n(\delta)}}, \mu_T\right) &\leq sep\left(U, \delta, d_T\right)\\
& \leq 2 +  |A_{n(\delta)}|^{L_{n(\delta)}(\delta)+T} sep\left(Y, \frac{\delta}{\gamma_{n(\delta)}}, \mu_T\right)\\
&\quad + \sum_{m =0}^{n(\delta)-1} |A_{m}|^{L_m(\delta)+T} \left(2 e^{h(Y, G)}\right)^T.
\end{split}
\end{equation*}

Thus if we choose the alphabets~$A_n$ in such a way that
\begin{equation}\label{eqn:An-size}
\forall n \in \bN \ \ |A_{n+1}| > 2|A_n|e^{h(Y,G)},
\end{equation}
then it follows from Corollary~\ref{corol:prod-sepspan} that
\begin{equation}\label{eqn:And-estimate}
\begin{split}
\liminf_{T \rightarrow \infty} \frac{\ln span (U, \delta, d_T)}{T} \ &= \ \ln |A_{n(\delta)}| + \liminf_{T \rightarrow \infty}
\frac{\ln span\left(Y, \frac{\delta}{\gamma_{n(\delta)}}, \mu_T\right)}{T},\\
\limsup_{T \rightarrow \infty} \frac{\ln span (U, \delta, d_T)}{T} \ &= \ \ln |A_{n(\delta)}| + \limsup_{T \rightarrow \infty}
\frac{\ln span\left(Y, \frac{\delta}{\gamma_{n(\delta)}}, \mu_T\right)}{T},\\
\liminf_{T \rightarrow \infty} \frac{\ln sep (U, \delta, d_T)}{T} \ &= \ \ln |A_{n(\delta)}| + \liminf_{T \rightarrow \infty}
\frac{\ln sep \left(Y, \frac{\delta}{\gamma_{n(\delta)}}, \mu_T\right)}{T},\\
\limsup_{T \rightarrow \infty} \frac{\ln sep (U, \delta, d_T)}{T} \ &= \ \ln |A_{n(\delta)}| + \limsup_{T \rightarrow \infty}
\frac{\ln sep \left(Y, \frac{\delta}{\gamma_{n(\delta)}}, \mu_T\right)}{T}.
\end{split}
\end{equation}

Note that for every~$\delta = \gamma_n$ we have $n(\delta) = n$ and~$\frac{\delta}{\gamma_{n(\delta)}} = \frac{\gamma_n}{\gamma_{n(\gamma_n)}} =  \frac{\gamma_n}{\gamma_{n}} = 1 = \eps$.
Let $sp$ stand for either ``span'' or ``sep.'' It follows from~\eqref{eqn:And-estimate} that the inequality
\begin{equation*}
\liminf_{T \rightarrow \infty} \frac{\ln sp (U, \delta, d_T)}{T} \ \leq  \limsup_{T \rightarrow \infty}
\frac{\ln sp(U, \delta, d_T)}{T}
\end{equation*}
will be strict for $\delta = \gamma_n$ whenever its counterpart for $\delta = \eps$ in~$(Y,G)$ is strict, and will turn into an equality for all~$\delta > 0$ whenever the same is true for its counterpart in~$(Y,G)$.
Recall that we chose $Y$ and $\mu$ as $X^-, W$, or $Z$ and~$D$ or~$\rho$ of Theorems~\ref{thm:main}--\ref{thm:main-spanonly}, respectively.
If after  performing the above
construction we rename the corresponding space~$U$ that we constructed
back to the labels of the original structures that we started with,  then~\eqref{eqn:And-estimate}, or~\eqref{eqn:s-in-12} when $\delta > \gamma_0$, implies
the equations and inequalities that are referenced in Corollary~\ref{corol:main}.  $\Box$

\bigskip

In view of~\eqref{eqn:An-size}, the above construction always gives systems with infinite topological entropy.
Also, the systems that were constructed in this
proof of Corollary~\ref{corol:main} are not topologically transitive.  The following problem remains open at the time of this
writing:

\begin{question}\label{q:sme-h-fin}
Are there examples as in Corollary~\ref{corol:main} that:

\smallskip

\noindent
(a) Have finite topological entropy?

\smallskip

\noindent
(b) Are topologically transitive?
\end{question}

We will return to this question in~\cite{PartII}.

\section{Warm-up: Failure of subadditivity for $\ln sep\left(\eps, d_T\right)$ and
$\ln span\left(\eps, d_T\right)$}\label{sec:finite}

We will prove the following result.

\begin{lemma}\label{lem:no-subadd}
There exists  a positive constant~$R^*$ such that for each positive integer~$T$ there exist a
finite dynamical system $(X_0, F)$ with
 $F$ a bijection and  a metric~$d$ on~$X_0$ so that

\smallskip

\noindent
(i) Every point in~$X_0$ is a periodic point of~$F$ with minimal period~$3T$.

\smallskip

\noindent
(ii) There exist $0 < \delta_0 < \eps_0 < 2\delta_0 < 1$ so that~$d$ takes values only in the
set~$\{0, \delta_0, \eps_0\}$.

\smallskip

For all $\eps$ with $\delta_0 < \eps \leq \eps_0$ the system satisfies:

\smallskip

\noindent
(iii) $sep\left(\eps, d_{3T}\right) = span\left(\eps, d_{3T}\right) = 3T 2^T = |X_0|$.

\smallskip

\noindent
(iv) $span\left(\eps, d_{T}\right) \leq sep \left(\eps, d_{T}\right) \leq R^* T^2$.

\smallskip

Moreover, if $T$ is chosen sufficiently large, then

\smallskip

\noindent
(v) $\exists T_1, T_2 > 0 \ \ln sep \left(\eps, d_{T_1+T_2}\right) > \ln sep\left(\eps, d_{T_1}\right) + \ln sep\left(\eps, d_{T_2}\right)$.

\smallskip

\noindent
(vi) $\exists T_1, T_2 > 0 \ \ln span \left(\eps, d_{T_1+T_2}\right) > \ln span\left(\eps, d_{T_1}\right) + \ln span\left(\eps, d_{T_2}\right)$.

\end{lemma}

\noindent
\textbf{Proof:} Parts~(v) and~(vi) follow from parts~(iii) and~(iv). To see this, let us focus on the
separation numbers.  For sufficiently large~$T$ and~$\eps$ as above we will have
\begin{equation*}\label{eqn:no-3-subadd}
3\ln sep\left(\eps, d_T\right) \leq 3(\ln R^* + 2 \ln T) < T \ln 2 <
\ln   sep\left(\eps, d_{3T}\right),
\end{equation*}
and~(iv) must be satisfied for either $T_1 = T_2 = T$ or $T_1 = T$ and $T_2 = 2T$.

It remains to prove parts~(i)--(iv).
The particular argument and notation are a little more cumbersome than strictly necessary. They have been chosen so
that they match and illustrate important ingredients of the proofs of
Theorems~\ref{thm:main} and~\ref{thm:main-seponly}.

Let us fix two positive integers~$T < T^+$ such that~$T^+ = 3T$  and positive
reals $\delta_0 < \eps_0$ such that $\delta_0 < \eps_0 < 2\delta_0 < 1$. Moreover, let $Y \subset {}^\bZ \{0,1\}$
be the set of all two-sided sequences of zeros and ones that are periodic with period~$T$.
Let $X_0$ be the set of all triples~$x = (y, 0, k)$, where $y \in Y$ and~$k \in T^+ = \{0, 1, \dots , T^+ -1\}$.

To look ahead a bit:  The set~$X_0$ here is almost the same as the set~$X_0$ that we will define in
Section~\ref{sec:Xn-construct}, except that the first coordinates~$y$ of the latter will no longer assumed to be periodic.
The middle label only serves to make the distinction from elements of~$X_n$ for $n > 0$
and is not needed here, but kept for consistency of notation.

However, periodicity is important in the current proof. It makes~$X_0$ a finite set, of cardinality~$|X_0| = 3T2^{T}$.

We define the function~$F: X_0 \rightarrow X_0$ as follows:
\begin{equation}\label{eqn:def-F0}
F((y, 0,  k)) = (\sigma(y), 0, (k+1) \ mod \ T^+),
\end{equation}
where $\sigma$ denotes the shift operator, so that
\begin{equation*}
\sigma(y)(i) = y(i+1) \qquad \mbox{for all} \ \ i \in \bZ.
\end{equation*}
Note that $F$ is a bijection such that each~$x \in X$ is periodic with minimal period~$T^+ = 3T$.  Moreover,
since $X_0$ is finite, for any metric~$d$ on~$X_0$ we obtain a compact state space, and $F$  will be a
homeomorphism. This proves part~(i).

Towards the definition of our particular metric~$d$, we first partition the interval $T^+ = \{0, \dots , T^+ -1\}$ into
three consecutive
subintervals~$I_j^0$ of length~$T$ each, where~$j\in \{1, 2, 3\}$.

Second, we associate with each element of~$X_0$ a function $\Phi((y, 0, k)) \in {}^{T^+}\{0,1\}$.
 These functions can be defined as
\begin{equation}\label{eqn:Phi-def-0}
\Phi((y, 0, k)) = (y(-k), y(-k+1), \dots , y(-k+T^+ -1)).
\end{equation}

By periodicity that we assumed here we will always have
\begin{equation}\label{eqn:Phi-invariant-F}
\Phi((y,  0, k)) = \Phi(F((y, 0,  k))).
\end{equation}

Without periodicity, \eqref{eqn:Phi-invariant-F} will sometimes fail, but a suitably modified version of  it will remain true in
the proofs of Theorems~\ref{thm:main} and \ref{thm:main-seponly} (see Proposition~\ref{prop:Phi-const} below).

Third, let ${}^{T^+}\{0,1\}$ denote the set of all functions with domain $\{0, 1, \dots , T^+ -1\}$ that
take values in the set $\{0,1\}$. For any set
$S$, let $[S]^2$ denote the set of all unordered pairs  of
different elements of~$S$. Recall the notation $[C] = \{1, 2, \dots, C\}$.
A \emph{coloring of~$[S]^2$ with~$C$ colors} is simply a
function $c: [S]^2 \rightarrow [C]$.
 We will fix a coloring of $\left[{}^{T^+}\{0,1\}\right]^2$ with $C = 3$ colors that has suitable properties
(to be specified later).   With~$c$ acting as a parameter,
we define a distance $d((y,0,k), (y',0,k'))$ as follows:
\begin{itemize}
\item[(d1)] If $k \neq k'$, then $d((y, 0, k), (y', 0, k')) = \eps_0$.
\item[(d2)] If $k = k'$ and $y = y'$, then $d((y, 0, k), (y', 0, k')) = 0$.
\item[(d3)] If $k = k'$ and $y \neq y'$, then
      \begin{itemize}
      \item[(d31)] If $y(0) = y'(0)$, then $ ((y,  0, k), (y', 0, k')) = \delta_0$.
      \item[(d32)] If $y(0) \neq y'(0)$, then we let $\varphi = \Phi((y,0,k))$ and $\psi = \Phi((y',0, k'))$
      and define:
            \begin{itemize}
            \item $d ((y, 0, k), (y',  0, k')) = \eps_0$ if $k  \in I_j$ and $c(\varphi, \psi) = j$.
            \item $d((y, 0, k), (y', 0, k')) = \delta_0$ if $k \in I_j$ and $c(\varphi, \psi) \neq j$.
            \end{itemize}
     \end{itemize}

\end{itemize}

This function~$d$ is similar to what we will call a \emph{\JXn-metric} in
Section~\ref{sec:Xn-construct}, except that there we will need a different version
of~(d31).

Note that~$d$ takes only values in
the set $\{0, \delta_0, \eps_0\}$ and is a metric on~$X_0$ for any choice of the coloring~$c$.
Reflexivity and symmetry of~$d$ are immediate from the definition; the Triangle Inequality follows from our assumption that $2\delta_0 > \eps_0$.
Thus~$d$  satisfies part~(ii) of the lemma.

 Now assume $(y, 0, k) \neq (y', 0, k')$ and~$k = k'$,
so that clause~(d32) of the definition of~$d$ applies.  As each relevant~$y$ is periodic with period~$T$,
for each~$j \in \{1, 2, 3\}$ there will be an
$i \in I_j$ with~$y(i) \neq y'(i)$. Since $k$ periodically shifts under the action of~$F$, but~$\varphi$
and~$\psi$ remain fixed in view of~\eqref{eqn:Phi-invariant-F}, for some $t$ with $0 \leq t < T^+ -1$ we must then have
$d(F^t((y, 0, k)), F^t((y', 0, k'))) = \eps_0$ according to clause~(d32).  This, together with clause~(d1) implies
that the entire set~$X_0$ is $(T^+, \eps_0)$ separated.  Since $T^+ = 3T$,  part~(iii) of
the lemma follows.

Note that no special property of the coloring~$c$ was used in the derivation of~(iii).

Now let $\eps > \delta$ and consider a $(T, \eps)$-separated subset $Z\subset X_0$ such that for some fixed~$k$ all elements
of~$Z$ are of the form~$(y, 0, k)$. Then for all $(y,0, k) , (y', 0, k) \in Z$ with $y \neq y'$ there must be
some~$t < T$ such that the inequality
\begin{equation*}
\eps \leq d(F^t((y, 0, k)), F^t((y', 0, k))) = \eps_0
\end{equation*}
is witnessed by clause~(d32).
As $\varphi$ and $\psi$ remain constant under the action of~$F$ in view of~\eqref{eqn:Phi-invariant-F}, and
as the third coordinates $(k+t) \mod T^+$ of $F^t((y, 0, k))$ can take values in only two of the three
intervals~$I_j$ while~$t$ ranges from~$0$ to~$T-1$, this in
turn implies that the restriction of~$c$ to the
set $\left[\{\Phi((y, 0, k)): \, (y, 0, k) \in Z\}\right]^2$ can take at most two of the three possible values.  In other
words, the set   $\{\Phi((y, 0, k)): \, (y, 0, k) \in Z\}$ must be \emph{$\leq 2$-chromatic} for~$c$.  Since
all elements of~$Z$ have the same third coordinate, the
restriction of $\Phi$ to~$Z$ is one-to-one. Thus we can conclude that~$|Z|$ cannot exceed the maximum size of a
$\leq 2$-chromatic set for~$c$.

Now we make use of the following fact. Here we assume $n \geq 2$ to avoid degenerate meanings of the term ``coloring.''

\begin{proposition}\label{lem:c-exists-basic}
Let $n \geq 2$ and let $R = \frac{1}{\ln \sqrt{3} - \ln \sqrt{2}}$.  Then there exists a
coloring $c: [n]^2 \rightarrow [3]$ for which every $\leq 2$-chromatic set has size at most $R \ln n$.
\end{proposition}

We let $T$ be a positive integer, choose~$R$ as in Proposition~\ref{lem:c-exists-basic},  and a coloring $c:{}^{3T}\{0,1\} \rightarrow [3]$ without
a~$\leq 2$-chromatic subset of size $> 3TR \ln 2$.
The argument that immediately precedes this proposition shows that if $Z$ is a $(T, \eps)$-separated
subset of~$X_0$, then for every fixed~$k < T^+ = 3T$, the set $Z$ can contain at most $R T \ln 2$ elements
of the form~$(y,k)$.
Thus $Z$ itself  can have at most $3R T^2 \ln 2$ elements.  Part~(iv) of Lemma~\ref{lem:no-subadd} then
follows for the choice $R^* = 3R\ln 2$.

Analogues of Proposition~\ref{lem:c-exists-basic} for more sophisticated colorings will be derived in later
sections.  In order to illustrate how these arguments work, we include here the more basic proof of
 Proposition~\ref{lem:c-exists-basic}.

\medskip

\noindent
\textbf{Proof of Proposition~\ref{lem:c-exists-basic}:} Our statement of the proposition is a consequence of its following version:

\begin{proposition}\label{prop:color-3}
Let $m$ be a positive integer.
Then for every $2 \leq n \leq \left(\frac{\sqrt{3}}{\sqrt{2}}\right)^{m-1}$ there exists a coloring $c: [n]^2 \rightarrow [3]$ that does not have a $\leq 2$-chromatic set of size~$\geq m$.
\end{proposition}

\noindent
\textbf{Proof:} Fix $m$ as in the assumption and let $2 \leq n \leq \left(\frac{\sqrt{3}}{\sqrt{2}}\right)^{m-1}$. Since every subset of a $\leq 2$-chromatic set is $\leq 2$-chromatic, we only need to show existence of a coloring $c: [n]^2 \rightarrow [3]$ without a $\leq 2$-chromatic set of size exactly $m$.

Let us consider the set of all possible colorings $c: [n]^2 \rightarrow [3]$ with the uniform distribution.  Let $F_{i, j, p}$ denote the event that $c(i, j) = p$ when we draw  coloring $c$ randomly from this distribution. Each of the events $F_{i, j, p}$ will then have probability~$\frac{1}{3}$, and the events are independent for different $\{i, j\} \in [n]^2$.

Let $M \subset [n]$ be of size exactly~$m$, that is, $M \in [n]^m$. Consider the r.v. (random variable) $\xi_M$ that for a randomly drawn~$c$ takes the value~$1$ if $M$ is $\leq 2$-chromatic and takes the value~0 otherwise.  Then the expected value of $\xi_M$ is given by

\begin{equation}\label{eqn:ExiM-2}
E(\xi_M) = P(\xi_M = 1) = 3 \left(\frac{2}{3}\right)^{m(m-1)/2}.
\end{equation}

Now let $\xi$ be the r.v. that counts the number of $\leq 2$-chromatic subsets of size~$m$:

\begin{equation}\label{eqn:xi-def-3}
\xi = \sum_{M \in [n]^m} \xi_m.
\end{equation}

Then
\begin{equation}\label{eqn:Exi-est-3}
\begin{split}
E(\xi) &= \sum_{M \in [n]^m} E(\xi_m) = 3\binom{n}{m} \left(\frac{2}{3}\right)^{m(m-1)/2} < n^m\left(\frac{2}{3}\right)^{m(m-1)/2} \\
&\leq
\left(\left(\frac{3}{2}\right)^{(m-1)/2}\left(\frac{2}{3}\right)^{(m-1)/2}\right)^m = 1.
\end{split}
\end{equation}

Since the expected value of~$\xi$ is less than~1, we must have $\xi(c) = 0$ for at least one coloring~$c$,
which witnesses the result  claimed in the proposition. $\Box$ $\Box$ $\Box$

\bigskip

\begin{remark}\label{rem:afterLem-finite}
Our constructions for
Theorems~\ref{thm:main} and~\ref{thm:main-seponly}  are based on infinite products of
systems~$(X_n, F_n)$ that are somewhat similar to the ones in the proof of
Lemma~\ref{lem:no-subadd}.
However, we will need to drop the assumption of periodicity  to have enough
candidates for inclusion in large separated subsets of the product space. Similarly, letting $T^+ = 3T$ will no longer work;
we need to give ourselves more flexibility by
picking~$T^+(n) = C(n)T(n)$ for some carefully chosen positive integers~$C(n)$ and  colorings~$c_n$
with~$C(n)$ colors. Moreover, the construction of large $(T^+, \eps_0)$-separated sets in our proof for
Lemma~\ref{lem:no-subadd}(iii) was based on finding
one~$t < T^+$ for which the coloring~$c$ ``takes the right value.'' But when we consider direct
products of infinitely many such systems, we will need to make sure that there is one~$t$ where all
 colorings~$c_n$ that define the metrics on the coordinates  ``take the right value'' all
 \emph{simultaneously.}  Our work in Subsection~\ref{subsec:lower-sep-T+n} shows how this
 can be achieved; the results of Subsection~\ref{subsec:lower-span-T+n} serve an analogous purpose for the spanning
 numbers.
\end{remark}

\section{Outline of the constructions for Theorems~\ref{thm:main} and~\ref{thm:main-seponly}}\label{sec:overview}

The construction proceeds as follows:

\begin{itemize}
\item In Section~\ref{sec:Xn-construct} we construct a sequence of dynamical systems $(X_n,  F_n)$ that we call \emph{\JXn-systems} and  metrics~$D^n$ on~$X_n$ that
we call~\emph{\JXn-metrics.}
The sequence of systems~$(X_n, F_n)$ is identical for the proofs of
Theorems~\ref{thm:main} and~\ref{thm:main-seponly}, but the metrics~$D^n$ will be chosen slightly differently.
\item Our constructions rely on a number of parameters. In Section~\ref{sec:earlysetup} we describe these
parameters, list their required properties, and prove the existence of parameters with these properties.
In particular, the parameters include sequences of times $T(n)$ and $T^+(n)$
with $T(n) < T^+(n) < T(n+1) $
for all~$n$ and the diameters $\eps_n$ of the spaces $(X_n, D^n)$.  These parameters will be identical for
the proofs of Theorems~\ref{thm:main} and~\ref{thm:main-seponly}. They will be described in
Subsections~\ref{subsec:T(n)} and~\ref{subsec:epsn-deltan}. In Subsection~\ref{subsec:choose-colors} we
prove the existence  of certain
colorings~$c_n$ of finite sets of pairs of functions. These colorings will determine when $D^n(x_n, x_n')$
can attain the maximum value~$\eps_n$ for $x_n, x_n' \in X_n$, and will have slightly different
properties for the proofs of Theorems~\ref{thm:main} and~\ref{thm:main-seponly}.
\item In Section~\ref{sec:X} we construct what we call \emph{\JX-systems} $(X, F)$ as products of \JXn-systems $(X_n, F_n)$.  In particular,
we let the \emph{\JX-space} $X$ consist of all
sequences $x = (x_n)_{n \in \bN}$ such that $x_n \in X_n$ for each~$n \in \bN$. Metrics~$D$  on~$X$ that we call \emph{\JX-metric} are
defined by
$D(x, x') = \sum_{n \in \bN} D^n(x_n, x'_n)$, where the \JXn-metrics~$D^n$ in the terms of this sum are slightly different in the proofs of Theorems~\ref{thm:main} and~\ref{thm:main-seponly}.  Since $\eps_n$ will always be the diameter of~$(X^n, D^n)$ in both constructions,
 for $\eps = \sum_{n \in \bN} \eps_n$ and any~$T \in \bN$ we will have
\begin{equation*}
D_T(x, x') = \eps \ \Leftrightarrow \ \exists 0 \leq t < T \, \forall n \in \bN \ D^n(F_n^t(x_n), F_n^t(x'_n)) = \eps_n.
\end{equation*}
\item In Section~\ref{sec:X-} we first choose certain subsets~$\Ymn$ of the sets of
functions ${}^{T^+(n)}\{0,1\}$ that will be used in our constructions.
we then derive lower bounds on the sizes  $|\Ymn|$ of these sets (Claim~\ref{size}) and also a lower
bound on  the sizes $|W^n|$  of related sets~$W^n\subset X$ (Corollary~\ref{corol:size}). In Subsection~\ref{subsec:WX-}, we define the subspace~$W$ of~$X$ that will be
used in Theorem~\ref{thm:main-seponly} and we let the set~$X^-$ of Theorem~\ref{thm:main} be
the closure in~$X$ of the union of all sets $F^t(W^n)$ for $t\in\bZ$ and $n \in \bN$.
The subspaces~$W, X^-$ of~$X$ are compact and both backward and forward invariant under~$F$.
We also prove parts~(i) of Theorems~\ref{thm:main} and~\ref{thm:main-seponly} in this subsection.
In Subsection~\ref{subsec:propsWX-} we prove part~(iii) of Theorem~\ref{thm:main} and parts~(iii), (iv) of Theorem~\ref{thm:main-seponly}.
In Subsection~\ref{subsec:MoreTrans} we prove  part~(v) of Theorem~\ref{thm:main-seponly}.
\item  In Section~\ref{sec:sepsizes} we derive bounds on the separation and spanning numbers in our systems with respect to $D_T$ for certain choices of~$T$.
More specifically, in Subsection~\ref{subsec:lower-sep-T+n} we derive lower bounds on $sep(W, \eps, D_{T^+(n)})$ and $sep(X^-, \eps, D_{T^+(n)})$. This subsection
is relevant for the proofs of both Theorems~\ref{thm:main} and~\ref{thm:main-seponly}. Subsection~\ref{subsec:lower-span-T+n} is part of the proof Theorems~\ref{thm:main} only. Here we derive
lower bounds for the spanning numbers $span\left(X^-, \eps, D_{T^+(n)}\right)$ (Corollary~\ref{corol:low-span}).
In Subsection~\ref{subsec:upper-sep-2Tn} we derive upper bounds on $sep (X_n, \eps_n, D^n_{2T(n)})$
(Lemma~\ref{corol:low-separation-n}). This part
of the argument relies on properties of the colorings~$c_n$ that are used in the definition of the
metrics~$D^n$, but is common to the proofs of Theorems~\ref{thm:main} and~\ref{thm:main-seponly}.  The same
upper bounds remain valid for  $sep\left(Y, \eps, D_{2T(n)}\right)$ and $span\left(Y, \eps, D_{2T(n)}\right)$,
where $Y \in \{X^-, W, X\}$.
\item In Section~\ref{sec:prove-main} we wrap up the argument by comparing the lower and upper bounds derived in Section~\ref{sec:sepsizes} and show that the strict inequalities in Theorems~\ref{thm:main}(ii) and~\ref{thm:main-seponly}(ii) hold.
\end{itemize}

\section{Choosing suitable parameters}\label{sec:earlysetup}

Our constructions rely on certain sequences of mathematical objects that will be used as parameters.  Here we describe
these parameters, list their required properties, and prove the existence of parameters with these properties.
Throughout Sections~\ref{sec:Xn-construct}--\ref{sec:prove-spanonly} the standing assumption will be that the
parameters of the construction have the properties listed in the current section.

\subsection{Choosing $T(n)$ and $T^+(n)$}\label{subsec:T(n)}

As a first step, we fix  two sequences of positive integers $(T(n))_{n \in \bN}$ and $(T^+(n))_{n \in \bN}$ with
    \begin{equation*}\label{eqn:TT+}
    1 < T(0) < T^+(0)  < \dots < T(n) < T^+(n) < T(n+1) < T^+(n+1) < \dots
    \end{equation*}

These sequences will be defined in terms of two auxiliary sequences  $(C(n))_{n \in \bN}$ and $(K(n))_{n \in \bN}$ of positive integers so that  for all $n \in \bN$:
\begin{equation}\label{eqn:TT+(n)}
    T^+(n) = C(n)T(n) \qquad \mbox {and }\ \ \ \ \ \ \  T(n) = K(n)T^+(n-1).
    \end{equation}

The second part of~\eqref{eqn:TT+(n)} makes sense for $n = 0$ if we adopt the convention that
\begin{equation}\label{T+(-1)}
T^+(-1) =1.
\end{equation}
Then~\eqref{eqn:TT+(n)} implies that for all $n \in \bN$:
\begin{equation}\label{eqn:TT+(n)-explicit}
\begin{split}
    T(n) &= \prod_{i=0}^{n-1}C(i)\prod_{i=0}^n K(i)\\
    T^+(n) &= \prod_{i=0}^{n}C(i)K(i).
\end{split}
\end{equation}

We will choose these sequences so that for all~$n \in \bN$:
\begin{equation*}
            \begin{split}
                &\mbox{(PCn):}\ \ \ \prod_{i=0}^n (C(i) - 2) > 0.95 \prod_{i=0}^n C(i).\\
                &\mbox{(PKn1):}\ \ \ K(n)\ \mbox{is a positive integer multiple of}\ 100.\\
                &\mbox{(PKn2):}\ \ \ 2^{0.05T^+(n)} > \left(\frac{C(n)^2}{2}\right)\prod_{m=0}^{n-1}\left(\frac{C(m)^2}{2}\right)^{\prod_{i=m+1}^{n}(C(i)-2)K(i)}.\\
                &\mbox{(PKn3):}\ \ \ \left(2^{0.7T^+(n-1)K(n)}\right)! > \frac{C(n)^2}{2}.\\
                &\mbox{(PKn4):}\ \ \ \left(\log_2 \sqrt{\frac{3}{2}}\right)\left(2^{0.7T^+(n-1)K(n)}-1\right) > T^+(n-1) C(n)K(n).\\
                &\mbox{(PKn5):}\ \ \ 2^{0.01T(n)} = 2^{0.01K(n)T^+(n-1)} \geq C(n).
            \end{split}
        \end{equation*}

Note that (PKn1) implies that the exponents in (PKn2) through (PKn5) are integers, and that in
view of~\eqref{eqn:TT+(n)} we can write (PKn3) and (PKn4) equivalently as

\begin{itemize}
\item[(pKn3):] $\left(2^{0.7T(n)}\right)! > \frac{C(n)^2}{2}$.
\item[(pKn4):] $\left(\log_2 \sqrt{\frac{3}{2}}\right)\left(2^{0.7T(n)}-1\right) > C(n)T(n)$.
\end{itemize}

\begin{proposition}\label{prop:CnKn-can-be-chosen}
It is possible to choose sequences $(C(n))_{n \in \bN}$ and $(K(n))_{n \in \bN}$ so that properties (PCn) and (PKn1)--(PKn5) are satisfied for all~$n \in \bN$.
\end{proposition}

\noindent
\textbf{Proof:} For the proof of this proposition, it will be convenient to write (PCn) in the following equivalent form:
\begin{equation*}
            \mbox{(pcn):}\ \ \ \prod_{i=0}^n\left[1-\frac{2}{C(i)}\right] > 0.95
\end{equation*}

We prove the proposition by a recursive construction.

To get started, we first choose a positive integer $C(0)$ such that:

\medskip

$C(0)$ is large enough such that $1 - \frac{2}{C(0)} > 0.95$, which means (pc0) holds.

\medskip

Then we choose $K(0) = T(0)$ large enough such that:
        \begin{equation*}
            \begin{split}
                \mbox{(PK01): } &K(0) = T(0)\ \mbox{is a positive integer multiple of}\ 100,\\
                \mbox{(PK02):\ }&2^{0.05C(0)K(0)} = 2^{0.05C(0)T(0)} > \frac{C(0)^2}{2},\\
                \mbox{(PK03):\ }&\left(2^{0.7K(0)}\right)! = \left(2^{0.7T(0)}\right)! > \frac{C(0)^2}{2},\\
                \mbox{(PK04):\ }&\left(\log_2{\sqrt{\frac{3}{2}}}\right)\left(2^{0.7T(0)}-1\right) = \left(\log_2{\sqrt{\frac{3}{2}}}\right)\left(2^{0.7K(0)}-1\right)\\
                &> C(0)T(0) = C(0)K(0),\\
                \mbox{(PK05): } &2^{0.01K(0)} = 2^{0.01T(0)} \geq C(0).
            \end{split}
        \end{equation*}

Note that we use here the first line of~\eqref{eqn:TT+(n)-explicit} and the fact that the products \newline
$\prod_{i=0}^{n-1}C(i)$ in~\eqref{eqn:TT+(n)-explicit} and $\prod_{m=0}^{n-1}\left(\frac{C(m)^2}{2}\right)^{\prod_{i=m+1}^{n}(C(i)-2)K(i)}$ in (PKn2)

have no terms and thus are equal to~1 for $n = 0$.

\medskip

Now assume $n > 0$ and $C(m), K(m)$ (and hence $T(m), T^+(m)$) have been chosen for all $m <n$ so that the conditions (pcm), (PKm1)--(PKm5) are satisfied for all
$m < n$.

Next we choose a sufficiently large integer~$C(n)$ so that

\begin{equation*}
\mbox{(pcn):\ \ }\prod_{i=0}^{n}\left[1-\frac{2}{C(i)}\right] > 0.95.
\end{equation*}

This is possible, since by the inductive assumption we have
\begin{equation*}
\mbox{(pc(n-1)):\ \ }\prod_{i=0}^{n-1}\left[1-\frac{2}{C(i)}\right] > 0.95.
\end{equation*}

Now we need to choose $K(n)$ such that (PKn1)--(PKn5) will hold.

Clearly, (PKn1) can be easily satisfied, and (PKn3) and (PKn5) will hold for all sufficiently large~$K(n)$ since the right-hand
side of this inequality is already fixed by the choices we have made up to this point.  Similarly, (PKn4) will hold for
all sufficiently large~$K(n)$, since $K(n)$ enters the left-hand side of this inequality in the exponent, while it
enters the right-hand side as a multiplicative factor.

Condition (PKn2) is slightly more delicate.   By inductive assumption we have:
\begin{equation*}
\mbox{(PK(n-1)2):}\ \ \ 2^{0.05T^+(n-1)} > \left(\frac{C(n-1)^2}{2}\right)\prod_{m=0}^{n-2}\left(\frac{C(m)^2}{2}\right)^{\prod_{i=m+1}^{n-1}(C(i)-2)K(i)}.
\end{equation*}

We want to show
\begin{equation*}
\mbox{(PKn2):}\ \ \ 2^{0.05T^+(n)} > \left(\frac{C(n)^2}{2}\right)\prod_{m=0}^{n-1}\left(\frac{C(m)^2}{2}\right)^{\prod_{i=m+1}^{n}(C(i)-2)K(i)}.
\end{equation*}

Let $LHS(n-1), LHSn, RHS(n-1), RHSn$ denote the left-hand sides and righthand sides of (PK(n-1)2) and (PKn2), respectively. Then by~\eqref{eqn:TT+(n)} and the form of the left-hand and right-hand sides:
\begin{equation*}
\begin{split}
LHSn &= [LHS(n-1)]^{C(n)K(n)}\\
RHSn &= \frac{C(n)^2}{2} [RHS(n-1)]^{C(n)K(n)}[RHS(n-1)]^{-2K(n)}.
\end{split}
\end{equation*}
Since $C(n)$ has already been chosen, $RHS(n-1) \geq \frac{C(n-1)^2}{2} > 1$, and $LHS(n-1) > RHS(n-1)$, it follows that we can guarantee (PKn2) by choosing $K(n)$ sufficiently large. $\Box$

\bigskip

\subsection{Choosing $I^n_j, \eps_n, \delta_n, \lambda$}\label{subsec:epsn-deltan}
Once $T(n), C(n), T^+(n)$ are chosen, we partition the interval $[0, T^+(n)-1)$ into consecutive
subintervals $I^n_j$ of length~$T(n)$ each, where ~$j$ ranges from~$1$ to~$C(n)$.  More precisely,
we will treat each~$I^n_j$ as a sequence rather than a set of consecutive integers, so that
\begin{equation}\label{eqn:Injs}
\begin{split}
I^n_1 &= (0, \dots , T(n)-1),\\
I^n_2 &= (T(n), \dots , 2T(n)-1),\\
\dots &= \dots \\
I^n_j   &= ((j-1)T(n), \dots , jT(n)-1),\\
\dots &= \dots \\
I^n_{C(n)} &= ((C(n)-1)T(n), \dots , C(n)T(n)-1).
\end{split}
\end{equation}

Next we fix  $\eps > 0$ and sequences $(\eps_n)_{n \in \bN}$ and $(\delta_n)_{n \in \bN}$ such that
    \begin{itemize}
          \item[(P$\eps$)] $(\eps_n)_{n \in \bN}$ is a strictly decreasing sequence of positive real numbers and
          \newline $\eps := \sum_{n \in \bN} \eps_n < \infty$.
    \end{itemize}

          \noindent{Moreover, for each~$n \in \bN$:}
    \begin{itemize}
          \item[(P$\delta$1)] $0 < \delta_n < \eps_n <2\delta_n$.
          \item[(P$\delta$2)] $\sum_{n^+ > n} \eps_{n^+} < 0.5(\eps_{n} - \delta_n)$.
          \item[(P$\delta$3)] $(\eps_{n} - \delta_n) < \eps_n 3^{-2\lambda T^+(n)}$.
    \end{itemize}

Here $\lambda$ is a parameter that represents a positive integer.  Throughout this preprint we will set
\begin{equation}\label{eqn:r=1}
\lambda = 1.
\end{equation}

This makes~$\lambda$ redundant for the current purpose; however, in some of constructions in~\cite{PartII} we may
need the added flexibility offered by other choices for~$\lambda$. In order to avoid duplication of effort, we
include this parameter here and show explicitly that certain essential properties of our constructions do not depend
on the particular choice~\eqref{eqn:r=1} of its value.

Satisfaction of the above conditions can be assured by first choosing $\eps_0$, next  $\delta_0$ such that  (P$\delta$1) and (P$\delta$3) hold.
At step $n$, when $\eps_m, \delta_m$ have already been chosen for $m < n$, first choose $\eps_n$ small
enough such that
\begin{equation*}
\mbox{(EDn): \ \ \ }  \forall m < n \ \eps_{m+1} + \eps_{m+2} + \dots + \eps_n < 0.4(\eps_m - \delta_m).
\end{equation*}
Since (ED0) is vacuously true, this will be possible under the inductive assumption that (ED(n-1)) holds.   Next
choose $\delta_n$ so that (P$\delta$1) and (P$\delta$3) hold, and so on.  Then, (EDn) will hold for all~$n$
so that for any fixed $m \in \bN$,
\begin{equation*}
    \forall n > m\ \ \eps_{m+1} + \eps_{m+2} + \dots + \eps_n < 0.4(\eps_m - \delta_m),
\end{equation*}
and we get
\begin{equation*}
    \sum_{m^+ > m}\eps_{m^+} \leq 0.4(\eps_m - \delta_m) < 0.5(\eps_m - \delta_m).
\end{equation*}

\subsection{Choosing colorings $c_n$}\label{subsec:choose-colors}

Let ${}^{T^+(n)}\{0,1\}$ denote the set of all functions with domain $\{0, 1, \dots , T^+(n)-1\}$ that take values in the set $\{0,1\}$. For a subset
$S \subseteq {}^{T^+(n)}\{0,1\}$ let $[S]^2$ denote the set of all unordered pairs $\{\varphi, \psi\}$ of different functions from~$S$.
Moreover, let $[C(n)] = \{1, 2, \ldots , C(n)\}$.

For the purpose of our arguments, a \emph{coloring} will be a function \newline
$c_n: \left[{}^{T^+(n)}\{0,1\}\right]^2 \rightarrow [C(n)]$ for some~$n$.

For the proof of Theorem~\ref{thm:main-seponly},
for each $n \in \bN$ we choose a coloring
$c_n$ such that for all $n \in \bN$:

\begin{itemize}
\item[(cC1)] Assume $\varphi\restrict (0, \dots, T^+(n-1)-1) \neq \psi\restrict (0, \dots, T^+(n-1)-1)$ while $\varphi(i) = \psi(i)$ for all $i \in \{T^+(n-1), T^+(n-1)+1, \dots, T^+(n)-1\}$. \newline
     Then, $c_n(\varphi, \psi) = 1$.
\item[(cC2)] Let $C(\varphi, \psi) \subset [C(n)]$ denote the set of $j$ such that $\varphi \restrict I^n_j \neq \psi \restrict I^n_j$. \newline If $|C(\varphi, \psi)| \geq 3$, then $c_n(\varphi, \psi) \in C(\varphi, \psi)$.
\item[(cC3)] For every subset $S \subset {}^{T^+(n)}\{0,1\}$ of size~$|S| \geq  2^{0.75T(n)}$ the restriction of~$c_n$ to $[S]^2$ takes on at least three colors.
\end{itemize}

Similarly, for the proof of Theorem~\ref{thm:main},
for each $n \in \bN$ we choose a coloring
$c_n$
such that for all $n \in \bN$:

\begin{itemize}
\item[(cCi)] Assume $\varphi(i) = 0$ for all $T^+(n-1) \leq i \leq T^+(n)-1$. \newline
Then, $c_n(\varphi, \psi) = 1$ for all $\psi \neq \varphi \in {}^{T^+(n)}\{0,1\}$.
\item[(cC)] Assume there exist $T^+(n-1) \leq i, j \leq T^+(n)-1$ such that $\varphi(i) = \psi(j) = 1$.
 Let $C(\varphi, \psi) \subset [C(n)]$ denote the set of $j$ such that $\varphi \restrict I^n_j \neq \psi \restrict I^n_j$. \newline
If $|C(\varphi, \psi)| \geq 3$, then $c_n(\varphi, \psi) \in C(\varphi, \psi)$.
\item[(cC3)] For every subset $S \subset {}^{T^+(n)}\{0,1\}$ of size~$|S| \geq  2^{0.75T(n)}$ the
restriction of~$c_n$ to $[S]^2$ takes on at least three colors.
\end{itemize}

In our proofs, conditions~(cC1) and~(cCi) will play similar roles.  Note that conditions (cC2) and (cC) can be
understood as saying ``If the value of~$c_n(\varphi, \psi)$ is not already determined by the first condition,
that is, by~(cC1) or by~(cCi), and if  $|C(\varphi, \psi)| \geq 3$,
then $c_n(\varphi, \psi) \in C(\varphi, \psi)$.''

\begin{lemma}\label{lem:coloring-exists}
For every~$n \in \bN$ there exist a coloring $c_n$ that satisfies (cC1)--(cC3) and a coloring $c_n'$ that satisfies (cCi)--(cC3).
\end{lemma}

\noindent
\textbf{Proof:} Fix~$n \in \bN$. Consider $c_n$ that is randomly chosen from the uniform distribution on all functions
  $c: \left[{}^{T^+(n)}\{0,1\}\right]^2 \rightarrow [C(n)]$ that satisfy conditions~ (cC1) and (cC2), and
  $c_n'$ that  is randomly chosen from the uniform distribution of all such functions~$c$ that satisfy conditions~(cCi) and~(cC).

Note that we can obtain random objects from these distributions by randomly and independently assigning values~$c_n(\varphi, \psi)$ and $c_n'(\varphi, \psi)$ subject to the following conditions:

\begin{itemize}
\item[(C1)] If $\varphi\restrict (0, \dots, T^+(n-1)-1) \neq \psi\restrict (0, \dots, T^+(n-1)-1)$ while $\varphi(i) = \psi(i)$ for all $i \in \{T^+(n-1), T^+(n-1)+1, \dots, T^+(n)-1\}$, then we let $c_n(\varphi, \psi) = 1$.
\item[(Ci1)] If $\varphi \neq \psi$ and $\varphi(i) = 0$ for all $T^+(n-1) \leq i \leq T^+(n)-1$ or $\psi(i) = 0$ for all $T^+(n-1) \leq i \leq T^+(n)-1$, then we let $c_n'(\varphi, \psi) = 1$.
\item[(C2)] If $\varphi$ and $\psi$ differ on at most two intervals $I^n_j$, but they do not differ only on $(0, \dots, T^+(n-1)-1)$, then
we randomly choose $c_n(\varphi, \psi)$ from~$[C(n)]$ with the uniform distribution.
\item[(Ci2)] If $\varphi$ and $\psi$ differ on at most two intervals $I^n_j$, and there exist $T^+(n-1) \leq i , j \leq T^+(n)-1$ such that $\varphi(i) = 1$ and $\psi(j) = 1$, then
we randomly choose $c_n'(\varphi, \psi)$ from~$[C(n)]$ with the uniform distribution.

\item[(C3)] If $\varphi$ and $\psi$ differ on at least three intervals $I^n_j$, then we
randomly choose $c_n(\varphi, \psi)$ from~$C(\varphi, \psi)$ with the uniform distribution.

\item[(Ci3)] If $\varphi$ and $\psi$ differ on at least three intervals $I^n_j$, and there
exist $T^+(n-1) \leq i , j \leq T^+(n)-1$ such that $\varphi(i) = 1$ and $\psi(j) = 1$, then we randomly choose $c_n(\varphi, \psi)$ from~$C(\varphi, \psi)$ with the uniform distribution.
\end{itemize}

Note that if clause~(C2) or clause~(Ci2) is used for determining the value of $c_n(\varphi, \psi)$ or of
$c_n'(\varphi, \psi)$ , then
\begin{equation}\label{eqn:P-one-third-3}
\begin{split}
\forall i \in [C(n)] \ \ P(c_n(\varphi,\psi)=i) &= \frac{1}{C(n)} \ \leq \ \frac{1}{3},\\
\forall i \in [C(n)] \ \ P(c_n'(\varphi,\psi)=i) &= \frac{1}{C(n)} \ \leq \ \frac{1}{3}.
\end{split}
\end{equation}

More generally,  in all clauses where we have some choice, that is, where $c_n(\varphi, \psi)$ is not
already determined by~(cC1)
and $c_n'(\varphi, \psi)$ is not already determined by~(cCi),  we get
\begin{equation}\label{eqn:P-one-third}
\begin{split}
\forall i \in [C(n)] \ \ P(c_n(\varphi,\psi)=i) \  &\leq \ \frac{1}{3},\\
\forall i \in [C(n)] \ \ P(c_n'(\varphi,\psi)=i) \ &\leq \ \frac{1}{3}.
\end{split}
\end{equation}

Now consider any $S \subset {}^{T^+(n)}\{0,1\}$ and let $S^- \subset S$.

We call~$S^-$ a \emph{(cC1)-free set} if
for all $\varphi \neq \psi \in S^-$ we have
$\varphi\restrict (T^+(n-1), \dots, T^+(n)-1) \neq \psi\restrict (T^+(n-1), \dots, T^+(n)-1)$ so that
condition in (cC1) does not determine the value of  $c_n(\varphi,\psi)$.

Similarly,
we call~$S^-$ a \emph{(cCi)-free set} if
for all $\varphi \in S^-$, there exists $T^+(n-1)\leq i \leq T^+(n)-1$ such that $\varphi(i) =1$, so that
condition in (cCi) does not determine the value of  $c_n(\varphi,\psi)$ for any $\varphi \neq \psi \in S^-$.

\begin{proposition}\label{prop:cCi-}
Let $S \subset {}^{T^+(n)}\{0,1\}$ be such that $|S| \geq 2^{0.75T(n)}$. Then

\smallskip

\noindent
(a) There exists a (cC1)-free set $S^- \subset S$ with $|S^-| \geq 2^{0.7T(n)}$.

\smallskip

\noindent
(b) There exists a (cCi)-free set $S^- \subset S$ with $|S^-| \geq 2^{0.7T(n)}$.
\end{proposition}

\noindent
\textbf{Proof:}  Let $S$ be as in the assumption. We prove the slightly easier part~(b) first.
For this part, we get $S^- \subset S$ by
 removing all $\varphi \in S$ with $\varphi(i) = 0$ for all $T^+(n-1) \leq i \leq T^+(n)-1$.   Then,
\begin{equation*}
    \begin{split}
    |S^-| &\geq 2^{0.75T(n)} - 2^{T^+(n-1)}\\
    &= 2^{0.7T(n)}\left(2^{0.05T(n)}-2^{\left(\frac{1}{K(n)}-0.7\right)T(n)} \right)\\
    &\geq 2^{0.7T(n)}\left(2^{0.05T(n)} - 1\right)\\
    &> 2^{0.7T(n)}.
    \end{split}
\end{equation*}

For the proof of part~(a) arrange the elements of $S = \{\varphi_i\}_{i=1}^{|S|} = L(1)$ into a list.
We recursively construct lists $L(i)$ for $i = 2, \dots |S|+1$ by (possibly) removing some elements of $L(i)$ to obtain $L(i+1)$ as follows:
\begin{itemize}
\item If $\varphi_i \notin L(i)$, then we let $L(i+1) = L(i)$.
\item If $\varphi_i \in L(i)$, then we obtain $L(i+1)$ by  removing all $\varphi_j \in L(i)$ with $j > i$ from $L(i)$ that differ from $\varphi_i$ only on $(0, \dots, T^+(n-1)-1)$.
\end{itemize}

Let $S^-$ be the set of $\varphi_i$ that survived this procedure, that is, let $S^- = L(|S|+1)$. Since removal
is always conditioned on a prior decision to retain some $\varphi_i \in L(i)$, there are at
most $|S^-|$ steps $i$ where any removal took place, that is, where $L(i+1) \neq L(i)$.
Moreover, at each such step we removed
 at most  $(2^{T^+(n-1)}-1)$ elements.    Thus,
\begin{equation*}
    \begin{split}
        |S^-| & \geq |S| - (2^{T^+(n-1)}-1)|S^-|\\
        &> |S| - 2^{T^+(n-1)}|S^-|\\
        &\geq 2^{0.75T(n)} - 2^{T^+(n-1)}|S^-|.
    \end{split}
\end{equation*}
Therefore
\begin{equation*}
    \begin{split}
        |S^-| &\geq \frac{2^{0.75T(n)}}{2^{T^+(n-1)}+1}\\
        &> 2^{0.75T(n)-2T^+(n-1)}\\
        &= 2^{(0.75-\frac{2}{K(n)})T(n)}\\
        &> 2^{0.7T(n)},
    \end{split}
\end{equation*}
where the last two lines follow from~\eqref{eqn:TT+(n)} and~(PKn1).
$\Box$

\bigskip

Now it remains to show the existence of a colorings $c \in \{c_n, c_n'\}$ that are randomly drawn as above and that take at
least three colors on~$[S^-]^2$ for all (cC$*$)-free sets~$S^-$  of size $2^{0.7T(n)}$, where ``$*$'' should be
replaced by ``$1$'' for the proof of existence of~$c = c_n$ and by ``i''
for the proof of existence of $c = c_n'$.

We say a subset $A \subset {}^{T^+(n)}\{0,1\}$ is $\leq 2$-chromatic for $c$ if $c$ takes on at most two colors on $[A]^2$.

Let $F_{\varphi,\psi,i}$ denote the event that $c(\varphi,\psi) = i$. It follows from the above description that for fixed~$i$ all these events are independent. Moreover,
by~\eqref{eqn:P-one-third}, when $\{\varphi, \psi\} \in [S^-]^2$ for some (cC$*$)-free set~$S^-$,
then for all~$i, j$:
\begin{equation}\label{eqn:PFi}
P\left(F_{\varphi,\psi,i}\right) \leq \frac{1}{3} \qquad \mbox{and} \qquad P\left(F_{\varphi,\psi,i} \ \cup F_{\varphi,\psi,j} \right) \leq \frac{2}{3}.
\end{equation}

For each (cC$*$)-free set $M \subset {}^{T^+(n)}\{0,1\}$ of size exactly $2^{0.7T(n)}$ let
 $\xi_M$ be the r.v. that, for a randomly drawn $c$, takes the value 1 if $M$ is $\leq 2$-chromatic and takes the value 0 otherwise.  Since there are $C(n) \choose 2$ possible subsets~$\{i,j\}$ of~$[C(n)]$,  it follows from~\eqref{eqn:PFi} and independence that the expected value of~$\xi_M$ satisfies
\begin{equation}\label{eqn:ExiM}
    E(\xi_M) = P(\xi_M = 1) \leq {C(n) \choose 2}\left(\frac{2}{3}\right)^{\frac{2^{0.7T(n)}\left(2^{0.7T(n)}-1\right)}{2}}.
\end{equation}

Let
\begin{equation*}
\cM = \{M \subset {}^{T^+(n)}\{0,1\}:\  |M| = 2^{0.7T(n)}\quad \mbox{and} \ M \ \mbox{is a (cC$*$)-free set}\}.
\end{equation*}

Then
\begin{equation}\label{eqn:size-cM}
|\cM| \leq  {2^{T^+(n)} \choose 2^{0.7T(n)}} < \frac{\left(2^{C(n)T(n)}\right)^{2^{0.7T(n)}}}{\left(2^{0.7T(n)}\right)!}.
\end{equation}

Let $\xi$ be the r.v. that counts the number of $\leq 2$-chromatic (cC$*$)-free subsets of size~$2^{0.7T(n)}$.  That is, let
\begin{equation*}
    \xi = \sum_{M \in \cM} \xi_M.
\end{equation*}
Then we get the following estimate:
\begin{equation}\label{eqn:Exi}
    \begin{split}
        E(\xi) &= \sum_{M \in \cM} E(\xi_M)\\
        &\leq {C(n) \choose 2}\left(\frac{2}{3}\right)^{\left(\frac{2^{0.7T(n)}(2^{0.7T(n)}-1)}{2}\right)}|\cM|\\
        &< \frac{C(n)(C(n)-1)}{2}\left(\frac{2}{3}\right)^{\left(\frac{2^{0.7T(n)}(2^{0.7T(n)}-1)}{2}\right)}\left[\frac{\left(2^{C(n)T(n)}\right)^{2^{0.7T(n)}}}{\left(2^{0.7T(n)}\right)!}\right]\\
        &< \left(\frac{2}{3}\right)^{\left(\frac{2^{0.7T(n)}(2^{0.7T(n)}-1)}{2}\right)}\left(2^{C(n)T(n)}\right)^{2^{0.7T(n)}}\\
        &< \left[\left(\frac{2}{3}\right)^{\left(\frac{2^{0.7T(n)}-1}{2}\right)} 2^{\left[\left(\log_2 \sqrt{\frac{3}{2}}\right)\left(2^{0.7T(n)}-1\right)\right]}\right]^{2^{0.7T(n)}}\\
        &= \left[\left(\frac{2}{3}\right)^{\left(\frac{2^{0.7T(n)}-1}{2}\right)}\left(\frac{3}{2}\right)^{\left(\frac{2^{0.7T(n)}-1}{2}\right)}\right]^{2^{0.7T(n)}}\\
        &= 1.
    \end{split}
\end{equation}

The first inequality in~\eqref{eqn:Exi} follows from~\eqref{eqn:ExiM}, the second follows from~\eqref{eqn:size-cM}, the third follows from~(pKn3),
and the fourth follows from~(pKn4).

Since the expected value of $\xi$ is less than 1, we must have $\xi(c) = 0$ for at least one coloring $c$, which witnesses the result claimed here.
$\Box$

\section{Construction of \JXn-systems~$(X_n, F^n)$ and \JXn-metrics  $D_n$}\label{sec:Xn-construct}

For each~$n \in \bN$ we define the \emph{\JXn-system} $(X_n, F_n)$ as follows:
\begin{itemize}
\item The set $X_n$ consists of all triples $(y, n, k)$, where $y \in {}^\bZ \{0, 1\}$ and $k \in \{0, 1, \dots , T^+(n)-1\}$.

\item The function~$F_n$ is defined by
\begin{equation}\label{eqn:def-Fn}
\begin{split}
F_n((y, n, k)) &= (\sigma(y), n, F_n(k)), \ \mbox{where}\\
\sigma(y)(i) &= y(i+1) \ \mbox{ for all } i,\\
F_n(k) &= (k+1) \ mod \ T^+(n).
\end{split}
\end{equation}
\end{itemize}

Thus the \JXn-system $(X_n, F_n)$ is uniquely determined by~$n$ and is essentially the product of the full subshift system~$({}^\bZ \{0,1\}, \sigma)$ with a cyclic permutation of $T^+(n)$. The acronym ``\JX'' can be read, for example, as ``$\eps$-coloring'' or ``$\eps$-coding.''
The symbol~$F_n$ does double duty here and denotes both the cyclic permutation and the product with the subshift operator. This
will be useful in Section~\ref{sec:prove-spanonly} and  should not lead to confusion, as we will always specify the argument
of~$F_n$.

Let $X_n^k = \{(y, n, k') \in X_n: \ k' = k\}$.

The sets $X_n^k$ are pairwise disjoint and $X_n = \bigcup_{0\leq k < T^+(n)} X^k_n$.

\smallskip

Next we define an auxiliary function $\Phi: \bigcup_{n \in\bN} X_n \rightarrow {}^{T^+(n)}\{0,1\}$:
\begin{equation}\label{eqn:Phi-def}
\Phi((y, n, k)) = (y(-k), y(-k+1), \dots , y(-k+T^+(n)-1)) \in {}^{T^+(n)}\{0,1\}.
\end{equation}

Let us make a crucial observation that follows immediately from the definitions
of the functions $\Phi((y, n, k))$ and $F_n$:

\begin{proposition}\label{prop:Phi-const}
Let $0 \leq k < T^+(n) - 1$.  Then $\Phi((y, n, k)) = \Phi(F_n((y, n, k)))$.

On the other hand, for $k = T^+(n)-1$ we may have $\Phi((y, n, k)) \neq \Phi(F_n((y, n, k)))$.
\end{proposition}

\noindent
\textbf{Proof:}  If $0 \leq k < T^+(n)-1$, then $1 \leq k+1 < T^+(n)$ and\\ $(k+1)\ mod\ T^+(n)  = k+1$.  Thus,
\begin{equation*}
\begin{split}
&\Phi(F_n((y,n,k))) = \Phi((\sigma(y),n,(k+1)\ mod\ T^+(n)))\\
&= \Phi((\sigma(y),n,k+1))\\
&= (\sigma(y)(-k-1), \sigma(y)(-k), \dots, \sigma(y)(-k-1+T^+(n)-1))\\
&= (y(-k), y(-k+1), \dots, y(-k+T^+(n)-1))\\
&= \Phi((y,n,k)).
\end{split}
\end{equation*}
On the other hand, if $k = T^+(n)-1$, then $k+1 = T^+(n)$ and \\
$(k+1)\ mod\ T^+(n) = 0$.  Thus,
\begin{equation*}
    \begin{split}
    \Phi(F_n((y,n,k))) &= \Phi((\sigma(y),n,(k+1)\ mod\ T^+(n)))\\
    &= \Phi((\sigma(y),n,0))\\
    &= (\sigma(y)(0), \sigma(y)(1), \dots, \sigma(y)(T^+(n)-1))\\
    & = (y(1), y(2), \dots, y(T^+(n))),
    \end{split}
\end{equation*}
and
\begin{equation*}
\begin{split}
    \Phi((y,n,k)) &= (y(-k), y(-k+1), \dots, y(-k+T^+(n)-1))\\
    &= (y(-T^+(n)+1), \dots, y(0)).
\end{split}
\end{equation*}
Therefore, in this case, we may have $\Phi((y,n,k)) \neq \Phi(F_n((y,n,k)))$.  $\Box$

\bigskip

We need two more auxiliary functions: the bijection $\#: \bZ \rightarrow \bN$ given by
\begin{equation}\label{eqn:numberi}
\#(0) = 0, \quad  \#(1) = 1, \quad \#(-1) = 2, \quad \#(2) = 3, \quad \#(-2) = 4, \dots
\end{equation}
and the  function $\Delta: \left({}^\bZ \{0, 1\}\right)^2 \rightarrow \bN \cup \{\infty\}$ that takes the value
$D(y,z) = \infty$ when $y = z$ and the value $\Delta(y,z) = \#(i)$ when
 $y \neq z$, where~$i$ is such that
\begin{equation}\label{deqn:define-Delta}
y(i) \neq z(i) \quad \mbox{and} \quad \forall j \in \bZ \ (\#(j) < \#(i) \ \Rightarrow \ y(j) = z(j)).
\end{equation}

Let $\beta_n \in \{\eps_n , \delta_n\}$. We define a \emph{\JXn-metric} as any function $D^n: (X_n)^2 \rightarrow [0, \infty)$ that satisfies the following conditions:

\smallskip

\begin{itemize}
\item[(Dn1)] If $k \neq k'$, then $D^n((y, n, k), (y', n, k')) = \beta_n$.
\item[(Dn2)] If $k = k'$ and $y = y'$, then $D^n((y, n, k), (y', n, k')) = 0$.
\item[(Dn3)] If $k = k'$ and $y \neq y'$, then
      \begin{itemize}
      \item[(Dn31)] If $0 < \Delta(y,y') < \infty$, then $D^n ((y, n, k), (y', n, k')) = \eps_n 3^{-\Delta(y, y')}$.
      \item[(Dn32)] If $\Delta(y,y') = 0$, then the value\\
      $D^n ((y, n, k), (y', n, k')) = D^n ((y', n, k'), (y, n, k))\in \{\eps_n, \delta_n\}$ \\
      and may depend only on~$k = k'$,\\
      $y\restrict (-\lambda T^+(n)+1, \dots, 0, \dots, \lambda T^+(n)-1)$, and \\
      $y'\restrict (-\lambda T^+(n)+1, \dots, 0, \dots, \lambda T^+(n)-1)$.
\end{itemize}
\end{itemize}

\smallskip

Recall from~\eqref{eqn:r=1} that in the constructions presented here we will always set $\lambda =1$. However, the general results \JXn- and \JX-systems with
\JXn- and \JX-metrics that we will derive in this preprint apply when the parameter~$\lambda$ is chosen as any positive integer. This may be useful
for our work in~\cite{PartII}.

For the proof of Theorem~\ref{thm:main-seponly} we will choose $\beta_n = \delta_n$ for all~$n \in \bN$ and work with \JXn-metrics
that satisfy:

\smallskip

\begin{itemize}
\item[(Dn1d)] If $k \neq k'$, then $D^n((y, n, k), (y', n, k')) = \delta_n$.
\end{itemize}

In contrast, for the proofs of Theorems~\ref{thm:main}
and~\ref{thm:main-spanonly} we will choose~$\beta_n = \eps_n$ for all~$n \in \bN$ and work with \JXn-metrics
that satisfy:

\smallskip

\begin{itemize}
\item[(Dn1e)] If $k \neq k'$, then $D^n((y, n, k), (y', n, k')) = \eps_n$.
\end{itemize}

In the proofs of both Theorems~\ref{thm:main} and~\ref{thm:main-seponly},  clause~(Dn32) will take the following form for some colorings~$c_n$:

\smallskip

\begin{itemize}
  \item[(Dn32c)] If $\Delta(y,y') = 0$, then we let $\varphi = \Phi((y,n,k))$ and $\psi = \Phi((y',n, k'))$
      and define:
            \begin{itemize}
            \item $D^n ((y, n, k), (y', n, k')) = \eps_n$ if $k  \in I^n_j$ and $c_n(\varphi, \psi) = j$.
            \item $D^n ((y, n, k), (y', n, k')) = \delta_n$ if $k \in I^n_j$ and $c_n(\varphi, \psi) \neq j$.
            \end{itemize}
     \end{itemize}

In the proof of Theorem~\ref{thm:main} we will use
colorings~$c_n$ that satisfy conditions~(cCi), (cC), and~(cC3), while in the proof of Theorem~\ref{thm:main-seponly},
we will use
colorings~$c_n$ that satisfy conditions~(cC1)--(cC3).

\begin{proposition}\label{prop:metric-n}
Let $D^n$ be a \JXn-metric.   Then

\smallskip

\begin{itemize}

\item[(i)] The function $D^n$ is a metric on~$X_n$.

\smallskip

\item[(ii)] The systems~$(X_n, D^n, F_n)$ have the following properties:

\smallskip

            \begin{itemize}
            \item[(PDn1)] $\max\{ D^n(z, z'): \, z, z' \in X_n\} \leq \eps_n$.

            Moreover, when the definition of~$D^n$ includes clause~(Dn1e) or \\
            clause~(Dn32c),
            then $\max\{ D^n(z, z'): \, z, z' \in X_n\} = \eps_n$.

            \smallskip

            \item[(PDn2)] $D^n(z, z') < \eps_n \ \Rightarrow \ D^n(z, z') \leq \delta_n$.

            \smallskip

            \item[(Pnc)] Each $D^n$ is a metric on~$X_n$ that induces the topology of a compact Hausdorff space.

            \smallskip

            \item[(PFn)] Each $F_n : X_n \rightarrow X_n$ is a homeomorphism.
            \end{itemize}

\end{itemize}
\end{proposition}

\noindent
\textbf{Proof:}
For easier flow of the argument, we prove (PDn1) and (PDn2) in (ii) first.

\smallskip

\noindent
(ii)-(PDn1):  Let $z=(y,n,k)$ and $z'=(y',n,k') \in X_n$. Note that the definition of~$D^n$ allows only for values $D^n(z,z') \in \{\beta_n, 0, \eps_n, \delta_n, \eps_n 3^{-\Delta(y,y')}\}$, where $\beta_n \in \{\delta_n, \eps_n\}$.  Thus,  $\max\{ D^n(z, z'): \, z, z' \in X_n\} \leq \eps_n$.

\medskip

Now we show that the supposed maximum value~$\eps_n$ of~$D^n$ is actually attained under the specified conditions.

When (Dn1e) applies, as in the proofs of Theorems~\ref{thm:main} and~\ref{thm:main-spanonly}, let
$z = (y,n,k)$ and $z' = (y',n,k')$ be elements of $X_n$ with $k \neq k'$. Then $D^n(z,z') = \eps_n$
 according to this clause.

When the definition of $D^n$ is based on clause~(Dn32c), as in the proof of Theorem~\ref{thm:main-seponly},  let $y = {}^\bZ\{0\}$, $y' = {}^\bZ\{1\}$,
and let $\varphi = {}^{T^+(n)}\{0\}$ and $\psi = {}^{T^+(n)}\{1\}$.  Suppose $c_n(\varphi, \psi) = c \in [C(n)]$.
Take any $k \in I^n_c$. Let $z = (y,n,k)$ and let $z' = (y',n,k)$.  They are both elements of $X_n$.
Since $y(0) = 0 \neq 1 = y'(0)$, clause~(Dn32c) is used to  determine the value of~$D^n(z,z')$,
which must be equal to~$\eps_n$ for these choices according to this clause.

\smallskip

\noindent
(ii)-(PDn2):  If $D^n(z,z') < \eps_n$ is defined by clause (Dn1), then  $D^n(z,z') = \beta_n = \delta_n$.

If $D^n(z,z') < \eps_n$ is defined by clause (Dn2), then  $D^n(z,z') = 0 \leq \delta_n$.

If clause~(Dn31) applies, then $D^n(z,z') = \eps_n 3^{-\Delta(y,y')} \leq \frac{1}{3}\eps_n \leq \delta_n$ by (P$\delta$1).

Finally, clause~(Dn32) allows only $D^n(z,z') \in \{\eps_n, \delta_n\}$.

Thus, $D^n(z, z') < \eps_n \ \Rightarrow \ D^n(z, z') \leq \delta_n$.

\bigskip

\noindent
(i) We need to verify the defining properties of a metric.

\begin{itemize}
\item Reflexivity: Let $z = (y,n,k) \in X_n$.  With $k = k$ and $y = y$, we have $D^n(z, z) = 0$ by (Dn2).
\item Positive definiteness: Let $z = (y,n,k)$ and  $z' = (y',n,k')$ in $X_n$ be such that $z \neq z'$.   If $k \neq k'$, then $D^n(z, z') = \beta_n > 0$ by (Dn1).  If $k = k'$, then $y \neq y'$.   In this case, by (Dn3) and (P$\delta$1), we have $D^n(z, z')  = \eps_n 3^{-\Delta(y,y')} >0$ if $\Delta(y,y') > 0$; and $D^n(z, z')  \geq \delta_n > 0$ if $\Delta(y,y') = 0$.
\item Symmetry: Let $z, z' \in X_n$. The equality $D^n(z, z') = D^n(z',z)$ follows directly by the definition of $D^n$ and the fact that
$\Delta(y, y') = \Delta(y', y)$ and symmetry is explicitly built into clause~(Dn32).
\item The Triangle Inequality: Let $z_1 = (y_1, n, k_1)$, $z_2 = (y_2, n, k_2)$, and \newline
$z_3 = (y_3, n, k_3)$ be elements of~$X_n$.
We want to show that
\begin{equation}\label{eqn:Triangle}
D^n(z_1, z_2) + D^n(z_1, z_3) \geq D^n(z_2,z_3).
\end{equation}

We distinguish the following cases:

\smallskip

\begin{itemize}
\item[Case 1:]  $k_1$, $k_2$, and $k_3$ are pairwise distinct.  In this case, we have by~(Dn1):
\begin{equation*}
    D^n(z_1, z_2) + D^n(z_1, z_3) = \beta_n + \beta_n > \beta_n = D^n(z_2, z_3).
\end{equation*}

\item[Case 2:]  $k_1 = k_2 \neq k_3$ or $k_1 = k_3 \neq k_2$.  Again by~(Dn1):
\begin{equation*}
    D^n(z_1,z_2) + D^n(z_1,z_3) \geq \beta_n = D^n(z_2,z_3).
\end{equation*}
\item[Case 3:]  $k_2 = k_3 \neq k_1$. Then by~(Dn1), (P$\delta$1), and (PDn1):
\begin{equation*}
    D^n(z_1,z_2) + D^n(z_1,z_3) = \beta_n + \beta_n \geq 2\delta_n > \eps_n \geq D^n(z_2,z_3).
\end{equation*}

\item[Case 4:]  $k_1 = k_2 = k_3 = k$.  When two of the points $z_1, z_2, z_3$ are equal, then \eqref{eqn:Triangle}
follows from reflexivity and symmetry for the nonnegative function~$D^n$.  Thus we only need to focus on the case
where $z_1, z_2, z_3$ are pairwise distinct, so that also $y_1$, $y_2$ and $y_3$ are pairwise distinct.

\smallskip

It suffices to consider the following three subcases:

\smallskip

\begin{itemize}
\item[Case 4-1:]  $\Delta(y_1, y_2) = 0$ and $\Delta(y_1, y_3) = 0$.

In this case $\Delta(y_2,y_3) > 0$, so that (Dn3) together with (P$\delta$1) implies:
\begin{equation*}
    D^n(z_1, z_2) + D^n(z_1, z_3) \geq \delta_n + \delta_n > \eps_n > \eps_n 3^{-\Delta(y_2,y_3)} = D^n(z_2, z_3).
\end{equation*}
\item[Case 4-2:]  $\Delta(y_1, y_2) = 0$ and $\Delta(y_1, y_3) > 0$.  Then  $\Delta(y_2,y_3) = 0$.

If $\Delta(y_1, y_3) \leq 2\lambda T^+(n)$, then (Dn3), (P$\delta$1), and (P$\delta$3) imply that:
\begin{equation*}
    D^n(z_1,z_2) + D^n(z_1, z_3) \geq \delta_n + \eps_n 3^{-2\lambda T^+(n)} \geq \eps_n \geq D^n(z_2, z_3).
\end{equation*}
If $\Delta(y_1, y_3) > 2\lambda T^+(n)$, then it follows from~\eqref{eqn:numberi} and~\eqref{deqn:define-Delta} that
\begin{equation*}
\begin{split}
&y_1\restrict (-\lambda T^+(n)+1, \dots, 0, \dots, \lambda T^+(n)-1)\\
 = \ &y_3\restrict ( -\lambda T^+(n)+1, \dots, 0, \dots, \lambda T^+(n)-1).\\
\end{split}
\end{equation*}
In particular, $\Phi((y_1, n, k)) = \Phi((y_3, n, k))$, and in view of~(Dn32) we have $D^n(z_1, z_2) = D^n(z_2, z_3)$.  Thus:
\begin{equation*}
    D^n(z_1,z_2) + D^n(z_1, z_3) > D^n(z_1, z_2) = D^n(z_2, z_3).
\end{equation*}
\item[Case 4-3:]  $\Delta(y_1, y_2) > 0$ and $\Delta(y_1, y_3) > 0$.

In this case $y_1(0) = y_2(0) = y_3(0)$ so that $\Delta(y_2,y_3) > 0$ and  it is sufficient to consider the following two subsubcases:

\smallskip

\begin{itemize}
\item[Case 4-3-1:] $\Delta := \Delta(y_1, y_2) = \Delta(y_1,y_3)$.

Then $\Delta(y_2,y_3) > \Delta$ and hence:
\begin{equation*}
    D^n(z_1,z_2) + D^n(z_1,z_3) = \eps_n 3^{-\Delta} + \eps_n 3^{-\Delta} > \eps_n 3^{-\Delta(y_2, y_3)} = D^n(z_2, z_3).
\end{equation*}
\item[Case 4-3-2:] $\Delta(y_1, y_2) > \Delta(y_1,y_3)$.

Then $\Delta := \Delta(y_2,y_3) = \Delta(y_1,y_3)$ and hence:
\begin{equation*}
    D^n(z_1, z_2) + D^n(z_1, z_3) = \eps_n 3^{-\Delta(y_1,y_2)} + \eps_n 3^{-\Delta} > \eps_n 3^{-\Delta} = D^n(z_2,z_3).
\end{equation*}
\end{itemize}
\end{itemize}

\end{itemize}
\end{itemize}

We have shown that $D^n$ is a metric on $X_n$.

\bigskip

\noindent
(ii)-(Pnc):  As metric spaces are Hausdorff, it suffices to show that $(X_n, D^n)$ is compact.

Consider any infinite sequence $(z_m)_{m=1}^{\infty} = ((y_m, n, k_m))_{m=1}^{\infty}$ of elements of $X_n$.

By repeatedly applying the Pigeonhole Principle, we can recursively construct infinite subsets $M_r$ for $r \in \bN$ such that
\begin{itemize}
\item There exists $k \in \{0, 1, \dots, T^+(n)-1\}$ such that $k_m = k$ for all~ $z_m \in M_0$.
\item $M_{r+1} \subseteq M_r$.
\item For all $r \geq 1$, there exists $f_r \in {}^{\{-r, -r+1, \dots, r-1, r \}}\{0,1\}$ such that \\$y_m \restrict \{-r, -r+1, \dots , r-1, r\} = f_r$ for all $z_m \in M_r$.
\end{itemize}

Note that the first two items imply that $f_r \subset f_{r+1}$ for all $r \geq 1$.

To construct $M_0$, note that for any fixed $n$  there are only finitely many possibilities for the value of
$k_m$.  Then by the Pigeonhole Principle, there exist \newline
$k \in \{0, 1, \dots, T^+(n)-1\}$ and an infinite subset $M_0$ of $(z_m)_{m=1}^{\infty}$ such that $k_m = k$ for all~ $z_m \in M_0$.

To construct $M_1$, by the fact that for any $i \in \bZ$ we have $y_m\restrict \{i\} \in \{0, 1\}$ for all $z_m \in M_0$, again by the Pigeonhole Principle, there exists $f_1 \in {}^{\{-1,0,1\}}\{0,1\}$ and an infinite subset $M_1$ of $M_0$ such that $y_m \restrict \{-1,0,1\} = f_1$ for all $z_m \in M_1$.

Now suppose $r \geq 1$ and the sets $M_0, \dots, M_{r}$ that satisfy the conditions we spelled out above are all constructed.
Then there exists $f_r \in {}^{\{-r, -r+1, \dots, r-1, r\}}\{0,1\}$ such that  $y_m \restrict \{-r, -r+1, \dots, r-1, r\} = f_r$ for all $z_m \in M_r$.  That is,\\ $y_m \restrict \{-r-1, -r, \dots, r, r+1\} = (y_m(-r-1), f_r, y_m(r+1))$ for all $z_m \in M_r$ where $y_m(-r-1), y_m(r+1) \in \{0, 1\}$.  Thus, again by the Pigeonhole Principle, there exits $f_{r+1} \in {}^{\{-r-1, -r, \dots, r, r+1\}}\{0,1\}$ and an infinite subset $M_{r+1}$ of $M_r$ such that\\ $y_m \restrict \{-r-1, -r, \dots, r, r+1\} = f_{r+1}$ for all $z_m \in M_{r+1}$.

This completes the construction of $M_r$.

Now choose
$z_{m_r} \in M_r$, where ${m_r} \neq {m_s}$ for $r \geq 1$ and $s < r$, and consider the subsequence  $(z_{m_r})_{r=1}^{\infty} = ((y_{m_r}, n, k))_{r=1}^{\infty}$.  Let
$y = \bigcup_{r \geq 1} f_r$.

Then $z:= (y, n, k) \in X_n$.  Moreover, if $y_{m_r} \neq y$, then  $\Delta(y_{m_r}, y) > r > 0$ for all $r \geq 1$.  Thus for computing
$D^n(z_{m_r}, z)$ either clause~(Dn2) or clause~(Dn31) will be used.  In either case,
$D^n(z_{m_r}, z) \leq \eps_n3^{-\Delta(y_{m_r}, y)} < \eps_n 3^{-r}$. It follows that the subsequence~$(z_{m_r})_{r=1}^{\infty}$ converges to~$z$.
Thus every sequence in $(X_n, D^n)$ has a subsequence that converges to a point in $X_n$, which proves compactness.

\bigskip

\noindent
(ii)-(PFn): We need to show three properties of~$F_n$:

\begin{itemize}
\item $F_n$ is onto:

We want to show that for all $z = (y,n,k) \in X_n$, there exists $z' = (y', n, k') \in X_n$ such that $F_n(z') = z$.  In fact, these $y'$ and $k'$ can simply be chosen  so that $y'(i) = y(i-1)$ for all $i$ and  $k' = (k-1)\ mod \ T^+(n)$.  Then it follows from the definition of~$F_n$ that $\sigma(y') = y$ and $\sigma(k') = k$, and therefore we have $F_n(z') = z$.

\smallskip

\item $F_n$ is one-to-one:

Consider any $z = (y,n,k)$ and $z' = (y',n,k')$ with $F_n(z) = F_n(z')$.  That is, $\sigma(y) = \sigma(y')$ and $F_n(k) = F_n(k')$.
Then we have $y(i+1) = y'(i+1)$ for all $i \in \bZ$, which implies that $y = y'$.
Similarly, $(k+1) \ mod \ T^+(n) = (k'+1) \ mod \ T^+(n)$ implies $k = k'$, as $0 \leq k,k' < T^+(n)$.
We conclude that $z = z'$, and hence $F_n$ is one-to-one.

\smallskip

\item $F_n$ is continuous:

Fix any $z = (y,n,k) \in X_n$ and any $\gamma > 0$. We want to show that
there exists $\delta > 0$ such that for all $u\in X_n$
\begin{equation}\label{eqn:contFn}
D^n(z,u) < \delta \ \Rightarrow D^n(F_n(z), F_n(u)) < \gamma.
\end{equation}

  Take $\delta = \min\{\frac{\eps_n}{3^4}, \frac{\gamma}{3^3}\}$.  Suppose $u = (y',n,k') \in X_n$ is such that $D^n(z,u) < \delta$. If  $u = z$ there is nothing to prove, so assume $u \neq z$.  Then $D^n(z,u) < \frac{\eps_n}{3^4}$ by the choice of~$\delta$ and $\frac{\eps_n}{3^4} < \beta_n$.  Thus, $k=k'$ and we must have $y \neq y'$ and $\Delta(y,y')> 0$, since $D^n(z,u)$ cannot be defined in terms of clause~(Dn32).  Thus clause~(Dn31) applies, and  $D^n(z,u) = \eps_n 3^{-\Delta(y,y')} < \delta \leq \frac{\eps_n}{3^4}$. This in turn implies that $\Delta(y,y') > 4$ and, in view of~\eqref{deqn:define-Delta}
  and~\eqref{eqn:numberi}, that $\Delta(\sigma(y),\sigma(y')) \geq \Delta(y,y') -2 > 2 > 0$.  Thus with $F_n(k) = F_n(k')$, clause~(Dn31) is also used for computing
  $D^n(F_n(z),F_n(u))$. Moreover, since $D^n(z,u) = \eps_n 3^{-\Delta(y,y')} < \delta$, we have the inequality $3^{-\Delta(y,y')} < \frac{\delta}{\eps_n}$.  Therefore,
\begin{equation*}
    \begin{split}
        D^n(F_n(z),F_n(u)) &= \eps_n 3^{-\Delta(\sigma(y),\sigma(y'))} \leq \eps_n 3^{-(\Delta(y,y')-2)} \\
        &= 9\eps_n 3^{-\Delta(y,y')} < 9\eps_n\left(\frac{\delta}{\eps_n}\right)\\
        &= 9\delta \leq \frac{\gamma}{3} < \gamma.
    \end{split}
\end{equation*}
Thus the implication~\eqref{eqn:contFn} holds, and we conclude that $F_n$ is continuous.
\end{itemize}

\medskip

By (Pnc) and since a continuous bijection from a compact Hausdorff space~$X$ onto a compact Hausdorff space~$Y$ is automatically a homeomorphism (see, for example,  Proposition A.1.11  at page 705 of~\cite{KH}),
we conclude that $F_n$ is a homeomorphism.
$\Box$

\section{Construction of \JX-systems $(X, F)$ and \JX-metrics~$D$}\label{sec:X}

\begin{definition}\label{def:JX}
Let~$(X_n, F_n)_{n\in \bN}$ be a sequence of \JX-systems with \JX-metrics~$D^n$. Then the following construction defines a \emph{\JX-system~$(X,F)$} with
\emph{\JX-metric~$D$.}

\begin{itemize}
\item $X$:  Let $X = \prod_{n \in \bN} X_n$.  That is, we let $X$ consist of all sequences
\newline $x = (x_n)_{n \in \bN}$ such that $x_n \in X_n$ for each~$n \in \bN$.
\item $F$:  For  $x \in X$, define $F(x)_n = F_n(x_n)$ for all~$n \in \bN$.
\item $D$:  The function~$D: X^2 \rightarrow [0,\infty)$ is defined as:
\begin{equation}\label{eqn:define-D}
D(x, x') = \sum_{n \in \bN} D^n(x_n, x'_n).
\end{equation}
\end{itemize}
\end{definition}

Note that there is exactly one pair $(X, F)$ that can be an \JX-system, but for a variety of choices of \JX-metrics~$D$.
In the remainder of this note, we will  call~$X$ \emph{the \JX-space} and reserve the symbols~$X, F, D$ always for the objects that are constructed in
Definition~\ref{def:JX}. The notation chosen in the statements of Theorems~\ref{thm:main}--\ref{thm:main-spanonly}
also conforms to this convention.

\smallskip

\JX-metrics will indeed be metrics by Proposition~\ref{prop:D-OK} below.  Moreover, the \JXn-metrics and \JX-metrics used in the proofs of Theorems~2 and~3 will satisfy:
\begin{equation}\label{eqn:large-dist}
\begin{split}
diam(X_n, D^n) &= \eps_n \ \  \mbox{for all} \ \ n \in \bN,\\
diam(X,D) &= \eps,\\
D_T(x, x') &= \eps \ \Leftrightarrow \ \exists 0 \leq t < T \, \forall n \in \bN \quad D^n(F_n^t(x_n), F_n^t(x'_n)) = \eps_n,
\end{split}
\end{equation}
where the third line follows directly from the first line and Definition~\ref{def:JX}.

\begin{proposition}\label{prop:D-OK}
The function~$D$ as defined in~\eqref{eqn:define-D} is a metric on~$X$ that induces the product topology.
\end{proposition}

\noindent
\textbf{Proof:}  By (P$\eps$) and (PDn1), the series in \eqref{eqn:define-D} is convergent.  Then from
Theorem~4.2.2  at page 259 of \cite{RE}, and the statement following its proof,
we can deduce that $D$ is a metric on $X$ that induces the topology of the Cartesian product of the spaces $\{X_n\}_{n=0}^{\infty}$. $\Box$

\begin{proposition}\label{prop:X-OK}
The state space $(X,D)$ of an \JX-system is compact in the product topology, and $(X, F)$ is the product of the
\JX-systems~$(X_n, F_n)$.  In particular, $F$ is a homeomorphism.
\end{proposition}

\noindent
\textbf{Proof:} By Proposition \ref{prop:metric-n}, for each $n\in\bN$, the function $D^n$ is a metric on $X_n$ that induces
the topology of a compact Hausdorff space.  By Proposition \ref{prop:D-OK}, $D$ is a metric on $X$ that induces
the product topology.  Then by Tychonoff's theorem, $X = \prod_{n\in\bN}X_n$ is compact in the product
topology.

It follows directly from  the construction of $(X,D,F)$ that  $(X,F)$ is the product of the systems $(X_n,F_n)$.

It remains to show that $F$ is a homeomorphism.

First, $F$ is a bijection on $X$ as each $F_n$ is a bijection on $X_n$.

To show that $F$ is continuous, consider $x \in X$ and $\gamma > 0$.
We want to show that  there exists $\delta > 0$ such that
\begin{equation}\label{eqn:contF}
\forall x' \in X \ \left(x' \neq x \ \& \ D(x,x') < \delta \ \ \Rightarrow \ \ D(F(x), F(x')) < \gamma\right).
\end{equation}

Fix a $K(\gamma) \in \bN$ such that $\eps_{K(\gamma)} < \gamma$ and let  $x' \in X$.  Then by (P$\delta$2):
\begin{equation*}
\begin{split}
        \sum_{n = K(\gamma)+1}^{\infty} D^n(F_n(x_n),F_n(x'_n)) &\leq
        \sum_{n = K(\gamma)+1}^{\infty} \eps_n\\
         &< 0.5(\eps_{K(\gamma)}-\delta_{K(\gamma)}) \\
        &< 0.5\eps_{K(\gamma)}\\
        &< 0.5\gamma.
\end{split}
\end{equation*}

For each $0 \leq n \leq K(\gamma)$, by the continuity of $F_n$, there exists $\eta_n >0$ such that
\begin{equation*}
\forall x' \in X \ \left(x' \neq x \ \& \ D^n(x_n,x'_n) < \eta_n \ \ \Rightarrow \ \ D^n(F_n(x_n), F_n(x'_n)) < \frac{\gamma}{2K(\gamma)+2}\right).
\end{equation*}

Let $\delta = \min\{\eta_n: 0 \leq n \leq K(\gamma)\}$.  Then, for any $x' \in X$ with $x' \neq x$ and $D(x,x') < \delta$,
\begin{equation*}
    \begin{split}
        D(F(x), F(x')) &= \sum_{n\in\bN} D^n(F_n(x_n), F_n(x'_n))\\
        &= \sum_{n=0}^{K(\gamma)} D^n(F_n(x_n), F_n(x'(n))) + \sum_{n=K(\gamma)+1}^{\infty} D^n(F_n(x_n), F_n(x'(n)))\\
        &< (K(\gamma)+1)\frac{\gamma}{2K(\gamma)+2} + 0.5\gamma\\
        &= \gamma.
    \end{split}
\end{equation*}

Thus~\eqref{eqn:contF} holds, and we conclude that $F$ is continuous.

Since $X$ is compact and $F$ is a continuous bijection, the result quoted above implies that $F^{-1}$ is also
continuous, so that $F$ is a homeomorphism.
$\Box$

\section{Construction of $\Ymn, W$, and $X^-$}\label{sec:X-}

\subsection{Construction of $\Ymn$}\label{subsec:Xn-Yn-}

The state spaces $X^-, W$ of the systems in Theorems~\ref{thm:main} and~\ref{thm:main-seponly} will be subspaces~$Y$ of the~\JX-space~$X$.    The key to our arguments is deriving lower bounds on $sep\left(Y, \eps, D_{T^+(n)}\right)$ and $span\left(Y, \eps, D_{T^+(n)}\right)$ and then upper bounds on
$sep\left(X, \eps, D_{2T(n)}\right)$ and $span\left(X, \eps, D_{2T(n)}\right)$ for certain~\JX-metrics~$D$.  For the former, we will need the following notion.

\begin{definition}\label{def:Yn-}
An \emph{RY-sequence} $(\Ymn)_{n \in \bN}$ is a sequence of subsets $\Ymn \subset {}^{T^+(n)}\{0,1\}$ such that:
\begin{itemize}
\item For $n = 0$:
\begin{itemize}
    \item[(PY1)]  For all $\varphi,\psi \in \Ymz$ we have $\varphi \neq \psi \Rightarrow C(\varphi, \psi) > 2$ \newline
    (that is, $\varphi$ and $\psi$ differ on at least 3 intervals $I^0_j$).
    \item[(PY2)]  for all $\varphi \in \Ymz$, there exists $0 < t \leq T(0)-1$ such that $\varphi (t) = 1$.
\end{itemize}
\item For each~$n > 0$:
 \begin{itemize}
 \item[(PR1)] Each sequence in $\Ymn$ consists of consecutive blocks of length $T^+(n-1)$ in $\Ymnm$.
 \item[(PR2)] For all $\varphi,\psi \in \Ymn$ we have $\varphi \neq \psi \Rightarrow C(\varphi, \psi) > 2$ \newline
  (that is, $\varphi$ and $\psi$ differ on at least 3 intervals $I^n_j$).
\end{itemize}
\end{itemize}
\end{definition}

Note that the recursive construction of the sets $\Ymn$ from sets $\Ymnm, \dots , \Ymz$ guarantees the following generalization of condition~(PY2):

\begin{itemize}
\item[(PY2+)] For all $n \in \bN$, non-negative integer multiple~$\tau$ of~$T^+(0)$ with~$\tau < T^+(n)$, and all $\varphi \in \Ymn$, there exists $0 < i \leq T(0)-1$ such that $\varphi (\tau + i) = 1$.
\end{itemize}

\begin{remark}\label{rem:PY}
Property~(PY1) guarantees that if $c_0: \left[{}^{T^+(0)}\{0,1\}\right]^2 \rightarrow [C(0)]$ is any coloring that satisfies conditions~(cC1) and~(cC2),
then for $\varphi \neq \psi \in \Ymz$ and~$j = c_0(\varphi, \psi)$ we have $\varphi \restrict I_j^0 \neq \psi \restrict I_j^0$.
Property~(PY2) guarantees that the same will be true if~$c_0$ satisfies conditions~(cC) and~(cCi), as
in conjunction with~\eqref{T+(-1)} it guarantees that  the value of~$c_0(\varphi, \psi)$ will not be automatically
determined by condition~(cCi).
\end{remark}

For the remainder of this paper we fix an RY-sequence $(\Ymn)_{n \in \bN}$ that satisfies~\eqref{size:Y_n^-LowerBound} of the following result.

\bigskip

\begin{claim}\label{size}
There exists an RY-sequence $(\Ymn)_{n \in \bN}$ such that for all $n \in \bN$ the following inequality holds:
\begin{equation}\label{size:Y_n^-LowerBound}
|\Ymn|  \geq 2^{0.9T^+(n)}.
\end{equation}
\end{claim}

\noindent
\textbf{Proof:}
We show that one can recursively choose subsets $\Ymn \subset {}^{T^+(n)}\{0,1\}$   with properties~(PY1), (PY2),
(PR1), and (PR2) so  that for each $n \in \bN$  the following inequality holds
\begin{equation}\label{eqn:size-of-Xn-}
|\Ymn| \geq \frac{2^{[\prod_{i=0}^n (C(i)-2)][\prod_{i=0}^n K(i)]}}{{C(n) \choose 2} \prod_{m=0}^{n-1}\left[{C(m) \choose 2}^{\prod_{i=m+1}^{n}[(C(i)-2)K(i)]}\right]}
\geq 2^{0.9T^+(n)}.
\end{equation}

By \eqref{eqn:TT+(n)-explicit}, the second inequality in~\eqref{eqn:size-of-Xn-} is equivalent to
\begin{equation*}
    \frac{2^{[\prod_{i=0}^n (C(i)-2)][\prod_{i=0}^n K(i)]}}{2^{0.9[\prod_{i=0}^n C(i)][\prod_{i=0}^n K(i)]}} \geq {C(n) \choose 2} \prod_{m=0}^{n-1}\left[{C(m) \choose 2}^{\prod_{i=m+1}^{n}[(C(i)-2)K(i)]}\right].
\end{equation*}
By (PCn) and~\eqref{eqn:TT+(n)-explicit},
\begin{equation*}
    \begin{split}
    LHS &= 2^{\left[\prod_{i=0}^n K(i)\right]\left[\prod_{i=0}^n (C(i)-2) - 0.9\prod_{i=0}^n C(i)\right]}\\
    &> 2^{\left[\prod_{i=0}^n K(i)\right]\left[0.95\prod_{i=0}^n C(i) - 0.9\prod_{i=0}^n C(i)\right]}\\
    &= 2^{0.05\left[\prod_{i=0}^n C(i)\right]\left[\prod_{i=0}^n K(i)\right]}\\
    &= 2^{0.05T^+(n)}.
    \end{split}
\end{equation*}
Now (PKn2) implies that
\begin{equation*}
    2^{0.05T^+(n)} > RHS.
\end{equation*}

In the proof for~$n = 0$, we start with a list~$L(0,1)$
that will be a bijective
enumeration of the set of all $\varphi \in {}^{T^+(0)} \{0,1\}$ that satisfy property~(PY2):
\begin{equation*}
\begin{split}
L(0,1) &= {}^{T^+(0)}\{0,1\} \backslash \left[\left\{{}^{T(0)}\{0\}\times {}^{\{T(0), T(0)+1, \dots, T^+(0)-1\}}\{0,1\}\right\}\cup\{\varphi_0\}\right]\\
&= \{\varphi_i\}_{i=1}^{2^{T^+(0)}-2^{(C(0)-1)T(0)}-1},
\end{split}
\end{equation*}
where $\varphi_0$ is the function that takes the value~$\varphi(0) = 1$ and the value~$\varphi(t) = 0$
for~$t > 0$.

Thus the length~$\ell$ of the list~$L(0,1)$ is given by
\begin{equation}\label{eqn:ell-Z}
\ell = 2^{T^+(0)}-2^{(C(0)-1)T(0)}-1.
\end{equation}

We recursively construct lists $L(0, i)$ for $i = 2, 3, \dots$
by (possibly) removing some elements of $L(0,i)$ to obtain $L(0,i+1)$ as follows:
\begin{itemize}
\item If $\varphi_i \notin L(0, i)$, then we let $L(0, i+1) = L(0, i)$.
\item If $\varphi_i \in L(0, i)$, then we obtain $L(0, i+1)$ by  removing from $L(0,i)$ all $\varphi_j \in L(0,i)$ with $j > i$ that differ from $\varphi_i$ on at most two intervals $I^0_j$.
\end{itemize}

Let   $\Ymz ~= ~L(0, ~2^{T^+(0)}-2^{(C(0)-1)T(0)})$ be the set of $\varphi_i$ that survived this procedure.
By construction, this family has both properties~(PY1) and~(PY2).

Since removal is always conditioned on a prior decision to retain some~$\varphi_i$ in~$L(0,i)$,  there are at
most~$|\Ymz|$ steps $i$ where any removal took place, that is, where $L(0, i+1) \neq L(0, i)$.
Moreover, at each such step we could have removed
 at most $NR_0 = {C(0) \choose 1}\left(2^{T(0)}-1\right)+{C(0) \choose 2}\left(2^{T(0)}-1\right)^2$ elements.
 Thus
\begin{equation*}
    \begin{split}
        |\Ymz| &\geq \ell - NR_0|\Ymz| \\
        \ell &\leq \left(1+NR_0\right)|\Ymz|\\
        &< {C(0) \choose 2}\left(2^{2T(0)}- 2^{T(0)}\right)|\Ymz|,
    \end{split}
\end{equation*}
where the last inequality follows from the observation that
\begin{equation*}
     \begin{split}
     1 + NR_0 &= 1+C(0)\left(2^{T(0)}-1\right) + {C(0) \choose 2}2^{2T(0)} + {C(0) \choose 2} - {C(0) \choose 2}2^{T(0)+1}\\
     &= \left[{C(0) \choose 2}2^{2T(0)} - {C(0) \choose 2}2^{T(0)}\right]\\
     &\ \ \ \ + \left[1-C(0)\right]+ 2^{T(0)}\left[C(0)-{C(0) \choose 2}\right] + {C(0) \choose 2}.
     \end{split}
 \end{equation*}

By~\eqref{eqn:ell-Z} and \eqref{eqn:TT+(n)},
\begin{equation}\label{eqn:size-of-Y0-}
\begin{split}
|\Ymz| &> \frac{2^{T^+(0)}-2^{C(0)T(0)-T(0)}-1}{{C(0) \choose 2}\left[2^{2T(0)}-2^{T(0)}\right]}\\
|\Ymz| &\geq \frac{2^{T^+(0)}-2^{C(0)T(0)-T(0)}}{{C(0) \choose 2}\left[2^{2T(0)}-2^{T(0)}\right]}\\
&= \frac{2^{C(0)K(0)}-2^{C(0)K(0)-K(0)}}{{C(0) \choose 2}\left[2^{2K(0)}-2^{K(0)}\right]}\\
&= \frac{2^{C(0)K(0)}\left[1-2^{-K(0)}\right]}{{C(0) \choose 2}2^{2K(0)}\left[1-2^{-K(0)}\right]}\\
&= \frac{2^{(C(0)-2)K(0)}}{{C(0) \choose 2}}.
\end{split}
\end{equation}

Notice that for $n = 0$ the product $\prod_{m=0}^{n-1}\left[{C(m) \choose 2}^{\prod_{i=m+1}^{n}[(C(i)-2)K(i)]}\right]$
has no terms and is treated as equal to~1, so
that~\eqref{eqn:size-of-Y0-} is equivalent to the first inequality in~\eqref{eqn:size-of-Xn-} for the special case~$n = 0$.

\medskip

Now assume by induction that for a fixed $n \geq 0$ we have already constructed~$\Ymn$ so that, in particular,
\begin{equation}\label{eqn:Ymn-indass}
|\Ymn| \geq \frac{2^{[\prod_{i=0}^n (C(i)-2)][\prod_{i=0}^n K(i)]}}{{C(n) \choose 2} \prod_{m=0}^{n-1}\left[{C(m) \choose 2}^{\prod_{i=m+1}^{n}[(C(i)-2)K(i)]}\right]}.
\end{equation}

To obtain $\mathcal{Y}_{n+1}^-$ with properties~(PR1), (PR2) and
the desired lower bound for $|\mathcal{Y}_{n+1}^-|$, arrange those elements of ${}^{T^+(n+1)}\{0,1\}$ that consist of  blocks of length $T^+(n)$ in $\Ymn$ as $ \{\varphi_i\}_{i=1}^{|\Ymn|^{C(n+1)K(n+1)}} = L(n+1,1)$ into a list.
We recursively construct lists $L(n+1, i)$ for $i = 2, \dots , |\Ymn|^{C(n+1)K(n+1)}+1$ by (possibly) removing some elements of $L(n+1,i)$ to obtain $L(n+1,i+1)$ as follows:
\begin{itemize}
\item If $\varphi_i \notin L(n+1, i)$, then we let $L(n+1, i+1) = L(n+1, i)$.
\item If $\varphi_i \in L(n+1, i)$, then we obtain $L(n+1, i+1)$ by  removing all $\varphi_j \in L(n+1,i)$ with $j > i$ from $L(n+1,i)$ that differ from $\varphi_i$ on at most two intervals $I^{n+1}_j$.
\end{itemize}

Let $\mathcal{Y}_{n+1}^- = L(n+1,  |\Ymn|^{C(n+1)K(n+1)}+1)$ be the set of $\varphi_i$ that survived this procedure.

Note that again there are at most $|\mathcal{Y}_{n+1}^-|$ steps $i$ where any removal took place, that is, where $L(n+1, i+1) \neq L(n+1, i)$.
Moreover, by \eqref{eqn:TT+(n)} and the specification of our construction, at each such step we could have removed
 at most \\$NR_{n+1} = {C(n+1) \choose 1}\left(|\Ymn|^{K(n+1)}-1\right)+{C(n+1) \choose 2}\left(|\Ymn|^{K(n+1)}-1\right)^2$ elements.  Thus
\begin{equation}\label{eqn:Yn+1-est}
    \begin{split}
        |\mathcal{Y}_{n+1}^-| \geq \ &|\Ymn|^{C(n+1)K(n+1)}
        - NR_{n+1}|\mathcal{Y}_{n+1}^-| \\
        |\Ymn|^{C(n+1)K(n+1)} \leq &\left(1+NR_{n+1}\right)|\mathcal{Y}_{n+1}^-| \\
        &< {C(n+1) \choose 2}|\Ymn|^{2K(n+1)}|\mathcal{Y}_{n+1}^-|,
    \end{split}
\end{equation}
where the last inequality follows from the observation that
\begin{equation*}
    \begin{split}
        1 + NR_{n+1} &=1+C(n+1)\left(|\Ymn|^{K(n+1)}-1\right) + {C(n+1) \choose 2}|\Ymn|^{2K(n+1)}\\
        &\ \ \ +  {C(n+1) \choose 2} - 2{C(n+1) \choose 2}|\Ymn|^{K(n+1)}\\
        &=
        \left[1-C(n+1)\right]+ |\Ymn|^{K(n+1)}\left[C(n+1)-{C(n+1) \choose 2}\right]  \\
        &\ \ \ + {C(n+1) \choose 2}\left[1-|\Ymn|^{K(n+1)}\right] + {C(n+1) \choose 2}|\Ymn|^{2K(n+1)}.
    \end{split}
\end{equation*}

From inequality~\eqref{eqn:Yn+1-est} and the inductive assumption~\eqref{eqn:Ymn-indass} we infer
\begin{equation*}
    \begin{split}
       |\mathcal{Y}_{n+1}^-| &> \frac{|\Ymn|^{(C(n+1)-2)K(n+1)}}{{C(n+1) \choose 2}} \\
       &> \frac{\left[\frac{2^{[\prod_{i=0}^n (C(i)-2)][\prod_{i=0}^n K(i)]}}{{C(n) \choose 2} \prod_{m=0}^{n-1}\left[{C(m) \choose 2}^{\prod_{i=m+1}^{n}[(C(i)-2)K(i)]}\right]}\right]^{(C(n+1)-2)K(n+1)}}{{C(n+1) \choose 2}}\\
       &= \frac{2^{[\prod_{i=0}^{n+1} (C(i)-2)][\prod_{i=0}^{n+1} K(i)]}}{{C(n+1) \choose 2} \prod_{m=0}^{n}\left[{C(m) \choose 2}^{\prod_{i=m+1}^{n+1}[(C(i)-2)K(i)]}\right]}.
    \end{split}
\end{equation*}

This recursive construction gives an RY-sequence $(\Ymn)_{n \in \bN}$ for which~\eqref{eqn:size-of-Xn-} follows by induction. $\Box$

\subsection{Construction of $W, W^n$, and $X^-$}\label{subsec:WX-}
Let $(X,F)$ be an \JX-system. Here and in much of our subsequent work we adopt the following notation:
\begin{itemize}
\item For $x \in X$ and $n \in \bN$, the $n$-th coordinate of~$x$ will be denoted by \newline $x_n = (y_n, n, k_n)$.
\item For each~$n \in \bN$ and $\varphi \in {}^{T^+(n)}\{0,1\}$, we let $y_{\varphi} \in {}^{\bZ}\{0,1\}$ be such that\\
$y_{\varphi}\restrict(0,\dots,T^+(n)-1) = \varphi$ and $y_{\varphi}(i) = 0$ when $i \geq T^+(n)$ or $i \leq -1$.
\item We let $x^\varphi$ denote the element of~$X$ such that $x^{\varphi}_n = (y_{\varphi}, n, 0)$ for all $n \in \bN$.
\item $W^n := \{x^{\varphi}: \varphi \in \Ymn\} \subset X$.
\end{itemize}

There is a one-to-one correspondence between the
set~$W^n$ and the set~$\Ymn$.
Thus
Claim~\ref{size} implies:
\begin{corollary}\label{corol:size}
 For all $n \in \bN$ we have
$|W^n|  \geq 2^{0.9T^+(n)}$.
\end{corollary}

The sets~$X^-$ and~$W$ in the statements of Theorems~\ref{thm:main} and~\ref{thm:main-seponly} will be the following subsets of~$X$:
\begin{equation*}
\begin{split}
  X^- &= \overline{\bigcup_{t\in\bZ}F^t\left(\bigcup_{n\geq 0} W^n\right)},\\
    W &= \{x\in X:  \exists y \in {}^{\bZ}\{0,1\}\, \forall n \in \bN \ \ y_n = y\ \mbox{and} \ k_n = k_{n+1} \mod T^+(n)\}.
\end{split}
\end{equation*}

For simplicity, we will usually not make a notational distinction between $F, F\restrict X^-$, and $F \restrict W$.

\smallskip

Notice that
$X^-$ is by definition a closed subset of~$(X, D)$.  Moreover, $F$ is by definition forward and backward
invariant on the set $\bigcup_{t\in\bZ}F^t\left(\bigcup_{n\geq 0} W^n\right)$ and its closure~$X^-$ in~$X$.

\smallskip

Thus  part~(i) of Theorem~\ref{thm:main}  follows from Proposition~\ref{prop:X-OK}.  Similarly, part~(i) of Theorem~\ref{thm:main-seponly} is a consequence o the following observations about the set~$W$.

\begin{proposition}\label{prop:W}
Let~$W$ be defined as above, and let $D$ be any \JX-metric.  Then

\smallskip

\noindent
(i) Each~$x^\varphi \in W$.  In particular, $W \neq \emptyset$.

\smallskip

\noindent
(ii) The set~$W$ is closed in the space~$(X, D)$.

\smallskip

\noindent
(iii) The set~$W$ is both forward and backward invariant under~$F$.
\end{proposition}

\noindent
\textbf{Proof:} Part~(i) follows directly from the definitions.

\smallskip

For part~(ii), note that if $x \notin W$, then there must exist $n \in \bN$ such that  at least one of the following holds:

\smallskip

\begin{itemize}
\item[Case 1:] $y_n \neq y_{n+1}$.\\
Then $\Delta : = \Delta(y_n, y_{n+1}) < \infty$, and for each $x' \in W$ we must have\\
$\Delta(y_n, y'_{n}) \leq \Delta$ or $\Delta(y_{n+1}, y_{n+1}') \leq \Delta$, so that
\begin{equation*}
\begin{split}
D(x, x') &\geq D^n(x_n, x_n') + D^{n+1}(x_{n+1}, x_{n+1}')\\
&\geq \min \{\delta_{n+1} ,  \eps_{n+1} 3^{-\Delta}\}\\
&> \eps_{n+1} 3^{-\Delta - 1}.
\end{split}
\end{equation*}
\smallskip

\item[Case 2:] $k_n \neq k_{n+1} \mod T^+(n)$.\\
Then for each $x' \in W$ we must have
$k_n \neq k_n'$ or $k_{n+1} \neq k_{n+1}'$, so that
\begin{equation*}
\begin{split}
D(x, x') &\geq D^n(x_n, x_n') + D^{n+1}(x_{n+1}, x_{n+1}')\\
&\geq \delta_{n+1}.
\end{split}
\end{equation*}
\end{itemize}

In either case, we find an open ball around~$x$ that is disjoint from~$W$.  Thus $X \backslash W$ is open, and~$W$ is closed.

\smallskip

For part~(iii), recall the definition~\eqref{eqn:def-Fn} of the maps~$F_n$ on the coordinates of~$X$:
\begin{equation*}
\begin{split}
F_n((y, n, k)) &= (\sigma(y), n, F_n(k)), \ \mbox{where}\\
\sigma(y)(i) &= y(i+1) \ \mbox{ for all } i,\\
F_n(k) &= (k+1) \ mod \ T^+(n).
\end{split}
\end{equation*}

Here $\sigma(y)$ does not depend on~$n$, so that for $x \in W$ and $y \in {}^\bZ \{0,1\}$ such that $y = y_n$ for all~$n \in \bN$ we will have
$F(x)_n = (\sigma(y), n, F_n(k_n))$ for all~$n$; similarly for $F^{-1}(x)$. Thus $F(x)$ and $F^{-1}(x)$ retain the property of having the same $y$-component on all coordinates.

Similarly, if $k_n = k_{n+1} \mod T^+(n)$, then
\begin{equation*}
\begin{split}
F_n(k_n) &= (k_n + 1) \mod T^+(n) \\
&= (k_{n+1} + 1) \mod T^+(n)\\
&=  ((k_{n+1} + 1) \mod T^+(n+1)) \mod T^+(n)\\
&= F_{n+1}(k_{n+1}) \mod T^+(n).
\end{split}
\end{equation*}
The third of the above equalities follows from our choice of $T^+(n+1)$ as an integer multiple of~$T^+(n)$.
Thus the consistency property of the components $k_n$ in the definition of~$W$ is preserved by~$F$.  The analogous argument shows that it is also preserved by~$F^{-1}$, and we obtain part~(iii) of the proposition.
$\Box$

\bigskip

Note that:
\begin{equation}\label{eqn:X-subW}
\begin{split}
\forall n \in \bN \ \ W^n &\subset W \cap X^-,\\
X^- &\subseteq W.
\end{split}
\end{equation}

The first line of~\eqref{eqn:X-subW} follows immediately from our definitions; the second line then follows  from
the definition of~$X^-$ and  Proposition~\ref{prop:W}.

\subsection{Some properties of the systems $(X, F), (W, F\restrict W)$, and $(X^-, F\restrict X^-)$}\label{subsec:propsWX-}

Here we prove all parts of Theorems~\ref{thm:main} and~\ref{thm:main-seponly}, except parts~(i) that were already shown in the previous subsection and parts~(ii) that will be derived in the next two sections.

The following result proves Theorem~\ref{thm:main-seponly}(iiia) for $\delta^* := \sum_{n\in\bN} \delta_n$.

\begin{lemma}\label{lem:span0}
Let $(X, F)$ be an \JX-system with \JX-metric~$D$ that is constructed based on conditions~(Dn1d) for all components~$D^n$.  Let $Y \subset X$ be a closed subset that is invariant under~$F$. Then the corresponding \JX-system $(Y, F\restrict Y)$ satisfies:
\begin{equation*}
    \forall \delta > \sum_{n\in\bN} \delta_n\ \ \ \lim_{T\rightarrow\infty} \frac{\ln span(Y,\delta, D_T)}{T} = 0.
\end{equation*}
\end{lemma}

\noindent
\textbf{Proof:}  For $\delta > \sum_{n\in\bN}\delta_n$, we distinguish the following three cases:
\begin{itemize}
\item[(a)]:  $\delta > \eps$.
\item[(b)]:  $\delta = \eps$.
\item[(c)]:  $\sum_{n\in\bN} \delta_n < \delta < \eps$.
\end{itemize}

\smallskip

We will show that in all three cases we can find a fixed finite subset $S \subset Y$ that is $(T, \delta)$-spanning for all $T > 0$.  Then
\begin{equation*}\label{eqn:span-S-fin}
0 \leq \liminf_{T\rightarrow\infty} \frac{\ln span(Y,\delta, D_T)}{T} \leq \limsup_{T\rightarrow\infty} \frac{\ln span(Y,\delta, D_T)}{T} \leq \limsup_{T\rightarrow\infty} \frac{\ln |S|}{T} = 0,
\end{equation*}
and the result follows.

\smallskip

\noindent
(a):  In the case of $\delta > \eps$, choose any $x \in Y$, and let $S = \{x\} \subset Y$.  Then for all $x' \in Y$ and $t \geq 0$, by the definition of $D$ and conditions (P$\eps$) and (PDn1):
\begin{equation*}
    \begin{split}
    D(F^t(x), F^t(x')) &= \sum_{n\in\bN} D^n(F^t_n(x_n), F^t_n(x'_n))\\
    &\leq \sum_{n\in\bN}\eps_n\\
    &= \eps\\
    &< \delta.
    \end{split}
\end{equation*}

\smallskip

\noindent
(b):  In the case of $\delta = \eps$, fix any $x \in Y$ and let $S = \{x, F(x)\}$. Let $x' \in Y$. Then there exist $x'' \in S$ such that
for $x'_0 = (y'_0, 0, k'_0)$ and $x''_0 = (y''_0, 0, k''_0)$ the inequality~$k'_0 \neq k''_0$ holds. Thus by clause~(Dn1d) in the definition of $D^0$, together with conditions (P$\eps$), (PDn1), and (P$\delta$1):
\begin{equation*}
    \begin{split}
    D(F^t(x'), F^t(x'')) &= \sum_{n\in\bN} D^n(F^t_n(x'_n), F^t_n(x'_n))\\
    &= D^0(F^t_0(x'_0), F^t_0(x''_0)) + \sum_{n\geq 1} D^n(F^t_n(x'_n), F^t_n(x''_n))\\
    &\leq \delta_0 + \sum_{n\geq 1} \eps_n\\
    &= \delta_0 - \eps_0 + \sum_{n\in\bN} \eps_n \\
    &= \eps - (\eps_0 - \delta_0)\\
    &< \eps = \delta.
    \end{split}
\end{equation*}

\smallskip

\noindent
(c):  Suppose $\sum_{n\in\bN} \delta_n < \delta < \eps$.\\
There exists $K \in \bN$ such that $\sum_{n=K+1}^{\infty} \eps_n < \delta - \sum_{n\in\bN}\delta_n$.  Fix such a $K$, and let\\
$L = \{\bbk = (k_i)_{i = 0}^{K}:  \forall 0 \leq i \leq K\ \ k_i \in \{0, 1\} \}$.  \\
For each $\bbk = (k_i)_{i = 0}^{K} \in L$, pick $x(\bbk) \in Y$ with coordinates $x(\bbk)_i = (y(\bbk)_i, i , k_i)$ for all~$0 \leq i \leq K$ if possible; otherwise
let $x(\bbk)$ be an arbitrary element of~$Y$. \\
Now let $S = \{x(\bbk): \ \bbk \in L\}$. Then $|S| \leq |L| = 2^{K+1}$.

We show that $S$ is a $(T, \delta)$-spanning set in~$Y$ for all~$T>0$. Let $x' \in Y$, and let $\bbk \in L$ be the sequence
$\bbk = (F(k_i'))_{i=0}^K$, where $x'_i = (y'_i, i, k'_i)$ for all relevant~$i$. Let $x = x(\bbk) \in S$, with $x_i = (y_i, i, k_i)$ for all relevant~$i$.
Then $k_i + t \mod T^+(i) = k_i' + 1 + t \mod T^+(i) \neq k_i' + t \mod T^+(i)$ for all $0 \leq i \leq K$ and $t \geq 0$, and  clause~(Dn1d) will apply in the calculations
of~$D^i(F^t_i(x_i), F^t_i(x'_i))$.  Together with the choice of~$K$ and properties (P$\eps$) and (PDn1),  this implies for all $t \geq 0$:
\begin{equation*}
    \begin{split}
    D(F^t(x), F^t(x')) &= \sum_{n\in\bN} D^n(F^t_n(x_n), F^t_n(x'_n))\\
    &= \sum_{n=0}^K D^n(F^t_n(x_n), F^t_n(x'_n)) + \sum_{n=K+1}^{\infty} D^n(F^t_n(x_n), F^t_n(x'_n))\\
    &\leq \sum_{n=0}^K \delta_n + \sum_{n=K+1}^{\infty} \eps_n \\
    &< \sum_{n\in\bN} \delta_n + \left(\delta - \sum_{n\in\bN}\delta_n\right)\\
    &= \delta.
    \end{split}
\end{equation*}
$\Box$

\begin{lemma}\label{lem:span-lb}
Let $(X,F)$ be an \JX-system with \JX-metric~$D$.  Then for each $i \in \bN$ there exists $m \in \bN$ such that
\begin{equation}\label{eqn:spanbounds-delta-i}
 \forall\, T > 0 \quad    2^{(i+1)T} \leq span(X, \delta_i, D_T) \leq \left[\prod_{j=0}^i T^+(j)\right]2^{(i+1)(T+m-1)}.
\end{equation}
\end{lemma}

\begin{remark}\label{rem:span-upper-bd-redundant}
The upper bounds for $span(X, \delta_i, D_T)$ in~\eqref{eqn:spanbounds-delta-i} are not strictly needed for the proof of any parts of our theorems.
We included them here to round out the exposition.
\end{remark}

\noindent
\textbf{Proof of Lemma~\ref{lem:span-lb}:} Recall that the coordinates of $x, x' \in X$ are denoted by $x_n = (y_n, n, k_n)$ and $x'_n = (y'_n, n, k'_n)$. Instead of~$n$, we will use~$i$ or~$j$ as subscripts.

\smallskip

Throughout this proof, fix any $i \in \bN$.

\medskip

For the proof of the first inequality in~\eqref{eqn:spanbounds-delta-i},
assume towards a contradiction that for some $T > 0$ there exists a $(T, \delta_i)$-spanning subset $S \subset X$
with $|S| \leq 2^{(i+1)T} -1$.  Then there exists $x \in X$ such that for any $x' \in S$, there exists $0 \leq j \leq i$ with $y_j \restrict\{0, 1, \dots, T-1\} \neq y'_j \restrict\{0, 1, \dots, T-1\}$.  Then, if $k_j \neq k'_j$, clause (Dn1) applies in the definition of~$D^j(x, x')$ and we have
\begin{equation*}
\forall 0 \leq t \leq T-1 \ \ D^j(F^t_j(x_j), F^t_j(x'_j)) =  \delta_j \geq \delta_i.
\end{equation*}
If $k_j = k'_j$, then there exists $0 \leq t \leq T-1$ such that we get from clause~(Dn32):
\begin{equation*}
    D^j(F^t_j(x_j), F^t_j(x'_j)) \in \{\delta_j, \eps_j\} \geq \delta_j \geq \delta_i.
\end{equation*}
Therefore, $D_T(x, x') \geq D^j_T(x_j, x_j') \geq \delta_i$, which shows that $S$ is not a $(T, \delta_i)$-spanning set. This contradicts our assumption.

\medskip

For the proof of the first inequality in~\eqref{eqn:spanbounds-delta-i},
choose an odd number $m > 0$ large enough such that
\begin{equation}\label{eqn:mchoice}
(\eps_0+\eps_1 + \dots + \eps_i) 3^{-m} < \frac{3}{2}\delta_i - \frac{1}{2}\eps_i.
\end{equation}

Let $q \in \bZ$ be such that $\#(q) = m-1$. As~$m$ is odd, $q < 0$.
Then for any $T > 0$, we can choose a subset $S(T) \subset X$ of size
\begin{equation*}
|S(T)| = \left[\prod_{j=0}^i T^+(j)\right]2^{(i+1)(T+m-1)}
\end{equation*}

with the property that for all $x \in X$, there exists $x' \in S(T)$ such that

\begin{equation}\label{eqn:ks-eq}
(k'_j)_{j=0}^i = (k_j)_{j=0}^i
\end{equation}

and for all $0 \leq j \leq i$
\begin{equation}\label{eqn:ys-eq}
y'_j\restrict \{q, q+1, \dots, q+m+T-2\} = y_j \restrict \{q, q+1, \dots, q+m+T-2\}.
\end{equation}

Assume $x \in S(T)$ and $x' \in X$ are such that~\eqref{eqn:ks-eq} and~\eqref{eqn:ys-eq} hold.
Let $j \leq i$. Then~\eqref{eqn:ks-eq} implies that in the computation of~$D^j(F^t(x), F^t(x'))$ clause~(Dn1) will not be used for any~$t \in \bZ$. Similarly, in view
of~\eqref{eqn:ys-eq},
 for $t < T$ clause~(Dn32) will not apply either, and~$D^j(F^t(x), F^t(x')) \leq \eps_j 3^{-m}$ when clause~(Dn2) applies.

Therefore, for all $0 \leq t < T$,
\begin{equation*}
    \begin{split}
    D(F^t(x), F^t(x')) &= \sum_{j = 0}^i D^j(F^t_j(x_j), F^t_j(x'_j)) + \sum_{n = i+1}^{\infty} D^n(F^t_n(x_n), F^t_n(x'_n))\\
    &\leq \sum_{j = 0}^i \eps_j 3^{-m} + \sum_{n = i+1}^{\infty} \eps_n\\
    &< \sum_{j = 0}^i \eps_j 3^{-m} + \frac{1}{2}(\eps_i - \delta_i)\\
    &< \delta_i,
    \end{split}
\end{equation*}
where the second last inequality follows form~(P$\delta$2) and the last one from~\eqref{eqn:mchoice}.
Hence,
\begin{equation*}
    span(X, \delta_i, D_T) \leq |S(T)| = \left[\prod_{i=j}^i T^+(j)\right]2^{(i+1)(T+m-1)}.
\end{equation*}

$\Box$

\bigskip

The following result implies the first part of point~(iv) of Theorem~\ref{thm:main-seponly}.

\begin{corollary}\label{hinfty}
Let $(X,F)$ be an \JX-system with \JX-metric~$D$. Then
\begin{equation*}
    \lim_{T\rightarrow\infty} \frac{\ln{span(X, \delta_i, D_T)}}{T} = (i+1)\ln{2}\ \  \mbox{for all}\ i \in \bN.
\end{equation*}
In particular,  $h(X,  F) = \infty$.
\end{corollary}

\noindent
\textbf{Proof:} Fix $i \in \bN$. By Lemma~\ref{lem:span-lb},
\begin{equation*}
\forall\, T > 0  \quad span(X, \delta_i, D_T) \geq 2^{(i+1)T}.
\end{equation*}
Then
\begin{equation*}
\forall\,  T > 0 \quad    \frac{\ln{span(X, \delta_i, D_T)}}{T} \geq \frac{\ln{2^{(i+1)T}}}{T} = (i+1)\ln{2},
\end{equation*}
and hence
\begin{equation*}
    \liminf_{T\rightarrow\infty} \frac{\ln{span(X, \delta_i, D_T)}}{T} \geq (i+1)\ln{2}.
\end{equation*}

\medskip

Moreover, also by Lemma~\ref{lem:span-lb}, there exists $m \in \bN$ such that
\begin{equation*}
\forall\,  T > 0 \quad   span(X, \delta_i, D_T) \leq \left[\prod_{j=0}^i T^+(j)\right]2^{(i+1)(T+m-1)}.
\end{equation*}
Hence
\begin{equation*}
\forall\,  T > 0 \quad   \frac{\ln{span(X, \delta_i, D_T)}}{T} \leq \frac{\ln\left(\left[\prod_{j=0}^i T^+(j)\right]2^{(i+1)(T+m-1)}\right)}{T},
\end{equation*}
and
\begin{equation*}
    \limsup_{T\rightarrow\infty}\frac{\ln{span(X, \delta_i, D_T)}}{T} \leq \limsup_{T\rightarrow\infty}\frac{\ln \left(\left[\prod_{j=0}^i T^+(j)\right]2^{(i+1)(T+m-1)}\right)}{T} = (i+1)\ln{2}.
\end{equation*}

\medskip

It follows that
\begin{equation*}
    \lim_{T\rightarrow\infty} \frac{\ln{span(X, \delta_i, D_T)}}{T} = (i+1)\ln{2}.
\end{equation*}

\smallskip

Since $\lim_{i\rightarrow\infty} \delta_i = 0$, we have
$h(X,F) = \lim_{i\rightarrow \infty} (i+1)\ln{2} = \infty$.
$\Box$

\begin{lemma}\label{lem:h(Y)-for-Y-subspace-W}
Let $W$ be defined as in Subsection~\ref{subsec:WX-}, let~$D$ be an \JX-metric, let $Y \subseteq W$ be a closed subspace that is invariant under~$F$, and let
$\delta < \delta^* = \sum_{n\in\bN} \delta_n$.  Then
\begin{equation}\label{eqn:limspan-delta}
\begin{split}
    \lim_{T\rightarrow\infty} \frac{\ln{span(W, \delta, D_T)}}{T} &= \ln{2},\\
     \limsup_{T\rightarrow\infty} \frac{\ln{span(Y, \delta, D_T)}}{T} &\leq \ln{2}.
\end{split}
\end{equation}

In particular, $h(W, F) = \ln 2$ and $h(Y,F) \leq \ln 2$.
\end{lemma}

The first line of~\eqref{eqn:limspan-delta} give Theorem~\ref{thm:main-seponly}(iiib), and the last sentence of Lemma~\ref{lem:h(Y)-for-Y-subspace-W}
implies the second part of Theorem~\ref{thm:main-seponly}(iv).  Since~$X^-$ satisfies the assumptions on~$Y$ in this lemma in view of~\eqref{eqn:X-subW},   Theorem~\ref{thm:main}(iii) also follows.

\medskip

\noindent
\textbf{Proof of Lemma~\ref{lem:h(Y)-for-Y-subspace-W}:} First let us derive the last sentence from~\eqref{eqn:limspan-delta}.
By the definition of topological entropy,
\begin{equation*}
    \begin{split}
        h(W,F)& = \lim_{\delta \rightarrow 0^+}\limsup_{T\rightarrow\infty} \frac{\ln{span(W, \delta, D_T)}}{T}\\
        &= \lim_{\delta \rightarrow 0^+}\lim_{T\rightarrow\infty} \frac{\ln{span(W, \delta, D_T)}}{T}\\
        &= \lim_{\delta \rightarrow 0^+}\ln{2}\\
        &= \ln{2},\\
        h(Y,F)& = \lim_{\delta \rightarrow 0^+}\limsup_{T\rightarrow\infty} \frac{\ln{span(X^-, \delta, D_T)}}{T}\\
        &\leq \lim_{\delta \rightarrow 0^+}\ln{2}\\
        &= \ln{2}.
    \end{split}
\end{equation*}

\smallskip

Let $Y, D$ be as in the assumptions.
We first show that there exist $i, m \in \bN$  such that
\begin{equation}\label{eqn:span-upper-small-delta}
\forall T > 0  \quad span(Y, \delta, D_T) \leq  T^+(i) 2^{(T+m-1)}.
\end{equation}

We can choose $i, m$  with $m > 1$ odd  so  that
\begin{equation}\label{eqn:choose-im}
\begin{split}
\sum_{n=i+1}^{\infty}\eps_n &< \frac{\delta}{2},\\
\sum_{j=0}^i\eps_j 3^{-m} &< \frac{\delta}{2}.
\end{split}
\end{equation}

 Let $q \in \bZ$ be such that $\#(q) = m-1$.  Then   $q < 0$, as $m-1$ is even.

For any $T > 0$, we can choose a subset $S(T) \subset Y$ of size
\begin{equation*}
|S(T)| \leq T^+(i)2^{(T+m-1)}
\end{equation*}
with the property that for each $x \in Y$ there exists $x' \in S(T)$ such that

\begin{equation}\label{eqn:kj-eq}
k'_i = k_i
\end{equation}

and for all $0 \leq j \leq i$
\begin{equation}\label{eqn:ys-eq-copy}
y'_j\restrict \{q, q+1, \dots, q+m+T-2\} = y_j \restrict \{q, q+1, \dots, q+m+T-2\}.
\end{equation}

Here we use the assumption that $Y \subseteq W$ so that $y_j = y_i$ for all coordinates~$(y_j, j, k_j)$ of~$x$.

Assume $x \in S(T)$ and $x' \in Y$ are such that~\eqref{eqn:kj-eq} and~\eqref{eqn:ys-eq-copy} hold.
Let $j \leq i$. Since $Y \subseteq W$,  by~\eqref{eqn:kj-eq} and the definition of~$W$, also $k'_j =  k_j$ for all~$j \leq i$.
Thus in the computation of~$D^j(F^t(x), F^t(x'))$ clause~(Dn1) will not be used for any~$t \in \bZ$.
Similarly, in view of~\eqref{eqn:ys-eq-copy},
 for $t < T$ clause~(Dn32) will not apply either, and~$D^j(F^t(x), F^t(x')) \leq \eps_j 3^{-m}$ whenever clause~(Dn2) applies.

Therefore, for all $0 \leq t < T$ and $x, x'$ as above,
\begin{equation*}
    \begin{split}
    D(F^t(x), F^t(x')) &= \sum_{j = 0}^i D^j(F^t_j(x_j), F^t_j(x'_j)) + \sum_{n = i+1}^{\infty} D^n(F^t_n(x_n), F^t_n(x'_n))\\
    &\leq \sum_{j = 0}^i \eps_j 3^{-m} + \sum_{n = i+1}^{\infty} \eps_n\\
    &< \frac{\delta}{2} + \frac{\delta}{2}\\
    &= \delta,
    \end{split}
\end{equation*}
where the last inequality follows from~\eqref{eqn:choose-im}.
Hence,
\begin{equation*}
\forall T > 0 \quad   span(Y, \delta, D_T) \leq |S(T)| =  T^+(i) 2^{(T+m-1)}.
\end{equation*}

\medskip

Next we show that for any $T > 0$,
\begin{equation*}
    span(W, \delta, D_T) \geq 2^{T}.
\end{equation*}

\noindent
Assume towards a contradiction that for some $T > 0$, there exists a $(T, \delta)$-spanning subset $S \subset W$ with $|S| \leq 2^{T} -1$.  There exists $x \in W$ such that for any $x' \in S$ and $n \in \bN$,
we have  $y_n \restrict\{0, 1, \dots, T-1\} \neq y'_n \restrict\{0, 1, \dots, T-1\}$.
Then, if $k_n \neq k'_n$, clause (Dn1) applies in the definition of~$D^n(x, x')$ and we have
\begin{equation*}
\forall 0 \leq t \leq T-1 \ \ D^n(F^t_n(x_j), F^t_n(x'_n)) \geq \delta_n.
\end{equation*}
If $k_n = k'_n$, then there exists $0 \leq t \leq T-1$ such that we get from clause~(Dn32):
\begin{equation*}
    D^n(F^t_n(x_n), F^t_n(x'_n)) \in \{\delta_n, \eps_n\} \geq \delta_n.
\end{equation*}

Thus, there exists $0 \leq t \leq T-1$ such that
\begin{equation*}
    \begin{split}
        D(F^t(x), F^t(x')) &= \sum_{n = 0}^{\infty} D^n(F^t_n(x_n), F^t_n(x'_n))\\
        &\geq \sum_{n = 0}^{\infty} \delta_n\\
        &\geq \delta.
    \end{split}
\end{equation*}

Therefore, $D_T(x, x') \geq \delta$, which indicates that $S$ is not a $(T, \delta)$-spanning set, and this contradicts our assumption.

\medskip

We have shown that there exist $i \in \bN$ and $m \in N$ such that for all $T > 0$,
\begin{equation*}
\begin{split}
    2^T \leq span(W, \delta, D_T) &\leq T^+(i)2^{(T+m-1)},\\
   span(Y, \delta, D_T) &\leq T^+(i)2^{(T+m-1)},\\
   \ln{2} \leq \frac{\ln{span(W, \delta, D_T)}}{T} &\leq \frac{\ln T^+(j)}{T} + \left(\frac{T+m-1}{T}\right)\ln{2},\\
   \frac{\ln{span(Y, \delta, D_T)}}{T} &\leq \frac{\ln T^+(j)}{T} + \left(\frac{T+m-1}{T}\right)\ln{2},
\end{split}
\end{equation*}
and~\eqref{eqn:limspan-delta} follows. $\Box$

\bigskip

\subsection{Topological transitivity}\label{subsec:MoreTrans}

Recall that a dynamical system is \emph{topologically transitive} if there exists a dense forward orbit, or, equivalently, if for every nonempty open $U, V$ there exists $t \geq 0$ such that $V \cap F^t(U) \neq \emptyset$. The following result proves  part~(v) of
Theorem~\ref{thm:main-seponly}.

\begin{proposition}\label{prop:transitive}
Consider the \JX-system $(X,F)$ with any \JX-metric $D$, and let $W \subset X$ be as defined in
Subsection~\ref{subsec:WX-}. Then

\begin{itemize}
\item[(a)] The system $(X, F)$ is not topologically transitive.
\item[(b)] The system $(W, F \restrict W)$ is topologically transitive.
\end{itemize}
\end{proposition}

\noindent
\textbf{Proof:} For the proof of part~(a), consider $x, x' \in X$ with $k_0 = 0 = k_1$ and $k_0' = 0 \neq 1 = k_1'$.
Let
$V$ and $U$ denote the open balls of radius $\frac{\delta_1}{2}$ with centers $x, x'$, respectively.

Then for every
$t \in \bZ$ the first two coordinates of  $F^t(x')$ will be of the form $(\sigma^t(y_0'), 0, (t \mod T^+(0)))$ and
$(\sigma^t(y_1'), 1, (t + 1 \mod T^+(1)))$.  Since $T^+(1)$ is an integer multiple of~$T^+(0)$,  clause~(Dn1) of the definition of the metric~$D^n$ implies that
for all~$t$ we will have either $D^0(x_0, F_0^t(x_0')) \geq \delta_0$ or $D^1(x_1, F_1^t(x_1')) \geq \delta_1$. In both cases $D(x, F^t(x')) \geq \delta_1$.
Let $x'' \in U$. Then we must have $k_0'' = 0 \neq 1 = k_1''$,
and the same argument shows that  $D(x, F^t(x'')) \geq \delta_1$ for all $t \in \bZ$.

Thus  we will have
$V \cap F^t(U) = \emptyset$ for all $t \in \bZ$.

\smallskip

For the proof of part~(b), let us arrange into a sequence~$(Q^\ell)_{\ell \in \bN}$ all quadruples  of the form $Q = (m, \kappa, k, N)$, where
$\kappa \in {}^{[-m, m]}\{0, 1\}$ and $0 \leq k < T^+(N)$ for some~$m,N \in \bN$ with $m > T^+(N)$. Let $Q^\ell = (m^\ell, \kappa^\ell, k^\ell, N^\ell)$. Now we construct recursively a $y^* \in {}^\bZ \{0,1\}$ as follows: For every~$\ell \in \bN$, we pick $t(\ell)$ such that\\
 $t(\ell)  \mod T^+(N^\ell) = k^\ell$ and also pick an interval of positive integers\\
$J_\ell = [t(\ell) - m^\ell, \dots , t(\ell), \dots , t(\ell) + m^\ell]$.
 We choose these objects so that the intervals~$J_\ell$ will be pairwise disjoint. Then we choose~$y^*$ in such a way that $y^*(t(\ell) + i) = \kappa^\ell(i)$ for
all $-m^\ell \leq i \leq m^\ell$.  Finally, we let $x^* \in X$ be such that for all~$n$ we have $x^*_n = (y^*, n, 0)$.  Then $x^* \in W$.

Now consider any $x \in W$ and let $V$ be an open ball with center~$x$ and radius~$\gamma$ for some $\gamma > 0$.  Consider any  $x' \in W$.
Let $N \in \bN$ be such that
$\sum_{n = N+1}^\infty \eps_n < \frac{\gamma}{2}$.  Then $x' \in V$ whenever
\begin{equation}\label{eqn:x'VN}
\sum_{n = 0}^N D^n(x_n, x_n') < \frac{\gamma}{2}.
\end{equation}

Recall that in view of the definition of~$W$, there are~$y, y'$ such that $x_n = (y, n, k_n)$ and
$x'_n = (y', n, k_n')$ for all coordinates of~$x$ and~$x'$.

Now let us assume that  $k_N = k_N'$ and $n \leq N$. Then $k_n = k_n'$  by the definition of~$W$, so that~$D^n(x_n, x_n')$ will be determined by clauses~(Dn2) or~(Dn3).
By clause~(Dn31),   whenever  the restriction of~$y'$ to an
interval~$[-m, m] = \{-m, -m + 1, \dots , m -1, m\}$ is the same as the restriction of~$y$ to this interval so that $\Delta(y, y') > m$, then
$D^n(x_n, x'_n) \leq \eps_n 3^{-m}$.  If such~$m$ is chosen sufficiently large, the inequality~\eqref{eqn:x'VN} follows, and $x' \in V$.

It remains to show that for some~$t \geq 0$ the point $x' = F^t(x^*)$ will have the properties outlined above. Choose $N$ and then~$m$ sufficiently large
so that $m > T^+(N)$. Consider the quadruple~$Q = (m, \kappa, k_N, N)$, where $\kappa = y \restrict [-m , m]$. Then $Q = Q^\ell$ for some~$\ell$.
Let $x' = F^{t(\ell)}(x^*)$ with coordinates $x'_n = (\sigma^{t(\ell)}(y^*), n, k'_n)$. Then
\begin{equation*}\label{eqn:Qk}
k'_{N} \ = \  t(\ell) \! \! \mod T^+(N) \ = \ t(\ell) \!\! \mod T^+(N^\ell) \ = \ k^\ell \ = \ k_N,
\end{equation*}
and since~$x' \in W$, we also have $k_n' = k_n$ for all $n \leq N$.

Moreover, by the definition of~$F$ and our choice $y^*\restrict J_\ell = \kappa^\ell$, we also have
$\sigma^{t(\ell)}(y^*) \restrict [-m, m] = y\restrict [-m, m]$.  Thus by the choice of~$m, N$ we must have $x' = F^{t(\ell)}(x^*) \in V$.

We have shown that the forward orbit of~$x^*$ under~$F$ is dense in~$W$, and topological transitivity of~$(W, F\restrict W)$ follows.
$\Box$

\bigskip

It appears that the system $(X^-, F)$ of Theorem~\ref{thm:main} is not topologically transitive. However, it seems likely that it can be modified into a transitive system that still satisfies parts~(i)--(iii) of the theorem.  We will return to this issue in~\cite{PartII}.

\section{Bounds on separation and spanning numbers}\label{sec:sepsizes}

In this section we will derive lower bounds on $sep\left(W, \eps, D_{T^+(n)}\right)$, and then $span\left(X^-, \eps, D_{T^+(n)}\right)$, as well as upper bounds on
$sep\left(X, \eps, D_{2T(n)}\right)$.

\subsection{Lower bounds on $sep\left(W, \eps, D_{T^+(n)}\right)$ and $sep\left(X^-, \eps, D_{T^+(n)}\right)$}\label{subsec:lower-sep-T+n}

\begin{lemma}\label{prop:separated}
Let $(X, F)$ be an \JX-system with \JX-metric~$D$ that is based on condition~(Dn32c) with all
colorings~$c_n$ either satisfying conditions (cC1), (cC2) or conditions (cCi), (cC).
Let $n \in \bN$, and  let $\tau$ be any non-negative integer multiple of $T^+(n)$.  Suppose~$y,z \in {}^{\bZ} \{0,1\}$.
Let  $\varphi := \sigma^\tau (y) \restrict (0, \dots , T^+(n) -1)$ and \\ $\psi :=
\sigma^\tau (z) \restrict (0, \dots , T^+(n) -1)$ be such that

\begin{itemize}
\item[(P2-1)]   $|C(\varphi, \psi)| > 2$ and\\
 there exist $T^+(n-1) \leq i, j \leq T^+(n)-1$ such that $\varphi(i) = \psi(j) = 1$.
\end{itemize}

\smallskip

Then
\smallskip

\begin{itemize}
\item[(i)] There exists $0 \leq t \leq T^+(n)-1$ such that \newline $D^n_{T^+(n)}(F^{\tau}_n(y,n,0),F^{\tau}_n(z,n,0)) =  D^n (F_n^{t+\tau}(y, n, 0), F_n^{t+\tau}(z, n, 0)) = \eps_n$.

\smallskip

\item[(ii)] More precisely, there exists exactly one color~$j$ such that:

$D^n(F_n^{t+\tau}((y, n, 0)), F_n^{t+\tau}((z, n, 0))) = \eps_n$ for all $t \in I^n_j$ \newline
with $y(t+\tau) \neq z(t+\tau)$, and this inequality will hold for at least one $t \in I^n_j$.

\smallskip

\item[(iii)] Suppose  $\varphi \neq \psi$
are both elements  of~$\Ymn$ so that, in particular, (P2-1) holds.  Then there exists $0 \leq t < T^+(n)$ such that\\
$D^m(F_m^{t+\tau}((y, m, 0)), F_m^{t+\tau}((z, m, 0))) = \eps_m$ for all $0 \leq m \leq n$.
\end{itemize}
\end{lemma}

\noindent
\textbf{Proof:}   Let $y, z, \tau, \varphi, \psi$ be as in the assumptions, and let
$C(\varphi, \psi) = \{j_1,j_2,\dots, j_{\ell}\}$. Assume $|C(\varphi,\psi)| = \ell \geq 3$, so that, in particular,
$\varphi \neq \psi$.

\medskip
\noindent
(i) By the assumption (P2-1) of the lemma, condition (cC2) or (cC) applies and entails that $c_n(\varphi,\psi) \in C(\varphi,\psi)$.

To be specific,
let $c_n(\varphi,\psi) = j_r \in C(\varphi,\psi)$.  In view of the definition of~$C(\varphi, \psi)$
we must have $\varphi\restrict I_{j_r}^n \neq \psi\restrict I_{j_r}^n$,
so that we can choose $t \in I_{j_r}^n$ with $\sigma^\tau(y)(t) = y(t+\tau) \neq z(t+\tau) = \sigma^\tau(z)(t)$.
Then $\sigma^{t+\tau}(y)(0) \neq \sigma^{t+\tau}(z)(0)$.  Since $0 \leq t < T^+(n)$, it follows from our choice of $\varphi, \psi$,
from Proposition~\ref{prop:Phi-const}, and from the choice of~$\tau$ as an integer multiple of~$T^+(n)$ that $\Phi(F^{t+\tau}_n((y,n,0))) = \Phi((\sigma^{t+\tau}(y),n,t)) = \varphi$ and
$\Phi(F^{t+\tau}_n((z,n,0))) = \Phi((\sigma^{t+\tau}(z),n,t)) = \psi$.
Then by clause~(Dn32c) of the definition of $D^n$, we have:
\begin{equation*}
\begin{split}
    &D^n(F^{t+\tau}_n((y,n,0)), F^{t+\tau}_n((z,n,0)))\\
    &= D^n((\sigma^{t+\tau}(y),n,t),(\sigma^{t+\tau}(z),n,t))\\
    &= \eps_n.
\end{split}
\end{equation*}

\medskip

\noindent
(ii)  By Proposition~\ref{prop:Phi-const}, $\Phi(F^{t+\tau}_n((y,n,0))) = \Phi((\sigma^{t+\tau}(y),n,t)) = \varphi$ and\newline
$\Phi(F^{t+\tau}_n((z,n,0))) = \Phi((\sigma^{t+\tau}(z),n,t)) = \psi$
for all $0 \leq t < T^+(n)$.\newline
Thus  whenever clause~(Dn32c) applies, the value $c_n(\varphi,\psi) = j_r$ is used in the computation of
$D^n(F^{t+\tau}_n((y,n,0)), F^{t+\tau}_n((z,n,0)))$.  By our argument for part~(i), this will happen at least once
for~$t \in I^n_{j_r}$, and the same argument shows that $D^n(F^{t+\tau}_n((y,n,0)), F^{t+\tau}_n((z,n,0))) = \eps_n$
for all~$t \in I^n_{j_r}$ with $y(t+\tau) \neq z(t+\tau)$.

Conversely, if for some $j \neq j_r$ and $t \in I^n_j$ we have $y(t+\tau) \neq z(t+\tau)$, then clause~(Dn32c) will be applied in
the computation of \\$D^n(F^{t+\tau}_n((y,n,0)), F^{t+\tau}_n((z,n,0))) = D^n((\sigma^{t+\tau}(y),n,t),(\sigma^{t+\tau}(z),n,t))$, \newline and by Proposition~\ref{prop:Phi-const} again, \newline
$ c_n(\Phi(F^{t+\tau}_n((y,n,0))), \Phi(F^{t+\tau}_n((z,n,0)))) = c_n(\varphi,\psi) = j_r \neq j$. Thus in this case \newline
$D^n(F^{t+\tau}_n((y,n,0)), F^{t+\tau}_n((z,n,0))) = \delta_n \neq \eps_n$.

\medskip

\noindent
(iii)  Notice that if $\varphi \neq \psi$
are both elements  of~$\Ymn$, then  (P2-1) holds as a consequence of~(PY2+) and either (PY1) (for $n= 0$) or (PR2) (for $n > 0$).

Let us temporarily assume that $n$ is fixed.  Let (iii)[n] denote the assertion of part (iii) for this particular~$n$.  We prove by induction over~$n$ that (iii)[n] holds.

Note that (iii)[0] is simply an instance of~(i).

Now assume $n > 0$ and that (iii)[n-1] holds. Let $y, z, \tau, \varphi, \psi$ be as specified in the first paragraph of this proof of
Lemma~\ref{prop:separated}. Moreover, assume that  $\varphi, \psi \in \Ymn$. Let ~$j$ be such that the conclusion of part~(ii) holds.

Recall from~\eqref{eqn:Injs} and~\eqref{eqn:TT+(n)} that $I^n_j \subset (0, T^+(n)-1)$ is an interval of
length~$T(n) = K(n)T^+(n-1)$ that starts at an integer multiple of~$T(n)$. Thus by Definition~\ref{def:Yn-}, the restrictions of~$\varphi$ and~$\psi$ to~$I^n_j$ consist of $K(n)$
consecutive blocks of length~$T^+(n-1)$, each of which is in~$\Ymnm$.
Moreover, by the last sentence of~(ii), there must be some such block, call its domain~$B$, with the property that for some~$t \in B$ we have $\varphi(t) \neq \psi(t)$.  Also, since $B \subset I_j^n$, by the first part of point~(ii), we will have
\begin{equation}\label{eqn:all-t-B}
\forall t \in B \ (\varphi(t) \neq \psi(t) \ \Rightarrow \ D^n(F_n^{t+\tau}((y, n, 0)), F_n^{t+\tau}((z, n, 0))) = \eps_n).
\end{equation}

Now let $B = (t_0, t_0+1, \dots , t_0 + T^+(n-1)-1)$. Since the lengths of the intervals~$I^n_i$ and of the aforementioned blocks inside~$I^n_j$ are all
integer multiples of~$T^+(n-1)$, $t_0$ is an integer multiple of
$T^+(n-1)$. Thus it follows from the definition of~$F_{n-1}$ that $F^{t_0+\tau}_{n-1}((y,n-1,0)) = (\sigma^{t_0+\tau}(y),n-1,0)$ and \newline
$F^{t_0+\tau}_{n-1}((z,n-1,0)) = (\sigma^{t_0+\tau}(z),n-1,0)$. Let \newline
$\varphi^- = \sigma^{t_0+\tau}(y) \restrict \{0, \dots , T^+(n-1)-1\}$
and
$\psi^- = \sigma^{t_0+\tau}(z) \restrict \{0, \dots , T^+(n-1)-1\}$. \newline
As we already mentioned,
$\varphi^-, \psi^- \in \Ymnm$.  As we have chosen $B$ in such a way that~$\varphi^- \neq \psi^-$, we can infer from Definition~\ref{def:Yn-} that
$|C(\varphi^-, \psi^-)| > 2$, where~$C(\varphi^-, \psi^-)$ is now a subset of~$[C(n-1)]$ rather than of~$[C(n)]$. Then by properties~(PY2+) and either (PY1) (for $n -1 = 0$) or (PR2) (for $n -1 > 0$) the assumption~(P2-1) of Lemma~\ref{prop:separated} is satisfied if we substitute~$n-1$ for~$n$
and $(t_0+\tau)$ for ~$\tau$.
Thus by the inductive assumption there exists $0 \leq t_1 < T^+(n-1)$ such that
\begin{equation}\label{eqn:(iii)-}
D^m(F_m^{t_1+t_0+\tau}((y, m, 0)), F_m^{t_1+t_0+\tau}((z, m, 0))) = \eps_m \ \ \mbox{for all} \ \ 0 \leq m \leq n-1.
\end{equation}

Let $t = t_0 + t_1$.  Then $t \in I^n_j$, and it follows from~\eqref{eqn:(iii)-} that $y(t+\tau) \neq z(t+\tau)$. Now it follows from~\eqref{eqn:all-t-B} and~\eqref{eqn:(iii)-} that
\begin{equation*}
D^m(F_m^{t+\tau}((y, m, 0)), F_m^{t+\tau}((z, m, 0))) = \eps_m \qquad \mbox{for all} \ \ 0 \leq m \leq n.
\end{equation*}
Thus (iii)[n] holds, and point (iii) follows by mathematical induction. $\Box$

\begin{corollary}\label{corol:large-sep}
Let $(X, F)$ be an \JX-system with \JX-metric~$D$ that is based on condition~(Dn32c) with all
colorings~$c_n$ either satisfying conditions (cC1), (cC2) or conditions (cCi), (cC).
Then for each $n \in \bN$, the set~$W^n$ is $(T^+(n),\eps)$-separated.
\end{corollary}

\noindent
\textbf{Proof:} Recall the definition of~$y_\varphi \in {}^\bZ\{0,1\}$ for~$\varphi \in {}^{T^+(n)}\{0, 1\}$:

$y_{\varphi}\restrict(0,\dots,T^+(n)-1) = \varphi$ and $y_{\varphi}(i) = 0$ when $i \geq T^+(n)$ or $i \leq -1$.

Let~$n \in \bN$ and consider~$\varphi \neq \psi \in \Ymn$.
The assumptions of the above lemma are satisfied for
$y = y_\varphi, z = y_\psi$, and~$\tau = 0$. Thus by part ~(iii) of the lemma we can find $t < T^+(n)$ so that
\begin{equation}\label{eqn:(iii)-repeat}
\forall 0 \leq m \leq n \ \ D^m(F^t_m((y_{\varphi},m,0)),F^t_m((y_{\psi},m,0))) = \eps_m.
\end{equation}

When $m > n$, then
$c_m(y_{\varphi}\restrict(0,\dots,T^+(m)),y_{\psi}\restrict(0,\dots,T^+(m)))$ is constrained by
either condition~(cC1) or condition~(cCi) to take the value~1. Clause~(Dn32c) will be used in the computation of
$D^m(F^t_m((y_{\varphi},m,0)),F^t_m((y_{\psi},m,0))) $, which will evaluate to~$\eps_m$ as  $t \in I^m_1$.
It follows that
\begin{equation*}
\forall m \in \bN \quad D^m(F^t_m((y_{\varphi},m,0)),F^t_m((y_{\psi},m,0))) = \eps_m.
\end{equation*}

Then for any such~$t$ we must have
\begin{equation*}
    \begin{split}
        D(F^t(x^{\varphi}),F^t(x^{\psi})) &= \sum_{m\in\bN} D^m(F^t_m(x^{\varphi}_m), F^t_m(x^{\psi}_m))\\
        &= \sum_{m\in\bN} D^m(F^t_m((y_{\varphi},m,0)), F^t_m((y_{\psi},m,0)))\\
        &= \sum_{m\in\bN} \eps_m\\
        &= \eps.
    \end{split}
\end{equation*}

Since all elements of~$W^n$ are of the form $x^{\varphi}$ for some~$\varphi$ as above, it follows that each set~$W^n$ is $(T^+(n),\eps)$-separated. $\Box$

\bigskip

As   $W^n \subset X^- \subset W$
by~\eqref{eqn:X-subW}, Corollary~\ref{corol:size}  implies that we have

\begin{corollary}\label{corol:+sep}
Let $n \in \bN$.  Then for the relevant \JX-metrics $D$:
\begin{equation*}
\begin{split}
sep(W, \eps, D_{T^+(n)})\  &\geq \ |W^n| \ \geq \ 2^{0.9T^+(n)},\\
 sep(X^-, \eps, D_{T^+(n)})\  &\geq \ |W^n| \ \geq \ 2^{0.9T^+(n)}.
\end{split}
\end{equation*}
\end{corollary}

\subsection{Lower bounds on  $span\left(X^-, \eps, D_{T^+(n)}\right)$}\label{subsec:lower-span-T+n}

\begin{lemma}\label{corol:low-span-deleted}
Fix any $n \in \bN$.  Then for all $x \in X^-$ and all $u \neq v \in W^n$,
\begin{equation*}
    \max\{D_{T^+(n)}(x,u), D_{T^+(n)}(x,v)\} = \eps.
\end{equation*}
\end{lemma}

\noindent
\textbf{Proof:}
Fix any $n \in \bN$ and $u \neq v \in W^n$.  For $x \in X^-$, we distinguish the following seven cases:
\begin{itemize}
\item[Case 1:] $x \in W^n$.
\item[Case 2:] $x \in W^m$ for some $0 \leq m < n$ (if $n \neq 0$).
\item[Case 3:] $x \in W^m$ for some $m > n$.
\item[Case 4:] $x \in F^{\tau}(W^m)$ for some $\tau > 0$ and $0 \leq m \leq n$.
\item[Case 5:] $x \in F^{\tau}(W^m)$ for some $\tau > 0$ and $m > n$.
\item[Case 6:] $x \in F^{\tau}(W^m)$ for some $\tau < 0$ and $m \in \bN$.
\item[Case 7:] $x \in X^-\backslash\left[\bigcup_{t \in \bZ}F^t(\bigcup_{n\geq 0}W^n)\right]$.
\end{itemize}

\medskip

Let us remark from the outset that Case~7 is not redundant, as $(X^-, F)$ has positive entropy while
$\bigcup_{t \in \bZ}F^t\left(\bigcup_{n\geq 0}W^n\right)$ is countable. We will reduce this case to the preceding ones
by a density argument. Case~1 is the base case.

\medskip

\noindent
Case 1: $x \in W^n$.

\smallskip

Since $u \neq v$, we may wlog assume that $x \neq u$.

By Corollary~\ref{corol:large-sep},  $W^n$ is $(T^+(n), \eps)$-separated.  Thus, $D_{T^+(n)}(x, u) \geq \eps$.

By~\eqref{eqn:large-dist}, we have $D_{T^+(n)}(x, u),\ D_{T^+(n)}(x, v) \leq \eps$.
Hence,
\begin{equation*}
    \max\{D_{T^+(n)}(x,u), D_{T^+(n)}(x,v)\} = \eps.
\end{equation*}

\medskip

In most of the Cases~2--6 and their subcases, we will be able to show that $D_{T^+(n)}(x,u) = \eps$; in some cases
we will need to rely on the inequality $u \neq v$ and pick the element of $\{u,v\}$ that has distance~$\eps$ from~$x$.
Without loss of generality we will name this element~$u$. We then consider functions
$y_x, y_u \in {}^\bZ\{0,1\}$ such that for all $s \geq 0$ the coordinates of~$x$ and~$u$ will be of the form
\begin{equation}\label{eqn:xs-us-general}
    x_s = (y_x, s, k_s),\ \ \ \ u_s = (y_u, s, 0).
\end{equation}

The detailed arguments for particular cases will use slightly different notations for separate instances
of~\eqref{eqn:xs-us-general}. In particular, we will use~$\varphi^-, \varphi, \varphi^+$ as subscripts for~$y$
to suggest the relative magnitudes of~$m$ and $n$ when
$x$ is a shifted version of an element of~$W^m$ and $u \in  W^n$.

The proof then boils down to finding $t$ with $0 \leq t < T^+(n)$ such that
\begin{equation}\label{eqn:Ds-large}
D^s(F^t_s(x_s), F^t_s(u_s)) =  D^s((\sigma^t(y_x), s, F_s^t(k_s)), (\sigma^t(y_u), s, F_s^t(0))= \eps_s,
\end{equation}
where $F_s^t(k_s) = k_s + t \mod T^+(s)$ and  $F_s^t(0) = t \mod T^+(s)$. Depending on the particular (sub)case, this will be accomplished
by relying on condition (PY2+) or on Lemma~\ref{prop:separated}(iii).

Recall from~\eqref{eqn:large-dist} that~\eqref{eqn:Ds-large} must hold for all~$s \in \bN$ simultaneously for the
same~$t$ so that we can deduce~$D_{T^+(n)}(x,u) = \eps$.
Once a suitable~$t$ is identified, we will then derive~\eqref{eqn:Ds-large}, by arguments that may be different for different values of~$s$.
When $F_s^t(k_s) = k_s + t \mod T^+(s) \neq F_s^t(0) = t \mod T^+(s)$, then we can simply invoke  clause (Dn1e) of the definition of~$D^s$. If not,
then~\eqref{eqn:Ds-large} will follow either from Lemma~\ref{prop:separated}(iii), or from the following observation:

\begin{proposition}\label{prop:cCiprop}
Let $x_s, u_s$ be as in~\eqref{eqn:xs-us-general} and let $t \in I_1^s$. Assume that
$F^{t}_s(x_s) = (\sigma^{t} (y_x), s, t)$ and $F^{t}_s(u_s) = (\sigma^{t} (y_u), s, t)$,
with~$\Delta(\sigma^{t} (y_x), \sigma^{t} (y_u)) = 0$.  If for at least one~$w \in \{x, u\}$ we have
$\Phi(\sigma^{t}(y_w), s, t)(i) = 0$ for all $T^+(s-1) \leq i \leq T^+(s)-1$, then~\eqref{eqn:Ds-large} holds.
\end{proposition}

\noindent
\textbf{Proof:} Under the assumptions of the proposition, $D^s(F^t_s(x_s), F^t_s(u_s))$ will be computed according to clause~(Dn32c), for a coloring that
satisfies condition~$(cCi)$, so that $c_n(\Phi(\sigma^{t}(y_x), s, t),\Phi(\sigma^{t}(y_u), s, t)) = 1$. $\Box$

\newpage

\noindent
Case 2: $x \in W^m$ for some $0 \leq m < n$ (if $n \neq 0$).

\smallskip

It suffices to prove that $D_{T^+(n)}(x, u) = \eps$.

There exist $\varphi \in \mathcal{Y}_m^-$ and $\varphi^+ \in \Ymn$ such that for all $s \geq 0$, we have
\begin{equation*}
    x_s = (y_{\varphi}, s, 0),\ \ \ \ u_s = (y_{\varphi^+}, s, 0).
\end{equation*}

Choose the smallest $T^+(m) \leq t \leq T^+(m)+T(0)-1 < T^+(n)-1$ such that

$y_{\varphi^+}(t) =1$.  The
existence of such a $t$ follows from (PY2+).

Our goal is to show that $D^s(F^t_s(x_s), F^t_s(u_s)) = \eps_s$ for all $s \in \bN$.
For $s \leq m$,
\begin{equation*}
    F^t_s(x_s) = (\sigma^t(y_{\varphi}), s, t - T^+(m)),\ \ \ \ F^t_s(u_s) = (\sigma^t(y_{\varphi^+}), s, t - T^+(m)).
\end{equation*}
\begin{equation*}
    \sigma^t(y_{\varphi})(0) = 0,\ \ \ \ \sigma^t(y_{\varphi^+})(0) = 1,
\end{equation*}

Hence $\Delta(\sigma^t(y_{\varphi}), \sigma^t(y_{\varphi^+})) = 0$.

Moreover, $\Phi(\sigma^t(y_{\varphi}), s, t-T^+(m)) = {}^{T^+(s)}\{0\}$.

Thus Proposition~\ref{prop:cCiprop} applies, and we get
$D^s(F^t_s(x_s), F^t_s(u_s)) = \eps_s$.

For $s > m$,
\begin{equation*}
    F^t_s(x_s) = (\sigma^t(y_{\varphi}), s, t),\ \ \ \ F^t_s(u_s) = (\sigma^t(y_{\varphi^+}), s, t).
\end{equation*}

Here $\Phi(\sigma^t(y_{\varphi}), s, t)\restrict (0, \dots, T^+(m)-1) = \varphi$
and $\Phi(\sigma^t(y_{\varphi}), s, t)(i) = 0$ for all

$T^+(m) \leq T^+(s-1) \leq i \leq T^+(s)-1$.

Thus Proposition~\ref{prop:cCiprop} applies again, and we also get
$D^s(F^t_s(x_s), F^t_s(u_s)) = \eps_s$.

\medskip

\noindent
Case 3: $x \in W^m$ for some $m > n$.

\smallskip

There exist $\varphi^+ \in Y_m^-$ and  $\varphi, \varphi' \in \Ymn$ such that for all $s \geq 0$, we have
\begin{equation*}
    x_s = (y_{\varphi^+}, s, 0),\ \ \ \ u_s = (y_{\varphi}, s, 0),\ \ \ \ v_s = (y_{\varphi'},s,0).
\end{equation*}

As $u \neq v$, we have $\varphi \neq \varphi'$.  Wlog, we can assume that

$\varphi^+\restrict (0,\dots, T^+(n)-1) \neq \varphi$.

Now it suffices to show that $D_{T^+(n)}(x,u) = \eps$.
From the
the fact that $\varphi^+ \in Y_m^-$

it follows from property~(PR1) by induction
that $\varphi^+\restrict (0,\dots, T^+(n)-1)$

and  $\varphi$ are both elements of $\Ymn$.

Thus, by Lemma~\ref{prop:separated}(iii), there exists $t < T^+(n)$ such that

$D^s((F_s^t(y_{\varphi^+}, s, 0), F_s^t(y_{\varphi}, s, 0)) = \eps_s$ for all $0 \leq s \leq n$.

This implies that $\sigma^t(y_{\varphi^+})(0) \neq \sigma^t(y_{\varphi})(0)$ for such a $t$.

Fix $t$ as above.
We still need to show that $D^s(F^t_s(y_{\varphi^+}, s,0), F^t_s(y_{\varphi}, s, 0)) = \eps_s$

for all $s > n$.
Given any $s > n$,
\begin{equation*}
    F^t_s(y_{\varphi^+}, s,0) = (\sigma^t(y_{\varphi^+}), s, t),\ \ \ \ F^t_s(y_{\varphi}, s,0) = (\sigma^t(y_{\varphi}), s, t).
\end{equation*}

Then $\sigma^t(y_{\varphi^+})(0) \neq \sigma^t(y_{\varphi})(0)$ and hence  $\Delta(\sigma^t(y_{\varphi^+}), \sigma^t(y_{\varphi})) = 0$.

Moreover, $\Phi(\sigma^t(y_{\varphi}), s, t)\restrict (0,\dots, T^+(n)-1) = \varphi$ and $\Phi(\sigma^t(y_{\varphi}), s, t)(i) = 0$

for all $T^+(n) \leq T^+(s-1) \leq i \leq T^+(s)-1$.  Thus Proposition~\ref{prop:cCiprop} applies,

and we  get
$D^s(F^t_s(y_{\varphi^+}, s,0), F^t_s(y_{\varphi}, s, 0)) = \eps_s$.

\medskip

\noindent
Case 4: $x \in F^{\tau}(W^m)$ for some $\tau > 0$ and $0 \leq m \leq n$.

\smallskip

There exist $\varphi^- \in \mathcal{Y}_m^-$ and  $\varphi \in \Ymn$ such that for all $s \geq 0$, we have
\begin{equation*}
    x_s = (\sigma^{\tau}(y_{\varphi^-}), s, \tau\ \mbox{mod}\ T^+(s)),\ \ \ \ u_s = (y_{\varphi}, s, 0).
\end{equation*}

In this case, we distinguish the following three subcases for the value of $\tau$:

\smallskip

Case 4-1:  $\tau \geq T^+(m)$.

Case 4-2:  $0 < \tau < T^+(m)$ is not a positive integer multiple of $T^+(0)$.

Case 4-3:  $0 < \tau < T^+(m)$ is a positive integer multiple of $T^+(0)$.

\smallskip

We will show that in each of of these subcases $D_{T^+(n)}(x,u) = \eps$.

\medskip

\noindent
Case 4-1: $\tau \geq T^+(m)$.

\smallskip

In this case, $\sigma^{\tau}(y_{\varphi^-})(i) = 0$ for all $s \in \bN$ and $i \geq 0$. Choose the smallest

$0 \leq t \leq T(0)-1$ with $\varphi(t) = 1$.  The existence of such a $t$ follows from (PY2+).

We are going to show that $D^s(F^t_s(y_{\varphi},s,0), F^t_s(\sigma^{\tau}(y_{\varphi^-}),s, \tau\ \mbox{mod}\ T^+(s))) = \eps_s$

for all $s \in \bN$.  Here,
\begin{equation*}
\begin{split}
    F^t_s(y_{\varphi},s,0) &= (\sigma^t(y_{\varphi}), s, t),\\ F^t_s(\sigma^{\tau}(y_{\varphi^-}),s, \tau\ \mbox{mod}\ T^+(s)) &= (\sigma^{t+\tau}(y_{\varphi^-}), s, (t+\tau)\ \mbox{mod}\ T^+(s)).
\end{split}
\end{equation*}

If $t \neq (t+\tau)\ \mbox{mod}\ T^+(s)$,

then $D^s(F^t_s(y_{\varphi},s,0), F^t_s(\sigma^{\tau}(y_{\varphi^-}),s, \tau\ \mbox{mod}\ T^+(s))) = \eps_s$ by (Dn1e).

If $t = (t+\tau)\ \mbox{mod}\ T^+(s)$, with $\sigma^t(y_{\varphi})(0) = 1 \neq 0 = \sigma^{t+\tau}(y_{\varphi^-})(0)$, we have

$\Delta(\sigma^t(y_{\varphi}), \sigma^{t+\tau}(y_{\varphi^-})) = 0$ and
$\Phi(\sigma^{t+\tau}(y_{\varphi^-}), s, t) = {}^{T^+(s)}\{0\}$.

Then Proposition~\ref{prop:cCiprop} implies
$D^s(F^t_s(y_{\varphi},s,0), F^t_s(\sigma^{\tau}(y_{\varphi^-}),s, \tau\ \mbox{mod}\ T^+(s))) = \eps_s$.

\medskip

\noindent
Case 4-2: $0 < \tau < T^+(m)$ is not a positive integer multiple of $T^+(0)$.

\smallskip

In this case, $\tau\ \mbox{mod}\ T^+(s) \neq 0$ for all $s \in \bN$.  Then by (Dn1e), we have

$D^s((\sigma^{\tau}(y_{\varphi^-}), s, \tau\ \mbox{mod}\ T^+(s)), (y_{\varphi}, s, 0)) = \eps_s$ for all $s \in \bN$.  Thus, $D(x, u) = \eps$.

\medskip

\noindent
Case 4-3:  $0 < \tau < T^+(m)$ is a positive integer multiple of $T^+(0)$.

\smallskip

Choose the smallest $T^+(m)-\tau \leq t \leq T^+(m)-\tau+T(0)-1$ with $\varphi(t) \neq 0$.

The existence of such a $t$ follows from (PY2+).

We will show that
\begin{equation*}
    D^s(F^t_s(y_{\varphi}, s, 0), F^t_s(\sigma^{\tau}(y_{\varphi^-}), s, \tau\ \mbox{mod}\ T^+(s))) = \eps_s\ \ \mbox{for all}\ s \in \bN.
\end{equation*}

Here,
\begin{equation*}
\begin{split}
    F^t_s(y_{\varphi}, s, 0) &= (\sigma^t(y_{\varphi}), s, t\ \mbox{mod}\ T^+(s)),\\
    F^t_s(\sigma^{\tau}(y_{\varphi^-}), s, \tau\ \mbox{mod}\ T^+(s)) &= (\sigma^{t+\tau}(y_{\varphi^-}), s, (\tau + t)\ \mbox{mod}\ T^+(s)).
\end{split}
\end{equation*}

For $s \in \bN$ such that $\tau$ is not a positive integer multiple of $T^+(s)$, we have

$(t\ \mbox{mod}\ T^+(s)) \neq (\tau  + t)\ \mbox{mod}\ T^+(s)$, and thus
\begin{equation*}
    D^s(F^t_s(y_{\varphi}, s, 0), F^t_s(\sigma^{\tau}(y_{\varphi^-}), s, \tau\ \mbox{mod}\ T^+(s))) = \eps_s
\end{equation*}

by (Dn1e).

For $s \in \bN$ such that $\tau$ is a positive integer multiple of $T^+(s)$, we must have

$s < m$, and $(t\ \mbox{mod}\ T^+(s)) = (\tau + t)\ \mbox{mod}\ T^+(s)$.  Furthermore, there exists

$0 \leq t_1 \leq T(0)-1$ such that $t = T^+(m)-\tau + t_1$.  Then, $(t\ \mbox{mod}\ T^+(s)) = t_1$.

By the choice of $t$ we have $\sigma^t(y_{\varphi})(0) = 1$ and  $\sigma^{t+\tau}(y_{\varphi^-})(0) = 0$.  Thus
\begin{equation*}
    \Delta((\sigma^t(y_{\varphi}), \sigma^{t+\tau}(y_{\varphi^-})) = 0.
\end{equation*}

Moreover, $\Phi(\sigma^{t+\tau}(y_{\varphi^-}), s, t_1) = {}^{T^+(s)}\{0\}$.

Thus Proposition~\ref{prop:cCiprop} applies, and we get:
\begin{equation*}
    D^s((\sigma^t(y_{\varphi}), s, t_1), (\sigma^{t+\tau}(y_{\varphi^-}), s, t_1)) = \eps_s.
\end{equation*}

\medskip

\noindent
Case 5: $x \in F^{\tau}(W^m)$ for some $\tau > 0$ and $m > n$.

\smallskip

 There exist $\varphi^+ \in \mathcal{Y}_m^-$ and  $\varphi, \varphi' \in \Ymn$ such that for all $s \geq 0$ we have
\begin{equation*}
    x_s = (\sigma^{\tau}(y_{\varphi^+}), s, \tau\ \mbox{mod}\ T^+(s)),\ \ \ \ u_s = (y_{\varphi}, s, 0),\ \ \ \ v_s = (y_{\varphi'},s,0).
\end{equation*}
In this case, for the value of $\tau$, we distinguish the following five subcases:

\smallskip

Case 5-1: $\tau \geq T^+(m)$.

Case 5-2:  $\tau < T^+(m)$ is not a positive integer multiple of $T^+(0)$.

Case 5-3:  $T^+(m)-T^+(n) < \tau < T^+(m)$ is a positive integer multiple of $T^+(0)$.

Case 5-4:  $0 < \tau \leq T^+(m)-T^+(n)$ is a positive integer multiple of $T^+(n)$.

Case 5-5:  $0 < \tau \leq T^+(m)-T^+(n)$ is a positive integer multiple of $T^+(0)$,

\qquad \qquad \ \  but not a positive integer multiple of $T^+(n)$.

\medskip

\noindent
Case 5-1: $\tau \geq T^+(m)$.

\smallskip

In this case, $\sigma^{\tau}(y_{\varphi^+})(i) = 0$ for all $s \in \bN$ and $i \geq 0$. Choose the smallest

$0 \leq t \leq T(0)-1$ with $\varphi(t) = 1$.  The existence of such a $t$ follows from (PY2+).

We want to show that $D^s(F^t_s(y_{\varphi},s,0), F^t_s(\sigma^{\tau}(y_{\varphi^+}),s, \tau\ \mbox{mod}\ T^+(s))) = \eps_s$

for all $s \in \bN$.  Here,
\begin{equation*}
\begin{split}
    F^t_s(y_{\varphi},s,0) &= (\sigma^t(y_{\varphi}), s, t),\\ F^t_s(\sigma^{\tau}(y_{\varphi^+}),s, \tau\ \mbox{mod}\ T^+(s)) &= (\sigma^{t+\tau}(y_{\varphi^+}), s, (t+\tau)\ \mbox{mod}\ T^+(s)).
\end{split}
\end{equation*}

If $t \neq (t+\tau)\ \mbox{mod}\ T^+(s)$,

then $D^s(F^t_s(y_{\varphi},s,0), F^t_s(\sigma^{\tau}(y_{\varphi^+}),s, \tau\ \mbox{mod}\ T^+(s))) = \eps_s$ by (Dn1e).

If $t = (t+\tau)\ \mbox{mod}\ T^+(s)$, with $\sigma^t(y_{\varphi})(0) = 1 \neq 0 = \sigma^{t+\tau}(y_{\varphi^+})(0)$,

we have $\Delta(\sigma^t(y_{\varphi}), \sigma^{t+\tau}(y_{\varphi^+})) = 0$ and
$\Phi(\sigma^{t+\tau}(y_{\varphi^+}), s, t) = {}^{T^+(s)}\{0\}$.

Thus Proposition~\ref{prop:cCiprop} applies, and we get

$D^s(F^t_s(y_{\varphi},s,0), F^t_s(\sigma^{\tau}(y_{\varphi^+}),s, \tau\ \mbox{mod}\ T^+(s))) = \eps_s$.

\medskip

\noindent
Case 5-2: $\tau < T^+(m)$ is not a positive integer multiple of $T^+(0)$.

\smallskip

In this case, $\tau\ \mbox{mod}\ T^+(s) \neq 0$ for all $s \in \bN$.  Then by (Dn1e) we have

$D^s((\sigma^{\tau}(y_{\varphi^+}), s, \tau\ \mbox{mod}\ T^+(s)), (y_{\varphi}, s, 0)) = \eps_s$ for all $s \in \bN$.  Thus, $D(x, u) = \eps$.

\medskip

\noindent
Case 5-3: $T^+(m)-T^+(n) < \tau < T^+(m)$ is a positive integer multiple of $T^+(0)$.

\smallskip

Choose the smallest $T^+(m)-\tau \leq t \leq T^+(m)-\tau+T(0)-1$ with $\varphi(t) \neq 0$.

The existence of such a $t$ follows from (PY2+).

We will show that
\begin{equation*}
    D^s(F^t_s(y_{\varphi}, s, 0), F^t_s(\sigma^{\tau}(y_{\varphi^+}), s, \tau\ \mbox{mod}\ T^+(s))) = \eps_s\ \ \mbox{for\ all}\ s \in \bN.
\end{equation*}

Here,
\begin{equation*}
\begin{split}
    F^t_s(y_{\varphi}, s, 0) &= (\sigma^t(y_{\varphi}), s, t\ \mbox{mod}\ T^+(s)),\\
    F^t_s(\sigma^{\tau}(y_{\varphi^+}), s, \tau\ \mbox{mod}\ T^+(s)) &= (\sigma^{t+\tau}(y_{\varphi^+}), s, (\tau + t)\ \mbox{mod}\ T^+(s)).
\end{split}
\end{equation*}

For $s \in \bN$ such that $\tau$ is not a positive integer multiple of $T^+(s)$,

we have $(t\ \mbox{mod}\ T^+(s)) \neq (\tau + t)\ \mbox{mod}\ T^+(s)$, and thus by (Dn1e):
\begin{equation*}
    D^s(F^t_s(y_{\varphi}, s, 0), F^t_s(\sigma^{\tau}(y_{\varphi^+}), s, \tau\ \mbox{mod}\ T^+(s))) = \eps_s.
\end{equation*}

For $s \in \bN$ such that $\tau$ is a positive integer multiple of $T^+(s)$,

we must have $s < m$, and $(t\ \mbox{mod}\ T^+(s)) = (\tau + t)\ \mbox{mod}\ T^+(s)$.

Furthermore, there exists $0 \leq t_1 \leq T(0)-1$ such that $t = T^+(m)-\tau + t_1$.

Hence, with both $T^+(m)$ and $\tau$ being positive integer multiples of $T^+(s)$,

we have $(t\ \mbox{mod}\ T^+(s)) = t_1$.

By the choice of $t$ we have $\sigma^t(y_{\varphi})(0) = 1$ and  $\sigma^{t+\tau}(y_{\varphi^+})(0) = 0$.  Thus
\begin{equation*}
    \Delta((\sigma^t(y_{\varphi}), \sigma^{t+\tau}(y_{\varphi^+})) = 0.
\end{equation*}

Moreover, $\Phi(\sigma^{t+\tau}(y_{\varphi^+}), s, t_1) = {}^{T^+(s)}\{0\}$.
Thus Proposition~\ref{prop:cCiprop} applies, and
\begin{equation*}
    D^s((\sigma^t(y_{\varphi}), s, t_1), (\sigma^{t+\tau}(y_{\varphi^+}), s, t_1)) = \eps_s.
\end{equation*}

\medskip

\noindent
Case 5-4:  $0 < \tau \leq T^+(m)-T^+(n)$ is a positive integer multiple of $T^+(n)$.

\smallskip

As $u \neq v$, we have $\varphi \neq \varphi'$.

Wlog, we can assume that $\sigma^{\tau}(y_{\varphi^+})\restrict (0,\dots, T^+(n)-1) \neq \varphi$.

Now it suffices to show that $D_{T^+(n)}(x,u) = \eps$.

By property~(PR1), the assumption that $\varphi \in \mathcal{Y}_m^-$, and our assumptions for

this subcase, $\sigma^{\tau}(y_{\varphi^+})\restrict (0,\dots, T^+(n)-1)$ and  $\varphi$ are both elements of $\Ymn$.

Thus, by Lemma~\ref{prop:separated}(iii), there exists $t < T^+(n)$ such that
\begin{equation*}
D^s(F_s^t(\sigma^{\tau}(y_{\varphi^+}), s, 0), F_s^t(y_{\varphi}, s, 0)) = \eps_s \quad \mbox{for all} \ \ 0 \leq s \leq n.
\end{equation*}

Since $\tau\ \mbox{mod}\ T^+(s) = 0$ for any $0 \leq s \leq n$, for such a $t$ we have
\begin{equation*}
D^s(F_s^t(\sigma^{\tau}(y_{\varphi^+}), s, \tau\ \mbox{mod}\ T^+(s)), F_s^t(y_{\varphi}, s, 0)) = \eps_s \quad  \mbox{for all} \ \ 0 \leq s \leq n.
\end{equation*}

This implies that $\sigma^t(\sigma^{\tau}(y_{\varphi^+}))(0) \neq \sigma^t(y_{\varphi})(0)$ for such a $t$.

Fix $t$ as above.
Given any $s > n$,
\begin{equation*}
\begin{split}
    F_s^t(\sigma^{\tau}(y_{\varphi^+}), s, \tau\ \mbox{mod}\ T^+(s)) &= (\sigma^{t+{\tau}}(y_{\varphi^+}), s, (t+\tau)\ \mbox{mod}\ T^+(s)),\\
    F^t_s(y_{\varphi}, s,0) &= (\sigma^t(y_{\varphi}), s, t).
\end{split}
\end{equation*}

If $(t + \tau)\ \mbox{mod}\ T^+(s) \neq t$, then by (Dn1e):
\begin{equation*}
D^s(F_s^t(\sigma^{\tau}(y_{\varphi^+}), s, \tau\ \mbox{mod}\ T^+(s)), F_s^t(y_{\varphi}, s, 0)) = \eps_s.
\end{equation*}

Assume $(t + \tau)\ \mbox{mod}\ T^+(s) = t$.  By the choice of $t$ we have

$\sigma^t(\sigma^{\tau}(y_{\varphi^+}))(0) \neq \sigma^t(y_{\varphi})(0)$, hence
$\Delta(\sigma^t(\sigma^{\tau}(y_{\varphi^+})), \sigma^t(y_{\varphi})) = ~0$.

Moreover, $\Phi(\sigma^t(y_{\varphi}), s, t)\restrict (0,\dots, T^+(n)-1) = \varphi$
and $\Phi(F^t_s(y_{\varphi}), s, t)(i) = 0$

for all $T^+(n) \leq T^+(s-1) \leq i \leq T^+(s)-1$. Then by Proposition~\ref{prop:cCiprop},
\begin{equation*}
D^s(F_s^t(\sigma^{\tau}(y_{\varphi^+}), s, \tau\ \mbox{mod}\ T^+(s)), F_s^t(y_{\varphi}, s, 0)) = \eps_s.
\end{equation*}

\medskip

\noindent
Case 5-5: $0 < \tau \leq T^+(m)-T^+(n)$ is a positive integer multiple of $T^+(0)$,

\qquad \qquad \ \ but not a positive integer multiple of $T^+(n)$.

\smallskip

As $u \neq v$, we have $\varphi \neq \varphi'$.

Wlog, we can assume that $\varphi \neq \sigma^{\tau}(y_{\varphi^+})\restrict (0, \dots, T^+(n)-1)$.

That is, there exists $0 \leq r \leq T^+(n)-1$ such that $\varphi(r) \neq \sigma^{\tau}(y_{\varphi^+})(r)$.

Let $S$ be the largest number in $\bN$ such that $\tau$ is a positive integer multiple

of $T^+(S)$.  As $\tau$ is not a positive integer multiple of $T^+(n)$, this number $S$

must be less than~$n$. Since $T^+(n)$ is a positive integer multiple of $T^+(S)$,

there exist integers $r_1$ and $r_2$ such that
\begin{itemize}
\item $0 \leq r_1 \leq r \leq r_2 \leq T^+(n)-1$.
\item $r_1$ is an integer multiple of $T^+(S)$.
\item $r_2 - r_1 + 1 = T^+(S)$.
\end{itemize}

Moreover, $\sigma^{\tau}(y_{\varphi^+})\restrict (r_1, \dots, r_2)$ and $\varphi \restrict (r_1, \dots, r_2)$ are
distinct elements of $\mathcal{Y}_S^-$.

Thus Lemma~\ref{prop:separated}(iii) applies with~$r_1$ playing the role of~$\tau$ in condition~(P2-1).

Hence there exists $r_1 \leq t \leq r_2$ such that for all $0 \leq s \leq S$: {\small
\begin{equation*}
\begin{split}
D^s(F_s^t(\sigma^{\tau}(y_{\varphi^+}), s, \tau\ \mbox{mod}\ T^+(s)), F_s^t(y_{\varphi}, s, 0)) &= D^s(F_s^t(\sigma^{\tau}(y_{\varphi^+}), s, 0), F_s^t(y_{\varphi}, s, 0))\\
&= \eps_s.
\end{split}
\end{equation*}
}

This implies that $\sigma^t(\sigma^{\tau}(y_{\varphi^+}))(0) \neq \sigma^t(y_{\varphi})(0)$ for such a $t$.

Fix $t$ as above.

Now it remains to show that
\begin{equation}\label{eqn:Dseps-for-s>S}
    D^s(F_s^t(\sigma^{\tau}(y_{\varphi^+}), s, \tau\ \mbox{mod}\ T^+(s)), F_s^t(y_{\varphi}, s, 0)) = \eps_s \ \ \ \ \mbox{for\ all}\ s > S.
\end{equation}

Note that $\tau$ cannot be a positive integer multiple of $T^+(s)$ for any $s > S$.  Hence,
\begin{equation*}
    t\ \mbox{mod}\ T^+(s) \neq (\tau + t)\ \mbox{mod}\ T^+(s).
\end{equation*}

Thus, \eqref{eqn:Dseps-for-s>S} follows from (Dn1e).

\medskip

\noindent
Case 6: $x \in F^{\tau}(W^m)$ for some $\tau < 0$ and $m \in \bN$.

\smallskip

There exist $\varphi \in \mathcal{Y}_m^-$ and $\varphi' \in \Ymn$ such that for all $s \geq 0$, we have
\begin{equation*}
    x_s = (\sigma^{\tau}(y_{\varphi}), s, \tau\ \mbox{mod}\ T^+(s)),\ \ \ \  u_s = (y_{\varphi'}, s, 0).
\end{equation*}

In this case, we distinguish the following two subcases for the value of $\tau$:

\smallskip

Case 6-1:  $\tau$ is not an integer multiple of $T^+(0)$.

Case 6-2:  $\tau$ is an integer multiple of $T^+(0)$.

\medskip

\noindent
Case 6-1: $\tau$ is not an integer multiple of $T^+(0)$.

\smallskip

In this case, $\tau\ \mbox{mod}\ T^+(s) \neq 0$ for all $s \in \bN$.  Then by (Dn1e), we have

$D^s((\sigma^{\tau}(y_{\varphi}), s, \tau\ \mbox{mod}\ T^+(s)), (y_{\varphi'}, s, 0)) = \eps_s$ for all $s \in \bN$.  Thus, $D(x, u) = \eps$.

\medskip

\noindent
Case 6-2:  $\tau$ is an integer multiple of $T^+(0)$.

\smallskip

Choose the smallest $0 \leq t \leq T(0)-1$ with $\varphi'(t) \neq 0$.

The existence of such a $t$ follows from (PY2+).

We will show that
\begin{equation*}
    D^s(F^t_s(y_{\varphi'}, s, 0), F^t_s(\sigma^{\tau}(y_{\varphi}), s, \tau\ \mbox{mod}\ T^+(s))) = \eps_s\ \ \mbox{for all}\ s \in \bN.
\end{equation*}

Here,
\begin{equation*}
\begin{split}
    F^t_s(y_{\varphi'}, s, 0) &= (\sigma^t(y_{\varphi'}), s, t\ \mbox{mod}\ T^+(s)),\\
    F^t_s(\sigma^{\tau}(y_{\varphi}), s, \tau\ \mbox{mod}\ T^+(s)) &= (\sigma^{t+\tau}(y_{\varphi}), s, (\tau + t)\ \mbox{mod}\ T^+(s)).
\end{split}
\end{equation*}

For $s \in \bN$ such that $\tau$ is not an integer multiple of $T^+(s)$, we have

$(t\ \mbox{mod}\ T^+(s)) \neq (\tau  + t)\ \mbox{mod}\ T^+(s)$, and thus by (Dn1e):
\begin{equation*}
    D^s(F^t_s(y_{\varphi'}, s, 0), F^t_s(\sigma^{\tau}(y_{\varphi}), s, \tau\ \mbox{mod}\ T^+(s))) = \eps_s
\end{equation*}

For $s \in \bN$ such that $\tau$ is an integer multiple of $T^+(s)$, we will have $T^+(s) \leq |\tau|$,

and $(t\ \mbox{mod}\ T^+(s)) = (\tau + t)\ \mbox{mod}\ T^+(s) = t$.

Now we are going to show that in this case
\begin{equation}\label{eqn:Dst-eps-s}
    D^s((\sigma^t(y_{\varphi'}), s, t), (\sigma^{t+\tau}(y_{\varphi}), s, t)) = \eps_s.
\end{equation}

By $|\tau| \geq T^+(0)$ and the choice of $t$ we have $\sigma^t(y_{\varphi'})(0) = 1$
and  $\sigma^{t+\tau}(y_{\varphi})(0) = 0$.

Thus
\begin{equation*}
    \Delta((\sigma^t(y_{\varphi'}), \sigma^{t+\tau}(y_{\varphi})) = 0.
\end{equation*}

Moreover, $\Phi(\sigma^{t+\tau}(y_{\varphi}), s, t) = {}^{T^+(s)}\{0\}$,
and~\eqref{eqn:Dst-eps-s} follows from Proposition~\ref{prop:cCiprop}.

\medskip

\noindent
Case 7:  $x \in X^-\backslash\left[\bigcup_{t\in \bZ}F^t(\bigcup_{n\geq 0}W^n)\right]$.

\smallskip

In this case, there exists $\{x^j\}_{j=1}^{\infty} \subset \bigcup_{t\in\bZ} F^t\left(\bigcup_{n\geq 0}W^n\right)$ such that
\begin{equation*}
    \lim_{j\rightarrow \infty} D_{T^+(n)}(x^j, x) = 0.
\end{equation*}
Assume towards a contradiction that $D_{T^+(n)}(x,u) < \eps$ and that $D_{T^+(n)}(x,v) < \eps$.
Fix $\gamma > 0$ such that
\begin{equation*}
     D_{T^+(n)}(x,u) + \gamma < \eps,\ \ \  D_{T^+(n)}(x,v) + \gamma < \eps.
\end{equation*}
For this $\gamma$, there exists $N \in \bN$ such that for all $j > N$, we have $D_{T^+(n)}(x^j, x) < \frac{\gamma}{2}$.
Thus, for each $j > N$,
\begin{equation*}
    \begin{split}
    D_{T^+(n)}(x^j, u) &\leq D_{T^+(n)}(x^j, x) + D_{T^+(n)}(x, u)\\
    &< \frac{\gamma}{2} + \eps -\gamma \\
    &= \eps - \frac{1}{2}\gamma,
    \end{split}
\end{equation*}
\begin{equation*}
    \begin{split}
    D_{T^+(n)}(x^j, v) &\leq D_{T^+(n)}(x^j, x) + D_{T^+(n)}(x, v)\\
    &< \frac{\gamma}{2} + \eps -\gamma \\
    &= \eps - \frac{1}{2}\gamma.
    \end{split}
\end{equation*}
However, $\{x^j\}_{j=1}^{\infty} \subset \bigcup_{t\in \bZ} F^t\left(\bigcup_{n\geq 0}W^n\right)$.  By Cases 1--6,
\begin{equation*}
    \max\{D_{T^+(n)}(x^j, u),  D_{T^+(n)}(x^j, v)\} = \eps.
\end{equation*}
We arrived at a contradiction.  Therefore,
\begin{equation*}
    \max\{D_{T^+(n)}(x, u),  D_{T^+(n)}(x, v)\} = \eps.
\end{equation*}
$\Box$

\bigskip

\begin{corollary}\label{corol:low-span}
Let $n \in \bN$.  Then $span(X^-, \eps, D_{T^+(n)}) \geq 2^{0.9T^+(n)}$.
\end{corollary}

\noindent
\textbf{Proof:} By Corollary~\ref{corol:size}, $|W^n| \geq 2^{0.9T^+(n)}$.  Thus, it suffices to show that\\
 $span(X^-, \eps, D_{T^+(n)}) \geq |W^n|$.

Assume towards a contradiction that $span(X^-, \eps, D_{T^+(n)}) < |W^n|$.  Then there exists $A \subset X^-$ that is $(T^+(n), \eps)$-spanning with $|A| < |W^n|$.  Hence, by the Pigeonhole Principle, there exist $u \neq v \in W^n$ and $x \in A$ such that
\begin{equation*}
    D_{T^+(n)}(x,u) < \eps,\ \ \ \ D_{T^+(n)}(x,v) < \eps.
\end{equation*}
However, by Lemma~\ref{corol:low-span-deleted},
\begin{equation*}
    \max\{D_{T^+(n)}(x,u), D_{T^+(n)}(x,v)\} = \eps.
\end{equation*}
We arrived at a contradiction.
Therefore, $span(X^-, \eps, D_{T^+(n)}) \geq |W^n|$.  $\Box$

\bigskip

\subsection{Upper bounds on $sep(X, \delta, D_{2T(n)})$}\label{subsec:upper-sep-2Tn}

Here we prove the following result:

\begin{corollary}\label{corol:sep-upper-2Tn}
Let $(X, F)$ be an \JX-system with \JX-metric~$D$ that is based on condition~(Dn32c) with all
colorings~$c_n$ satisfying condition~(cC3).  Let $(Y, F\restrict Y)$ be a subsystem of~$(X,F)$.
Then the following inequalities hold: {\small
\begin{equation}\label{eqn:sepY-upper}
span(Y, \eps, D_{2T(n)}) \leq  sep(Y, \eps, D_{2T(n)}) \leq sep(X, \eps, D_{2T(n)}) \leq T^+(n)2^{1.75T(n)}.
\end{equation} }
\end{corollary}

In the proof of Theorem~\ref{thm:main} we will use~\eqref{eqn:sepY-upper} for $Y = X^-$; in the proof of
Theorem~\ref{thm:main-seponly} we will use~\eqref{eqn:sepY-upper} for $Y = W$.
In both cases the metric~$D$ satisfies the assumptions of the corollary.

The first inequality in~\eqref{eqn:sepY-upper} is part
of~\eqref{ineq:covsepspan} of Lemma~\ref{lem:BasicFacts}.

The second inequality in~\eqref{eqn:sepY-upper} is true because we
assumed~$Y \subset X$.

For the proof of the third inequality in~\eqref{eqn:sepY-upper},
consider a subset $S \subset Y \subseteq X$  that is $(2T(n), \eps)$-separated, and let  $x \neq x' \in S$.
Then there exists $0 \leq t \leq 2T(n)-1$ such that $D(F^t(x),F^t(x')) \geq \eps$. By~\eqref{eqn:large-dist}, this implies that for each~$n \in \bN$
\begin{equation*}
    D^n(F_n^t(x_n),F_n^t(x'_n)) = \eps_n > \delta_n,
\end{equation*}
and if we choose~$\delta = \eps_n$ in the following lemma,
then we obtain the inequality $|S| \leq T^+(n)2^{1.75T(n)}$.

\begin{lemma}\label{corol:low-separation-n}
Let $n \in \bN$, and let $(X_n, F_n)$ be an \JXn-system with \JXn-metric~$D^n$ that is based on condition~(Dn32c) with a
coloring~$c_n$ satisfying condition~(cC3).  Then:
\begin{equation}\label{eqn:liminf-n}
\forall \delta > \delta_n \quad sep(X_n, \delta, D^n_{2T(n)}) \leq  T^+(n)2^{1.75T(n)}.
\end{equation}
\end{lemma}

\noindent
\textbf{Proof:}  Fix $n \in \bN$ and $\delta > \delta_n$.  If there were a subset $B \subset X_n = \cup_{k = 0}^{T^+(n)-1}X_n^k$
that is $(2T(n),\delta)$-separated of size $|B| > T^+(n)2^{1.75T(n)}$, by the Pigeonhole Principle, there
would exist $A \subset B$ with $|A| > 2^{1.75T(n)}$ and $0 \leq k \leq T^+(n)-1$ such that $A \subset X_n^k$.
Note that $A$ would still be $(2T(n),\delta)$-separated so that the inequality\\
$sep(X_n^k, \delta, D^n_{2T(n)}) > 2^{1.75T(n)}$ would hold.

Hence, it suffices to show that for all $0 \leq k \leq T^+(n)-1$
\begin{equation*}
    sep(X_n^k, \delta, D^n_{2T(n)}) \leq 2^{1.75T(n)}.
\end{equation*}

\medskip

Consider any $0 \leq k \leq T^+(n)-1$ and subset $A \subset X_n^k$ that is $(2T(n),\delta)$-separated.

For all $(y,n,k) \neq (y',n,k) \in A$, there exists $0 \leq t \leq 2T(n)-1$ such that
$D^n(F_n^t((y,n,k)),F_n^t((y',n,k))) \geq \delta > \delta_n$. By (PDn1) and (PDn2),
\begin{equation}\label{eqn:epsn-somewhere}
    \begin{split}
        \eps_n &= D^n(F_n^t((y,n,k)),F_n^t((y',n,k))) \\
        &= D^n((\sigma^t(y),n,(t+k)\ mod \ T^+(n)),(\sigma^t(y'),n,(t+k) \ mod \ T^+(n))).
    \end{split}
\end{equation}

Let us introduce some new notation.  Consider $J \subset \{0, 1, \dots , 2T(n)-1\}$, and let $y^- \in {}^J \{0,1\}$. Define
\begin{equation*}
A(y^-) = \{(y,n, k) \in A: \ y\restrict J = y^-\}.
\end{equation*}

Note that we can suppress the parameter~$J$ in this notation since it must be the domain of~$y^-$.
Let us make a few observations:

\smallskip

\begin{itemize}
\item[(Union)] $A = \bigcup_{y^- \in {}^J \{0,1\}} A(y^-)$.
\end{itemize}

\smallskip

Note that (Union)  implies that
\begin{equation*}
|A| \leq \sum_{y^- \in {}^J \{0,1\}} |A(y^-)|,
\end{equation*}
which by the Pigeonhole Principle implies, in particular,

\smallskip

\begin{itemize}
\item[(Size)] If $|J| = T(n)$, then there exists $y^- \in {}^J \{0,1\}$ such that $|A(y^-)| \geq |A|2^{-T(n)}$.
\end{itemize}

\smallskip

We will only use the following consequence
of~(Size):

\smallskip

\begin{itemize}
\item[(UseJ)] Assume that $|A| > 2^{1.75T(n)} = 2^{T(n)}2^{0.75T(n)}$ and $|J| = T(n)$.\\
Then there exists $y^- \in {}^J \{0,1\}$
such that $|A(y^-)| > 2^{0.75T(n)}$.
\end{itemize}

\smallskip

Let us make one more observation that follows from the definitions of~$A(y^-)$, $F_n$, and~$D_n$.  Namely, if $t \in J$, and
$(y,n,k) \neq (y',n,k) \in A(y^-)$, then $y(t) = y'(t)$, which is equivalent to $\sigma^t(y)(0) = \sigma^t(y')(0)$.  Hence \newline
$\Delta(\sigma^t(y), \sigma^t(y')) > 0$, so that (Dn31) of the definition of~$D^n$ applies.
Note that \eqref{eqn:epsn-somewhere} implies that
for all $(y,n,k) \neq (y',n,k) \in A$, there exists $0 \leq t \leq 2T(n)-1$ such that $y(t) \neq y'(t)$.  Thus:

\smallskip

\begin{itemize}
\item[(tnotJ)] For any $y^- \in   {}^J \{0,1\}$ and $(y,n,k) \neq (y',n,k) \in A(y^-)$, there exists \newline
$0 \leq t \leq 2T(n)-1$ with $t \notin J$ such that $y(t) \neq y'(t)$ and \newline
$D^n(F_n^t((y,n,k)),F_n^t((y',n,k))) = \eps_n$.
\end{itemize}

\smallskip

Now assume towards a contradiction that~$|A| > 2^{1.75T(n)}$.

We distinguish two cases.

\medskip

\noindent
Case 1: $T(n) - 1 + k < T^+(n)$

\smallskip

In this case we let $J = \{T(n), T(n+1), \dots,  2T(n)-1\}$. By (UseJ) we can pick $y^- \in   {}^J \{0,1\}$ such that $|A(y^-)| > 2^{0.75 T(n)}$.
Let $A^- := A(y^-)$.

Consider $(y,n,k) \neq (y',n,k) \in A^-$. Our choice of~$J$ implies the following consequence of (tnotJ):

\smallskip

\begin{itemize}
\item[(tsm)]  There exists $0 \leq t \leq T(n)-1$ such that $y(t) \neq y'(t)$ and\\ $D^n(F_n^t((y,n,k)),F_n^t((y',n,k))) = \eps_n$.
\end{itemize}

\smallskip

For any such~$t$, the distance $D^n(F_n^t((y,n,k)),F_n^t((y',n,k)))$ is calculated according to
clause~(Dn32c).  Let $\varphi = (y(-k), y(-k+1), \dots , y(-k + T^+(n) - 1))$ and $\psi = (y'(-k), y'(-k+1),
\dots , y'(- k + T^+(n) - 1))$.  Then the defining property of Case~1 together with Proposition~\ref{prop:Phi-const} imply that we have \newline
$(0, 1, \dots, T(n) - 1) \subset
(-k, -k+1, \dots, -k + T^+(n)- 1)$, so that  \newline
$\Phi(F_n^t((y,n,k))) = \varphi$ and $\Phi(F_n^t((y',n,k))) = \psi$ for
all~$0 \leq t < T(n)$, \newline
which in turn implies together with (tsm):

\smallskip

\begin{itemize}
\item[(Diff)] $\Phi((y,n,k)) \neq
\Phi((y',n,k))$ whenever $(y,n,k) \neq (y',n,k) \in A^-$.
\end{itemize}

\smallskip

Let $j(k)$ be such that $k \in I^n_{j(k)}$, and let $j = c_n(\varphi,\psi) = c_n(\Phi((y,n,k)), \Phi((y',n,k)))$.

\smallskip

If for some $t < T(n)$
\begin{equation*}
    \begin{split}
    &D^n(F_n^t((y,n,k)),F_n^t((y',n,k))) \\ &= D^n((\sigma^t(y),n,(t+k)\ \mbox{mod}\ T^+(n)),(\sigma^t(y'),n,(t+k)\ \mbox{mod}\ T^+(n))) \\ &= \eps_n,
    \end{split}
\end{equation*}
then according to~(Dn32c) we must have
\begin{equation}\label{eqn:t+kInj}
(t+k) \ mod \ T^+(n) = t+k \in I^n_j.
\end{equation}

Since $I^n_j$ has length~$T(n)$ and $t < T(n)$, \eqref{eqn:t+kInj} can hold only if~$j = j(k)$ or if $j = j(k) + 1$.  As this observation does not
depend on the particular choice of~$t$ and of $(y,n,k), (y',n,k) \in A^-$, we conclude that
\begin{equation}\label{eqn:2-homo1}
\forall (y,n,k) \neq (y',n,k) \in A^- \ c_n(\Phi((y,n,k)), \Phi((y',n,k))) \in \{j(k), j(k)+1\}.
\end{equation}

Now define $S = \{\Phi((y,n,k)): (y,n,k) \in A^-\}$. It follows from~\eqref{eqn:2-homo1} that the restriction of $c_n$ to $[S]^2$ takes on at most two colors: ~$j(k)$
or~$j(k) + 1$.  Moreover, (Diff) implies that
$|S| = |A^-|$.  Thus by (cC3), $|A^-| = |S| < 2^{0.75T(n)}$. This contradicts  our assumption about the sizes of~$A$ and~$A^-$.

\newpage

\noindent
Case 2: $T(n) - 1 + k \geq T^+(n)$

\smallskip

Let $t_0 < T(n)$ be such that $t_0 + k = T^+(n)$ and choose
\begin{equation*}
J = \{0, \dots, t_0-1\} \cup \{t_0 + T(n), \dots, 2T(n)-1\}.
\end{equation*}

By (UseJ) we can pick $y^- \in   {}^J \{0,1\}$ such that $|A(y^-)| > 2^{0.75 T(n)}$.\newline
Let $A^- := A(y^-)$.

Consider $(y,n,k) \neq (y',n,k) \in A^-$.  Our choice of~$J$ implies the following consequence of (tnotJ):

\smallskip

\begin{itemize}
\item[(tlg)] There exists $t_0 \leq t \leq t_0 + T(n)-1 < 2T(n)$ such that $y(t) \neq y'(t)$ and\\ $D^n(F_n^t((y,n,k)),F_n^t((y',n,k))) = \eps_n$.
\end{itemize}

\smallskip

For any such~$t$, the distance $D^n(F_n^t((y,n,k)),F_n^t((y',n,k)))$ is calculated according to clause~(Dn32c).

Let $u = \sigma^{t_0}(y)$, $u' = \sigma^{t_0}(y')$ and
\begin{equation*}
    \begin{split}
    \varphi &= (u(0), u(1), \dots , u(T^+(n) -1)) = (y(t_0), \dots, y(t_0 + T^+(n)-1)),\\
    \psi &= (u'(0), u'(1), \dots , u'(T^+(n) -1)) = (y'(t_0), \dots, y'(t_0+T^+(n)-1)).
    \end{split}
\end{equation*}
Then the defining property of Case~2 together with Proposition~\ref{prop:Phi-const} imply that we have
\begin{equation*}
    \begin{split}
    \Phi(F_n^t((y,n,k))) &= \Phi(F^{t-t_0}_n((u,n,0))) = \varphi\ \mbox{and}\\
    \Phi(F_n^t((y',n,k))) &= \Phi(F^{t-t_0}_n((u',n,0))) = \psi
    \end{split}
\end{equation*}
for all~$t_0 \leq t < t_0 + T(n) < 2T(n)$.

Since
$(t_0, t_0+1, \dots, t_0+T(n) - 1) \subset
(t_0, t_0+1, \dots, t_0+T^+(n)-1)$, \newline by (tnotJ) we have

\smallskip

\begin{itemize}
\item[(Diff')] $\Phi(F_n^{t_0}((y,n,k))) \neq
\Phi(F_n^{t_0}((y',n, k)))$ whenever $(y,n,k) \neq (y',n,k) \in A^-$.
\end{itemize}

\smallskip

Let $j = c_n(\varphi,\psi) = c_n(\Phi(F^{t_0}((y,n,k))), \Phi(F^{t_0}((y',n,k))))$.

\smallskip

If for some $t_0 \leq t < t_0 + T(n)$
\begin{equation*}
    \begin{split}
    &D^n(F_n^t((y,n,k)),F_n^t((y',n,k))) \\
    &= D^n((\sigma^t(y),n,(t+k)\ \mbox{mod}\ T^+(n)),(\sigma^t(y'),n,(t+k)\ \mbox{mod}\ T^+(n)))\\
    &= \eps_n,
    \end{split}
\end{equation*}
then according to~(Dn32c) we must have
\begin{equation}\label{eqn:t0+kInj}
(t+k) \ mod \ T^+(n) = t - t_0 \in I^n_j.
\end{equation}

Note that $0 \in I^n_1$, the interval~$I^n_j$ has length~$T(n)$, and $t - t_0 < T(n)$.
Thus~\eqref{eqn:t0+kInj} can hold only if~ $j =  1$.  As this observation does not
depend on the particular choice of~$t$ and of $(y,n,k), (y',n,k) \in A^-$, we conclude that
\begin{equation}\label{eqn:2-homo}
\forall (y,n,k) \neq (y',n,k) \in A^- \ c_n(\Phi(F^{t_0}((y,n,k))), \Phi(F^{t_0}((y',n,k)))) = j = 1.
\end{equation}

Now define $S = \{\Phi(F^{t_0}((y,n,k))): (y,n,k) \in A^-\}$. It follows from~\eqref{eqn:2-homo} that the restriction of $c_n$ to $[S]^2$ takes on only one, and thus at most two, values.
 By (Diff'),
$|S| = |A^-|$.  Thus by (cC3), $|A^-| = |S| < 2^{0.75T(n)}$. Again, this contradicts  our assumption about the sizes of~$A$ and~$A^-$. $\Box$

\bigskip

\section{Proof of Theorem~\ref{thm:main}(ii) and  Theorem~\ref{thm:main-seponly}(ii)}\label{sec:prove-main}

Let the expression $N(Z, \eps, D_{T^+(n)})$ stand either for $sep(W,\eps,D_{T^+(n)})$,\newline $sep(X^-,\eps,D_{T^+(n)})$, or $span(X^-,\eps,D_{T^+(n)})$. Then
By Corollary~\ref{corol:+sep} and Corollary~\ref{corol:low-span}, for all $n \in \bN$,
\begin{equation*}
    \begin{split}
        N(Z,\eps,D_{T^+(n)}) &\geq 2^{0.9T^+(n)},\\
        \ln{N(Z,\eps,D_{T^+(n)})} &\geq \ln{\left(2^{0.9T^+(n)}\right)} = 0.9T^+(n)\ln{2},\\
        \frac{\ln{N(Z,\eps,D_{T^+(n)})}}{T^+(n)} &\geq 0.9\ln{2}.
    \end{split}
\end{equation*}

Thus,
\begin{equation}\label{eqn:not-eq+}
\begin{split}
\limsup_{T \rightarrow \infty} \frac{\ln sep (W, \eps, D_T)}{T} &\geq 0.9\ln 2,\\
\limsup_{T \rightarrow \infty} \frac{\ln sep (X^-, \eps, D_T)}{T} &\geq  \limsup_{T \rightarrow \infty} \frac{\ln span (X^-, \eps, D_T)}{T} \geq 0.9\ln 2.
\end{split}
\end{equation}

\medskip

Next we fix $n$. It follows from Corollary~\ref{corol:sep-upper-2Tn} that

\begin{equation*}
    \begin{split}
    \ln{sep(X,\eps, D_{2T(n)})} &\leq \ln{\left(T^+(n)2^{1.75T(n)}\right)} = \ln{T^+(n)} + 1.75T(n)\ln{2}\\
    &= \ln{C(n)} + \ln{T(n)} + 1.75T(n)\ln{2},\\
        \frac{\ln{sep(X,\eps, D_{2T(n)})}}{2T(n)} &\leq  \frac{\ln{C(n)}}{2T(n)} + \frac{\ln{T(n)}}{2T(n)} + 0.875\ln{2}.
    \end{split}
\end{equation*}

By (PKn5), $C(n) \leq 2^{0.01T(n)}$.  Then
\begin{equation*}
    \begin{split}
    \frac{\ln{sep(X,\eps, D_{2T(n)})}}{2T(n)} &\leq \frac{\ln{2^{0.01T(n)}}}{2T(n)} + \frac{\ln{T(n)}}{2T(n)} + 0.875\ln{2}\\
    &= \frac{\ln{T(n)}}{2T(n)} + 0.88\ln{2}.
    \end{split}
\end{equation*}

Thus for $X^- \subset X$ we get

\begin{equation}\label{eqn:not-eq-2}
\liminf_{T \rightarrow \infty} \frac{\ln span (X^-, \eps, D_T)}{T} \leq \liminf_{T \rightarrow \infty} \frac{\ln sep (X^-, \eps, D_T)}{T} \leq 0.88 \ln 2.
\end{equation}

Finally, we get from~\eqref{eqn:not-eq+} and~\eqref{eqn:not-eq-2}: {\small
\begin{equation}\label{eqn:ineqs-final}
\begin{split}
\liminf_{T \rightarrow \infty} \frac{\ln sep (X, \eps, D_T)}{T} \leq 0.88 \ln 2 &< 0.9\ln{2} \leq \limsup_{T \rightarrow \infty} \frac{\ln sep (W, \eps, D_T)}{T}, \quad \\
\liminf_{T \rightarrow \infty} \frac{\ln sep (X^-, \eps, D_T)}{T} \leq 0.88 \ln 2 &< 0.9\ln{2} \leq \limsup_{T \rightarrow \infty} \frac{\ln sep (X^-, \eps, D_T)}{T}, \quad \\
\liminf_{T \rightarrow \infty} \frac{\ln span (X^-, \eps, D_T)}{T} \leq 0.88 \ln 2 &< 0.9\ln{2} \leq \limsup_{T \rightarrow \infty} \frac{\ln span (X^-, \eps, D_T)}{T}. \quad
\end{split}
\end{equation}
}

Since each~$(T,\eps)$-separated subset of~$W$ is also~$(T,\eps)$-separated in~$X$, the first line of~\eqref{eqn:ineqs-final} implies Theorem~\ref{thm:main-seponly}(ii). The other lines of~\eqref{eqn:ineqs-final} imply Theorem~\ref{thm:main}(ii).
$\Box$

\medskip

\section{Proof of Theorem~\ref{thm:main-spanonly}}\label{sec:prove-spanonly}

\subsection{A general observation}\label{subsec:spanonly-general} The construction relies on the following observation:

\begin{lemma}\label{lem:duplicate}
Let $(X^-, D), (Y, d)$ be two compact metric spaces with $X^- \cap Y = \emptyset$, and let $F: X^- \rightarrow X^-$,
$G: Y \rightarrow Y$ be homeomorphisms. Let $H = F \cup G$ (where functions are treated as sets of ordered pairs). Assume, moreover, that  $f: X^- \rightarrow Y$ is a conjugacy of the systems $(X^-, F)$ and
$(Y,G)$ such that
\begin{equation}\label{eqn:D<d}
\forall x , x' \in X^- \ D(x, x') \leq d(f(x), f(x')).
\end{equation}

Let~$diam(Y,d) > \alpha > 0.5\max \{diam(X^-, D), diam(Y,d)\} = 0.5 diam(Y,d)$ and consider the function
 $\rho$ on $X^- \cup Y$ that is defined by the following conditions:

 \smallskip

\begin{itemize}
\item[(RD)] $\forall x, x' \in X^- \ \rho(x,x') = D(x,x')$.

\smallskip

\item[(Rd)] $\forall y, y' \in Y \ \rho(y,y') = d(y,y')$.

\smallskip

\item[(R2)] $\forall x \in X^- \, \forall y \in Y \ \rho(x, y) = \rho(y, x) = \max \{\alpha, D(x,f^{-1}(y))\}$.
\end{itemize}

\smallskip

Then

\smallskip

\begin{itemize}
\item[(i)] $(X^- \cup Y, \rho)$ is compact, \\
$diam(X^- \cup Y, \rho) = diam(Y, d)$, and \\
$H$ is a homeomorphism with respect to $\rho$.

\item[(ii)] For all  $T > 0$:
\begin{equation}\label{eqn:sepbd}
\forall \delta > 0 \quad sep(Y, \delta, d_T) \leq sep(X^- \cup Y, \delta, \rho_T) \leq 2sep(Y, \delta, d_T),
\end{equation}

\begin{equation}\label{eqn:spanbd}
\forall \delta > \alpha \quad span(X^-, \delta, D_T) = span(X^- \cup Y, \delta, \rho_T).
\end{equation}
\end{itemize}
\end{lemma}

Note that in~\eqref{eqn:sepbd} and~\eqref{eqn:spanbd} the metrics $D_T, d_T, \rho_T$ are computed for $F, G, H$, respectively.

\medskip

We can derive Theorem~\ref{thm:main-spanonly} if we choose the ingredients of Lemma~\ref{lem:duplicate} as follows:

\begin{itemize}
\item $X^-, D, F, \eps$ will be as constructed in the proof of Theorem~\ref{thm:main}.
\item $Y, d, G, f$ satisfy the assumptions of Lemma~\ref{lem:duplicate} for these choices of~$X^-, D, F$,
and~$\eps$, with $diam(Y, d) = \eps$.
\item  $\alpha$ will be any positive real that satisfies the inequalities $0.5\eps < \alpha < \eps$.
\item Moreover, we will choose
$Y, d, G$ in such a way that
\begin{equation}\label{eqn:sep=inY}
\forall \delta > 0 \ \liminf_{T \rightarrow \infty} \frac{\ln sep (Y, \delta, d_T)}{T} =\limsup_{T \rightarrow \infty} \frac{\ln sep (Y, \delta, d_T)}{T}.
\end{equation}
\end{itemize}

The space $Z$ in Theorem~\ref{thm:main-spanonly} will then be $X^- \cup Y$, and $\rho, H$ will be the objects guaranteed by Lemma~\ref{lem:duplicate}. 
Theorem~\ref{thm:main-spanonly} will follow for these choices.

More specifically, point~(i) of Theorem~\ref{thm:main-spanonly} will follow from point~(i) of Lemma~\ref{lem:duplicate}.

The inequality~\eqref{eqn:not-eq-spanonly} in point~(ii) of Theorem~\ref{thm:main-spanonly} will follow from the analogous inequality in
Theorem~\ref{thm:main} and ~\eqref{eqn:spanbd}.

Equality~\eqref{eqn:eqsep} in point~(iii) of
Theorem~\ref{thm:main-spanonly} will follow  from~\eqref{eqn:sepbd} and~\eqref{eqn:sep=inY}.

Finally, point~(iv) of Theorem~\ref{thm:main-spanonly} will follow by our construction from point~(iii) of Theorem~\ref{thm:main} and~\eqref{eqn:sepbd}.

\medskip

\noindent
\textbf{Proof of Lemma~\ref{lem:duplicate}:} (i) First we prove that this $\rho$ is a metric on $X^-\cup Y$.
\begin{itemize}
\item Reflexivity:  For all $z \in X^-\cup Y$, $\rho(z,z) = D(z,z) = 0$ if $z \in X^-$ and $\rho(z,z) = d(z,z) = 0$ if $z \in Y$.
\item Positive definiteness:  For all $z \neq z' \in Z$, we have $\rho(z,z') = D(z,z') > 0$ if $z, z' \in X^-$, $\rho(z,z') = d(z,z') > 0$ if $z, z' \in Y$ and $\rho(z, z') \geq \alpha > 0$ if $z \in X^-$, $z' \in Y$ or $z \in Y$, $z' \in X^-$.
\item Symmetry:  For all $z, z' \in X^-\cup Y$,
\begin{itemize}
\item if $z, z' \in X^-$, then $\rho(z, z') = D(z,z') = D(z', z) = \rho(z', z)$;
\item if $z, z' \in Y$, then $\rho(z, z') = d(z,z') = d(z', z) = \rho(z', z)$;
\item if $z \in X^-$ and $z' \in Y$, then $\rho(z, z') = \rho(z', z)$ by (R2).
\end{itemize}
\item Triangle Inequality:  For all $z, z', z'' \in X^-\cup Y$ that are pairwise distinct, we show that $\rho(z,z') + \rho(z,z'') \geq \rho(z', z'')$.  We distinguish the following cases:
\begin{itemize}
\item[Case 1:]  $z, z', z'' \in X^-$ or $z, z', z'' \in Y$.  In this case, $\rho(z,z') + \rho(z,z'') \geq \rho(z', z'')$ follows directly from the assumption that $D$ and $d$ are metrics on $X^-$ and $Y$, respectively.
\item[Case 2:]  $z, z' \in X^-$ and $z'' \in Y$.
\begin{itemize}
\item[Case 2-1:]  $\alpha \geq D(z, f^{-1}(z'')), D(z',f^{-1}(z''))$. \\ Then $\rho(z, z') + \rho(z, z'') = D(z, z') + \alpha \geq \alpha = \rho(z', z'')$.
\item[Case 2-2:]  $\alpha \leq D(z, f^{-1}(z'')), D(z',f^{-1}(z''))$.  \\ Then $\rho(z, z') + \rho(z, z'') = D(z, z') + D(z, f^{-1}(z'')) \geq D(z', f^{-1}(z'')) = \rho(z', z'')$.
\item[Case 2-3:]  $D(z',f^{-1}(z'')) \leq \alpha \leq D(z, f^{-1}(z''))$.  \\ Then $\rho(z, z') + \rho(z, z'') = D(z, z') + D(z, f^{-1}(z'')) \geq D(z,z') + \alpha  \geq \alpha = \rho(z', z'')$.
\item[Case 2-4:]  $D(z,f^{-1}(z'')) \leq \alpha \leq D(z', f^{-1}(z''))$.  \\ Then $\rho(z, z') + \rho(z, z'') = D(z, z') + \alpha \geq D(z,z') + D(z, f^{-1}(z'')) \geq D(z', f^{-1}(z'')) = \rho(z', z'')$.
\end{itemize}
\item[Case 3:]  $z', z'' \in X^-$ and $z \in Y$.  In this case, $\rho(z, z') + \rho(z, z'') \geq D(z', f^{-1}(z)) + D(z'', f^{-1}(z)) \geq D(z', z'') = \rho(z', z'')$.
\item[Case 4:]  $z, z' \in Y$ and $z'' \in X^-$.
\begin{itemize}
\item[Case 4-1:]  $\alpha \geq D(z'', f^{-1}(z)), D(z'',f^{-1}(z'))$. \\ Then $\rho(z, z') + \rho(z, z'') = d(z, z') + \alpha \geq \alpha = \rho(z', z'')$.
\item[Case 4-2:]  $\alpha \leq D(z'', f^{-1}(z)), D(z'',f^{-1}(z'))$.  \\ Then $\rho(z, z') + \rho(z, z'') = d(z, z') + D(z'', f^{-1}(z)) \geq D(f^{-1}(z), f^{-1}(z')) + D(z'', f^{-1}(z)) \geq D(z'', f^{-1}(z')) = \rho(z', z'')$.
\item[Case 4-3:]  $D(z'',f^{-1}(z)) \leq \alpha \leq D(z'', f^{-1}(z'))$.  \\ Then $\rho(z, z') + \rho(z, z'') = d(z,z') + \alpha \geq D(f^{-1}(z), f^{-1}(z')) + D(z'', f^{-1}(z)) \geq D(z'', f^{-1}(z'))  = \rho(z', z'')$.\\
Here we used~\eqref{eqn:D<d} and the Triangle Inequality for~$D$.
\item[Case 4-4:]  $D(z'',f^{-1}(z')) \leq \alpha \leq D(z'', f^{-1}(z))$.  \\ Then $\rho(z, z') + \rho(z, z'') = d(z, z') + D(z'', f^{-1}(z)) \geq d(z,z') + \alpha \geq \alpha = D(z', z'')$.
\end{itemize}
\item[Case 5:]  $z', z'' \in Y$ and $z \in X^-$.  In this case,  $\rho(z, z') + \rho(z, z'') \geq 2\alpha > \max\{diam(X^-,D), diam(Y,d)\} \geq d(z', z'') = \rho(z, z'')$.
\end{itemize}
\end{itemize}
Hence, $\rho$ is a metric on $X^- \cup Y$.

\smallskip

The compactness of $(X^-\cup Y, \rho)$ follows directly from~(RD) and (Rd) and the assumptions that $(X^-, D)$ and $(Y,d)$ are compact metric spaces.

Moreover, for any $z, z' \in X^-\cup Y$, we have \\
$\rho(z,z') \leq \max\{diam(X^-,D), diam(Y,d), \alpha\} \leq diam(Y, d)$,\\
so that $diam(X^- \cup Y, \rho) = diam(Y,d)$.

\smallskip

Similarly, the continuity of $H = F\cup G$ follows from the assumptions that $F$ and~$G$ are homeomorphisms and that for all $x \in X^-$ and $y \in Y$, we have $\rho(x,y) \geq \alpha$.

\medskip

\noindent
(ii)  Let~$f$ be as in the assumption and fix $T > 0$.

Then for every~$t\in \bN$ and $x, x' \in X^-$ we have
\begin{equation*}
D(F^t(x), F^t(x')) \leq d(f(F^t(x)), f(F^t(x'))) = d(G^t(f(x)), G^t(f(x'))),
\end{equation*}
and it follows that~\eqref{eqn:D<d} can be extended to
\begin{equation}\label{eqn:DT<dT}
\forall x , x' \in X^- \quad D_T(x, x') \leq d_T(f(x), f(x')).
\end{equation}

Similarly, from (RD) and (Rd) we get
\begin{equation}\label{eqn:RDdT}
\begin{split}
&\forall x , x' \in X^- \quad \rho_T(x, x') = D_T(x, x'),\\
&\forall y , y' \in Y \quad \rho_T(y, y') = d_T(y, y').
\end{split}
\end{equation}

The second line of~\eqref{eqn:RDdT} implies that every $(T, \delta)$-separated subset of~$(Y, d)$ remains a $(T, \delta)$-separated subset of~$(X^- \cup Y, \rho)$, which in turn implies the first inequality in~\eqref{eqn:sepbd}. Similarly, if $A$ is a $(T, \delta)$-separated subset of~$(X^- \cup Y, \rho)$, then
$A \cap X^-$ must be a $(T, \delta)$-separated subset of~$(X^-, D)$ with $|A\cap X^-| \leq sep(X^-,\delta, D_T)$ and $A \cap Y$ must be a $(T, \delta)$-separated subset of~$(Y, d)$, with $|A\cap Y| \leq sep(Y, \delta, d_T)$. Moreover, \eqref{eqn:DT<dT} then implies that $f(A\cap X^-) = \{f(x): \ x \in A \cap X^-\}$ must also be a $(T, \delta)$-separated subset of~$(Y, d)$, with $|f(A\cap X^-)| \leq sep(Y, \delta, d_T)$, and the second inequality in~\eqref{eqn:sepbd} follows.

Now let $\delta > \alpha$, and let $A \subset X^-$ be a $(T, \delta)$-spanning subset of~$(X^-, D)$. Then $A$ remains a spanning set of~$(X^- \cup Y, \rho)$ by~(R2), since
every value $f^{-1}(z)$ must have a distance in~$(X^-,D)$ of less than~$\delta$ from some~$x \in A$.
It follows that $span(X^-, \delta, D_T) \geq span(X^- \cup Y, \delta, \rho_T)$. Conversely,
if $A \subset X^- \cup Y$ is a $(T, \delta)$-spanning subset of~$(X^- \cup Y, \rho)$, then the set $B := \{(A \cap X^-) \cup f^{-1}(A \cap Y)\}$ has cardinality $|B| \leq |A|$ and is a $(T, \delta)$-spanning subset of~$(X^-, D)$ by~(R2).  This implies the inequality $span(X^-, \delta, D_T) \leq span(X^- \cup Y, \delta, \rho_T)$, and concludes the proof of~\eqref{eqn:spanbd}. $\Box$

\subsection{Choosing $Y^+, d, G^+$}\label{subsec:chooseY}

We will essentially construct $(Y, G)$ as an \JX-system with \JX-metric~$d$, except for using slightly different ingredients.
We rely on the same (relevant) parameter choices as in earlier sections and we use the same notation, with one exception: Here for all $n \in \bN$ the elements of the
coordinate spaces~$X_n$ of the \JX-space~$X$ will be denoted by $(v, n, k)$, where $v \in {}^{\bZ}\{0,1\}$, $k \in \{0, 1, \dots,  T^+(n)-1\}$.

Let us first define $Y_n,  G_n, d^n$
for $n\geq 0$.

\begin{itemize}
\item $Y_n$:  The set $Y_n$ consists of all triples $(u, n, k)$, where $u \in {}^\bZ \{2, 3\}$ and $k \in \{0, 1, \dots , T^+(n)-1\}$.

Let $Y_n^k = \{(u, n, k') \in Y_n: \ k' = k\}$.

Then the sets $Y_n^k$ are pairwise disjoint and $Y_n = \bigcup_{0\leq k < T^+(n)} Y^k_n$.

\smallskip

\item $G_n$:  We define the function~$G_n: Y_n \rightarrow Y_n$ as follows:
\begin{equation}\label{eqn:def-Gn}
\begin{split}
G_n((u, n, k)) &= (\sigma(u), n, G_n(k)), \ \mbox{where}\\
\sigma(u)(i) &= u(i+1) \ \mbox{ for all } i,\\
G_n(k) &= (k+1) \ mod \ T^+(n).
\end{split}
\end{equation}

\smallskip

\item $d^n$:
Now we define $d^n((u,n,k), (u',n,k'))$ as follows:
\begin{itemize}
\item[(dn1)] If $k \neq k'$, then $d^n((u, n, k), (u', n, k')) = \eps_n$.
\item[(dn2)] If $k = k'$ and $u = u'$, then $d^n((u, n, k), (u', n, k')) = 0$.
\item[(dn3)] If $k = k'$ and $u \neq u'$, then $d^n ((u, n, k), (u', n, k')) = \eps_n 3^{-\Delta(u, u')}$.
\end{itemize}

\end{itemize}

Notice that each $(Y_n, G_n)$ is a \JXn-system with \JXn-metric~$d^n$ except for some renaming of the ingredients.
Here we don't partition~(dn3) into
subclauses though.  When $\Delta(u, u') = 0$, then we always get $d^n ((u, n, k), (u', n, k')) = \eps_n$, which is the maximum value allowed
by clause~(Dn32) in the definition of an \JXn-metric~$D^n$.  Thus we automatically get the following instance of Proposition~\ref{prop:metric-n}:

\begin{proposition}\label{prop:metric-n-4}
Let $d^n$ be defined as above.  Then

\smallskip

\begin{itemize}

\item[(i)] The function $d^n$ is a metric on~$Y_n$.

\smallskip

\item[(ii)] The systems~$(Y_n, d^n, G_n)$ have the following properties:

\smallskip

            \begin{itemize}
            \item[{}]\emph{(PDn1)} $\max\{ d^n(y, y'): \, y, y' \in Y_n\} = \eps_n$.

            \smallskip

            \item[{}]\emph{(PDn2)} $d^n(y, y') < \eps_n \ \Rightarrow \ d^n(y, y') \leq \delta_n$.

            \smallskip

            \item[{}]\emph{(Pnc)} Each $d^n$ is a metric on~$Y_n$ that induces the topology of a compact Hausdorff space.

            \smallskip

            \item[{}]\emph{(PFn)} Each $G_n : Y_n \rightarrow Y_n$ is a homeomorphism.
            \end{itemize}

\end{itemize}
\end{proposition}

Now we define the following objects:

\begin{itemize}
\item $Y^+$:  Let $Y^+ = \prod_{n \in \bN} Y_n$.  That is, we let $Y^+$ consist of all sequences
\newline $y = (y_n)_{n \in \bN}$ such that $y_n \in Y_n$ for each~$n \in \bN$.
\item $G^+$:  For  $y \in Y^+$, define $G^+(y)_n = G_n(y_n)$ for all~$n \in \bN$.
\item $d$:  The function~$d: \left(Y^+\right)^2 \rightarrow [0,\infty)$ is defined as:
\begin{equation}\label{eqn:define-d}
d(y, y') = \sum_{n \in \bN} d^n(y_n, y'_n).
\end{equation}

\end{itemize}

\smallskip

The system $(Y^+, G^+)$ is an \JX-system with \JX-metric $d$ except for some renaming of the ingredients,
and we get the following instances of Propositions~\ref{prop:D-OK} and~\ref{prop:X-OK}:

\begin{proposition}\label{prop:d-OK}
The function~$d$ as defined in~\eqref{eqn:define-d} is a metric on~$Y^+$ that induces the product topology.
\end{proposition}

\begin{proposition}\label{prop:Y+-OK}
$Y^+$ is compact in the product topology, and $(Y^+, G^+)$ is the product of the systems~$(Y_n, G_n)$.  In particular, $G^+$ is a homeomorphism.
\end{proposition}

Moreover, it follows from~(PDn1) and maximality of~$d$ among \JX-metrics  that
\begin{equation}\label{eqn:large-dist-Y}
diam(Y, d) = \eps.
\end{equation}

\subsection{Choosing $f,Y, G$}\label{subsec:choosef}

For every $n \in \bN$ we define a function~$f_n: X_n \rightarrow Y_n$ as follows:

For $(v, n, k) \in X_n$ we let $f_n(v, n, k) = (f^*(v), n, k)$, where
\begin{equation}\label{eqn:f*}
\forall i \in \bZ \ f^*(v)(i) = v(i) + 2.
\end{equation}

\begin{proposition}\label{prop:fn}
For all $n \geq 0$ the following properties hold:

\smallskip

\noindent
(i) $f_n$ is a bijection between $X_n$ and~$Y_n$.

\smallskip

\noindent
(ii) $f_n \circ F_n = G_n \circ f_n$

\smallskip

\noindent
(iii) For all~$x, x' \in X_n$ we have $D^n(x, x') \leq d^n(f_n(x), f_n(x'))$.
\end{proposition}

\noindent
\textbf{Proof:} (i) For all $(v,n,k) \neq (v', n, k') \in X_n$, we have $f_n(v,n,k) = (f^*(v),n,k)$ and $f_n(v',n, k') = (f^*(v'), n, k')$.
\begin{itemize}
\item If $k \neq k'$, then $f_n(v,n,k) = (f^*(v),n,k) \neq (f^*(v'), n, k') = f_n(v',n, k')$.
\item If $k = k'$, then $v \neq v'$.  There exists $i \in \bZ$ such that $v(i) \neq v'(i)$.  Then $f^*(v)(i) = v(i) + 2 \neq v'(i) + 2 = f^*(v')(i)$.  Hence,  $f_n(v,n,k) = (f^*(v),n,k) \neq (f^*(v'), n, k') = f_n(v',n, k')$.
\end{itemize}

Therefore, $f_n$ is one-to-one.

\smallskip

On the other hand, for each $(u, n, k) \in Y_n$, let $v \in {}^{\bZ}\{0,1\}$ be such that $v(i) = u(i)-2$ for all $i \in \bZ$.  Then $(v,n,k) \in X_n$ and $(u,n,k) = f_n(v,n,k)$.  Thus, $f_n$ is onto.

\smallskip

We can conclude that $f_n$ is a bijection between $X_n$ and $Y_n$.

\medskip

\noindent
(ii) For all $(v,n,k) \in X_n$,
\begin{equation*}
    \begin{split}
        F_n((v,n,k)) &= (\sigma(v),n,F_n(k)) = (\sigma(v), n, (k+1) \ mod\ T^+(n)),\\
        f_n(F_n((v,n,k))) &= (f^*(\sigma(v)), n, (k+1)\ mod\ T^+(n)),
    \end{split}
\end{equation*}
and
\begin{equation*}
    \begin{split}
        f_n(v,n,k) &= (f^*(v),n,k),\\
        G_n(f_n(v,n,k)) & = G_n((f^*(v),n,k)) = (\sigma(f^*(v)),n,(k+1)\ mod\ T^+(n)).
    \end{split}
\end{equation*}
For all $i \in \bZ$,
\begin{equation*}
    \begin{split}
        f^*(\sigma(v))(i) &= \sigma(v)(i) + 2 = v(i+1) + 2,\\
        \sigma(f^*(v))(i) &= f^*(v)(i+1) = v(i+1) +2.
    \end{split}
\end{equation*}
Therefore, $f_n \circ F_n = G_n \circ f_n$.

\medskip

\noindent
(iii)  Note that $D^n(x, x') = d^n(f_n(x), f_n(x'))$ \emph{unless} $x = (v, n, k')$ and
$x' = (v', n, k')$ with $k = k'$ and $\Delta(v, v') = 0$. \newline
In the latter case, we always have $\Delta(f^*(v), f^*(v')) = 0$ and hence  \newline
$d^n(f_n(x), f_n(x')) = \eps_n$, which is the diameter of~$(X_n, D^n)$. $\Box$

\bigskip

Now we define a function~$f^+: X \rightarrow Y^+$ as follows:
\begin{equation}\label{eqn:define-f}
f^+(x_0, x_1, \dots) = (f_0(x_0), f_1(x_1), \dots).
\end{equation}

Moreover, we define $Y = f^+(X^-)$, where~$X^-$ is the subspace of~$X$ constructed in Subsection~\ref{subsec:WX-}.
We let $f = f^+\restrict X^-$, and $G = G^+ \restrict Y$.

\begin{proposition}\label{prop:f-is-conjugacy}
Let $Y, f, G$ be defined as above.  Then

\smallskip

\noindent
(i) $f$ is a conjugacy between~$(X^-, F)$ and $(Y, G)$.

\smallskip

\noindent
(ii) $D(x, x') \leq d(f(x), f(x'))$ for all $x, x' \in X^-$.
\end{proposition}

\noindent
\textbf{Proof:} (i)  $f$ is one-to-one since each $f_n$ is one-to-one by Proposition \ref{prop:fn}(i). The same result implies that~$f^+$ is onto, and it
 follows from the definition of $Y$ that $f$ is onto.  Hence $f$ is a bijection from $X^-$ to $Y$.

 To show that $f$ is a homeomorphism, by compactness of~$X^-$ and~$Y$ we just need to show that $f$ is continuous.

 For all $0< \gamma < \frac{\eps_0}{3}$, there exists $N \in \bN$ such that $\sum_{n=N+1}^{\infty}\eps_n < \frac{\gamma}{2}$.  Fix this $N$ and choose $0 < \nu < \min\{\frac{\gamma}{2}, \frac{\delta_N}{3}\}$.  Then for all $x, x' \in X^-$ with $D(x, x') < \nu$,
\begin{equation}\label{eqn:dDf-small}
    \begin{split}
        d(f(x),f(x')) &= \sum_{n\in\bN} d^n(f_n(x_n),f_n(x_n'))\\
        &= \sum_{n=0}^N d^n(f_n(x_n),f_n(x_n')) + \sum_{n=N+1}^{\infty}d^n(f_n(x_n),f_n(x_n'))\\
        &= \sum_{n=0}^N D^n(x_n,x_n') + \sum_{n=N+1}^{\infty}d^n(f_n(x_n),f_n(x_n'))\\
        &\leq \sum_{n=0}^{\infty} D^n(x_n,x_n') + \sum_{n=N+1}^{\infty}\eps_n\\
        &< \nu + \frac{\gamma}{2}\\
        &< \frac{\gamma}{2} + \frac{\gamma}{2} = \gamma.
    \end{split}
\end{equation}

Notice that under our assumption  $D(x, x') < \nu <  \frac{\delta_N}{3}$, for ~$x = (v, n, k')$ and
$x' = (v', n, k')$ we must have $\Delta(v, v') > 0$.  Thus the observation that we made in the first two lines of the proof of Proposition~\ref{prop:fn}(iii)
gives the third equality in~\eqref{eqn:dDf-small}.

We have shown that $f$ is continuous, and it follows that $f$ is a homeomorphism.

Now it is left to show that $f \circ F = G \circ f$.  For all $x \in X^-$ and $n \in \bN$,
\begin{equation*}
    \begin{split}
        F(x)_n &= F_n(x_n),\\
        (f \circ F(x))_n &= f_n(F_n(x_n))\\
        &= G_n(f_n(x_n))\\
        &= G(f(x))_n\\
        &= (G \circ f (x))_n \in Y.
    \end{split}
\end{equation*}
Therefore, $f$ is a conjugacy between $(X^-, F)$ and $(Y, G)$.

\medskip

\noindent
(ii) For all $x, x' \in X^-$, we have $D^n(x_n, x_n') \leq d^n(f_n(x_n), f_n(x_n'))$ for all $n \in \bN$ by Proposition ~\ref{prop:fn}(iii).  Then
\begin{equation*}
    \begin{split}
        D(x,x') &= \sum_{n\in \bN}D^n(x_n, x_n')\\
        &\leq \sum_{n\in\bN} d^n(f_n(x_n), f_n(x_n'))\\
        &= d(f(x), f(x')).
    \end{split}
\end{equation*}
$\Box$

\bigskip

\subsection{The separation numbers $sep(Y, \delta, d_T)$}\label{subsec:sepY}

It remains to prove the following result.

\begin{lemma}\label{lem:sepY}
The system $(Y, G)$ satisfies~\eqref{eqn:sep=inY}.
\end{lemma}

\noindent
\textbf{Proof:} By the definition of $Y$ as $F^+(X^-)$, and since $X^- \subset W$ by~\eqref{eqn:X-subW}, for each $y \in Y$, there exists $u \in {}^{\bZ}\{2, 3\}$ such that $y_n = (u, n, k_n)$ for all $n \in \bN$, where each $k_n \in \{0, 1,  \dots  T^+(n)-1\}$.

Let us define
\begin{equation*}
    Y^{*} := \{u:  \exists y \in Y\ \mbox{such that}\ y_m = (u, m, k_m)\ \mbox{for all} \ m \in \bN\}.
\end{equation*}

Moreover, for each $n \in \bN$ and $u, u' \in Y^{*}$, let
\begin{equation*}
    d^{n*}(u, u') := \eps_n 3^{-\Delta(u,u')}.
\end{equation*}

Then $Y^* \subset {}^\bZ \{2, 3\}$ is a subshift, and each of the functions~$d^{n*}$ is a standard subshift metric.
Thus for the subshift system $(Y^*, \sigma)$ and each~$n \in \bN$, Proposition~\ref{prop:subshift=} implies:
\begin{equation}\label{Y*}
\forall \gamma > 0\ \ \liminf_{T \rightarrow 0} \frac{\ln{sep(Y^{*}, \gamma, d^{n*}_T)}}{T} = \limsup_{T \rightarrow 0} \frac{\ln{sep(Y^{*}, \gamma, d^{n*}_T)}}{T}.
\end{equation}

Since $diam(Y, d) = \eps$ by~\eqref{eqn:large-dist-Y}, for $\delta > \eps$ we always have $sep(Y, \delta , d_T) = 1$ so that \eqref{eqn:sep=inY} is trivial.

For $0 < \delta \leq \eps$, we choose first $M \in \bN$ and then $N \in \bN$ such that
\begin{equation}\label{eqn:choose-NMeps}
    \frac{\eps}{3^{M+1}} + \frac{1}{2}\eps_N < \delta \leq \frac{1}{3^M}\eps.
\end{equation}

Fix such a pair of $M$ and $N$.   We claim that for all $T > 0$,
\begin{equation}\label{eqn:sep-sandwich}
    sep\left(Y^{*}, \frac{\eps_N}{3^M}, d^{N*}_T\right) \leq sep(Y, \delta, d_T) \leq \left[\prod_{n=0}^N T^+(n)\right] sep\left(Y^{*}, \frac{\eps_N}{3^M}, d^{N*}_T\right).
\end{equation}

Suppose $A \subset Y^{*}$ is a $\left(T, \frac{\eps_N}{3^M}\right)$-separated subset of $Y^{*}$ with respect to the metric~$d^{N*}$.

For each $u \in A$, choose a $y \in Y$ with $y_n = (u, n, k_n)\ \mbox{for all}\ n \in \bN$, and let the collection of them be $B \subset Y$.
Then $|B| = |A|$.

Notice also that the metrics $d^{n*}$ differ just by scaling factors from each other.  In particular,
\begin{equation*}\label{eqn:dnt-scale}
\forall T >0 \, \forall u, u' \in Y^* \quad \left(d^{N*}_T (u,u') \geq \frac{\eps_N}{3^M} \ \ \Leftrightarrow \ \ \forall n \in \bN \ \  d^{n*}_T(u,u') \geq \frac{\eps_n}{3^M}\right).
\end{equation*}

Thus by the definition of~$d$, for all $y \neq y' \in B$,
\begin{equation*}
    d_T(y, y') \geq \sum_{n=0}^\infty d^{n*}_T(u, u') \geq \sum_{n=0}^{\infty} \frac{\eps_n}{3^M} = \frac{\eps}{3^M} \geq \delta.
\end{equation*}

It follows that~$B$ is~$(T, \delta)$-separated in~$(Y, G)$ with respect to~$d$, so that
\begin{equation*}
    sep\left(Y^{*}, \frac{\eps_N}{3^M}, d^{N*}_T\right) \leq sep(Y, \delta, d_T).
\end{equation*}

For the second inequality, assume towards a contradiction that there exists
$B \subset Y$ with $|B| > \left[\prod_{n=0}^N T^+(n)\right] sep\left(Y^{*}, \frac{\eps_N}{3^M}, d^{N*}_T\right)$ that is $(T, \delta)$-separated.  By the
 Pigeonhole Principle there exist $\bbk^* = (k^*_n)_{n = 0}^{N}$ and  $B_{\bbk^*} \subset B$ such that
\begin{itemize}
\item[(ks)]  for all $y = ((u,n,k_n))_{n=0}^{\infty} \in B_{\bbk^*}$, we have $(k_n)_{n=0}^N = (k^*_n)_{n = 0}^{N}$,

\smallskip

\item[(nosep)] $|B_{\bbk^*}| > sep\left(Y^{*}, \frac{\eps_N}{3^M}, d^{N*}_T\right)$,

\smallskip

\item[(sep)] $B_{\bbk^*}$ is $(T,\delta)$-separated in~$(Y, G)$ with respect to~$d$.
\end{itemize}

\begin{remark}\label{rem:spanonly-proof}
Not all sequences $(k^*_n)_{n=0}^N$ with $k^*_n \in T^+(n)$ for all~$n \leq N$ actually play a role in this Pigeonhole Principle-based argument. Since we started
from~$X^- \subset W$, our construction implies that only the ones with $k^*_n = k^*_N \mod T^+(n)$ are relevant. Thus the upper bound
in~\eqref{eqn:sep-sandwich} could be improved by replacing the term~$\left[\prod_{n=0}^N T^+(n)\right]$ with~$T^+(N)$.  However, this is not needed for our purposes, and it may be useful in follow-up work to have an argument, as the one given here, that requires only the weaker assumption that the space~$X^-$ that we started with has the property that for all $x \in X^-$ there exists a~$v \in {}^\bZ \{0,1\}$ such that $x_n = (v, n, k_n)$ for all~$n \in \bN$.
\end{remark}

Then by~(nosep) there exist $y \neq y' \in B_{\bbk^*}$ such that
\begin{equation}\label{eqn:dtu}
    d_T^{N*}(u, u') < \frac{\eps_N}{3^M}.
\end{equation}
By~(ks), clause~(dn3) will apply in the computation of $d^n(G_n^t(y_n), G^t(y'_n))$ for all $n \leq N$ and $t \in \bZ$, and by~\eqref{eqn:dtu} we will have:
\begin{equation*}
\begin{split}
&\forall n \leq N \ \forall 0 \leq t < T \ \     d^n(G_n^t(y_n), G_n^t(y'_n))  < \frac{\eps_n}{3^M},\\
&\forall n \leq N \  \     d^n_T(y_n, y'_n)  < \frac{\eps_n}{3^M},\\
&\forall n \leq N \  \     d^n_T(y_n, y'_n)  \leq \frac{\eps_n}{3^{M+1}},
\end{split}
\end{equation*}
where the last line follows from the fact that $d^n$ cannot take any values strictly between~$\frac{\eps_n}{3^{M+1}}$ and~$\frac{\eps_n}{3^{M}}$. Then
by~\eqref{eqn:choose-NMeps} and property (P$\delta$2):
\begin{equation*}
    \begin{split}
        d_T(y, y') &< \sum_{n=0}^{N} \frac{\eps_n}{3^M} + \sum_{n=N+1}^{\infty}\eps_n,\\
        d_T(y, y') &\leq \sum_{n=0}^{N} \frac{\eps_n}{3^{M+1}} + \sum_{n=N+1}^{\infty}\eps_n,\\
        d_T(y, y') &< \frac{\eps}{3^{M+1}} + \sum_{n=N+1}^{\infty}\eps_n\\
        &< \frac{\eps}{3^{M+1}} + \frac{1}{2}\eps_N \\
        &< \delta,
    \end{split}
\end{equation*}
which contradicts our assumption~(sep).  Hence,
\begin{equation*}
    sep(Y, \delta, d_T) \leq \left[\prod_{n=0}^N T^+(n)\right] sep\left(Y^{*}, \frac{\eps_N}{3^M}, d^{N*}_T\right).
\end{equation*}
We have shown that~\eqref{eqn:sep-sandwich} holds. Then,
\begin{equation*}
\begin{split}
    sep\left(Y^{*}, \frac{\eps_N}{3^M}, d^{N*}_T\right) &\leq sep(Y, \delta, d_T) \leq \left[\prod_{n=0}^N T^+(n)\right] sep\left(Y^{*}, \frac{\eps_N}{3^M}, d^{N*}_T\right)\\
    \frac{\ln{sep\left(Y^{*}, \frac{\eps_N}{3^M}, d^{N*}_T\right)}}{T} &\leq \frac{\ln{sep(Y, \delta, d_T)}}{T} \leq \frac{\ln{\left[\prod_{n=0}^N T^+(n)\right]}}{T} + \frac{\ln{sep\left(Y^{*}, \frac{\eps_N}{3^M}, d^{N*}_T\right)}}{T}.
\end{split}
\end{equation*}
Since
\begin{equation*}
     \lim_{T\rightarrow \infty} \frac{\ln{\left[\prod_{n=0}^N T^+(n)\right]}}{T} = 0,
\end{equation*}
it follows from~\eqref{Y*} that
\begin{equation*}
    \liminf_{T\rightarrow \infty} \frac{\ln{sep(Y, \delta, d_T)}}{T} = \limsup_{T\rightarrow \infty} \frac{\ln{sep(Y, \delta, d_T)}}{T}.  \qquad \Box
\end{equation*}

\section*{Acknowledgements}

We wish to thank T. Downarowicz, B. Hasselblatt, B. Weiss and T. Young for detailed and valuable feedback on the status of Question~\ref{q:motivating} and T. Young for suggesting the use of infinite products in our constructions.

\newpage
\section*{Appendix: Index of property abbreviations and other important notation}

\begin{tabular}{ |p{4cm}||p{3.5cm}|p{3.5cm}|  }
 \hline
 $T^+(n)$ & \multicolumn{2}{|c|}{Often used as shorthand for $\{0, 1, \dots , T^+(n)-1\}$} \\
 \hline
 $T(n), T^+(n)$   & Subsection~\ref{subsec:T(n)}    & page 19--20\\
 \hline
 $C(n), K(n)$     & Subsection~\ref{subsec:T(n)}    & page 20\\
 \hline
 (PCn) & Subsection~\ref{subsec:T(n)} & page 20\\
 \hline
 (PKn1)--(PKn5) & Subsection~\ref{subsec:T(n)} & page 20\\
 \hline
 (pKn3), (pKn4) & Subsection~\ref{subsec:T(n)} & page 20\\
 \hline
 (pcn) & Subsection~\ref{subsec:T(n)} & page 20\\
 \hline
 $I_j^n$  & Subsection~\ref{subsec:epsn-deltan} & page 22\\
 \hline
 $\eps_n, \eps$ &  Subsection~\ref{subsec:epsn-deltan} & page 22\\
 \hline
 $\delta_n$  & Subsection~\ref{subsec:epsn-deltan} & page 22\\
 \hline
 (P$\eps$) & Subsection~\ref{subsec:epsn-deltan} & page 22\\
 \hline
 (P$\delta_1$), (P$\delta_2$), (P$\delta_3$)  &  Subsection~\ref{subsec:epsn-deltan} & page 22\\
 \hline
 coloring & Subsection~\ref{subsec:choose-colors} & page 23\\
 \hline
 $[C(n)]$ & Subsection~\ref{subsec:choose-colors} & page 23 \\
 \hline
 $[S]^2$ & Subsection~\ref{subsec:choose-colors} & page 23 \\
 \hline
 (cC1), (cCi), (cC2), (cC), (cC3) & Subsection~\ref{subsec:choose-colors} & pages 23 \\
 \hline
 $X_n, X_n^k $ & Section~\ref{sec:Xn-construct} & page 26--27\\
 \hline $F_n$ & Section~\ref{sec:Xn-construct} & page 26\\
 \hline $D^n$ & Section~\ref{sec:Xn-construct} & page 27--28\\
 \hline $\Phi$ & Section~\ref{sec:Xn-construct} & page 27\\
 \hline $\#$ & Section~\ref{sec:Xn-construct} & page 27\\
 \hline $\Delta$ & Section~\ref{sec:Xn-construct} & page 27\\
 \hline (Dn1), (Dn2), (Dn3), (Dn31), (Dn32) & Section~\ref{sec:Xn-construct} & page 27--28\\
 \hline (Dn1d), (Dn1e), (Dn32c) & Section~\ref{sec:Xn-construct} & page 28\\
 \hline (PDn1), (PDn2), (Pnc), (PFn) & Section~\ref{sec:Xn-construct} & page 28\\
 \hline $X$ & Section~\ref{sec:X} & page 32\\
 \hline $D$ & Section~\ref{sec:X} & page 32\\
 \hline $F$ & Section~\ref{sec:X} & page 32\\
 \hline $\Ymn$ & Subsection~\ref{subsec:Xn-Yn-} & page 33\\
 \hline (PY1), (PY2), (PY2+) & Subsection~\ref{subsec:Xn-Yn-} & page 33--34\\
 \hline (PR1), (PR2) & Subsection~\ref{subsec:Xn-Yn-} & page 33\\
 \hline (P2-1) & Subsection~\ref{subsec:lower-sep-T+n} & page 46\\
 \hline  $W$ & Subsection~\ref{subsec:WX-} & page~37\\
 \hline  $W^n$ & Subsection~\ref{subsec:WX-} & page~37\\
 \hline  $X^-$ & Subsection~\ref{subsec:WX-} & page~37\\
 \hline $y_\varphi$ & Subsection~\ref{subsec:WX-} & page~37\\
 \hline $x^\varphi$ & Subsection~\ref{subsec:WX-} & page~37\\
 \hline (RD) & Subsection~\ref{subsec:spanonly-general} & page~60\\
 \hline (Rd) & Subsection~\ref{subsec:spanonly-general} & page~60\\
 \hline (R2) & Subsection~\ref{subsec:spanonly-general} & page~60\\

 \hline
\end{tabular}

\end{document}